\documentclass[11pt]{amsart}
\usepackage[letterpaper, total={6.5in,9in}, includefoot, centering]{geometry}
\usepackage{amsfonts}
\usepackage{amsmath,amsthm,amssymb,amscd}
\usepackage[foot]{amsaddr}

\usepackage[mathscr]{eucal}
\usepackage{latexsym}
\usepackage[new]{old-arrows}
\usepackage{graphics}
\usepackage{verbatim}
\usepackage{upgreek}
\usepackage{overpic}
\usepackage[all,cmtip]{xy} 
\usepackage{enumitem}
\usepackage{tikz}
\usepackage{tikz-cd}
\usepackage[pagebackref]{hyperref}
\renewcommand*\backref[1]{\ifx#1\relax \else (Page #1) \fi}
\usepackage[font=scriptsize]{caption}
\usepackage{upgreek}
\makeatletter
\newsavebox{\@brx}
\newcommand{\llangle}[1][]{\savebox{\@brx}{\(\m@th{#1\langle}\)}%
  \mathopen{\copy\@brx\kern-0.5\wd\@brx\usebox{\@brx}}}
\newcommand{\rrangle}[1][]{\savebox{\@brx}{\(\m@th{#1\rangle}\)}%
  \mathclose{\copy\@brx\kern-0.5\wd\@brx\usebox{\@brx}}}
\makeatother
\graphicspath{ {./images/} }

\usepackage[authoryear,sort&compress]{natbib}
\usepackage{adjustbox}
\usepackage{float}
\bibpunct{[}{]}{;}{n}{,}{,}

\makeatletter \@addtoreset{equation}{section}

\DeclareMathOperator{\ext}{ext}

\DeclareMathOperator*{\argmin}{arg\,min}

\renewcommand{\theequation}{\arabic{section}.\arabic{equation}}
\makeatother

\newtheorem{remark}{Remark}[section]
\newtheorem{example}{Example}[section]
\newtheorem{prop}{Proposition}[section]

\begin{document}
\title{Geometric Methods for Adjoint Systems}
\author{Brian Kha Tran and Melvin Leok}
\address{Department of Mathematics, University of California, San Diego, 9500 Gilman Drive, La Jolla, CA 92093-0112, USA.}
\email{b3tran@ucsd.edu, mleok@ucsd.edu}
\allowdisplaybreaks

\begin{abstract}
Adjoint systems are widely used to inform control, optimization, and design in systems described by ordinary differential equations or differential-algebraic equations. In this paper, we explore the geometric properties and develop methods for such adjoint systems. In particular, we utilize symplectic and presymplectic geometry to investigate the properties of adjoint systems associated with ordinary differential equations and differential-algebraic equations, respectively. We show that the adjoint variational quadratic conservation laws, which are key to adjoint sensitivity analysis, arise from (pre)symplecticity of such adjoint systems. We discuss various additional geometric properties of adjoint systems, such as symmetries and variational characterizations. For adjoint systems associated with a differential-algebraic equation, we relate the index of the differential-algebraic equation to the presymplectic constraint algorithm of \citet{GoNeHi1978}. As an application of this geometric framework, we discuss how the adjoint variational quadratic conservation laws can be used to compute sensitivities of terminal or running cost functions. Furthermore, we develop structure-preserving numerical methods for such systems using Galerkin Hamiltonian variational integrators (\citet{LeZh2011}) which admit discrete analogues of these quadratic conservation laws. We additionally show that such methods are natural, in the sense that reduction, forming the adjoint system, and discretization all commute, for suitable choices of these processes. We utilize this naturality to derive a variational error analysis result for the presymplectic variational integrator that we use to discretize the adjoint DAE system. Finally, we discuss the application of adjoint systems in the context of optimal control problems, where we prove a similar naturality result.
\end{abstract}

\maketitle

\tableofcontents

\section{Introduction}

\subsection{Applications of the Adjoint Equations}

The solution of many nonlinear problems involves successive linearization, and as such variational equations and their adjoints play a critical role in a variety of applications. Adjoint equations are of particular interest when the parameter space is significantly higher dimension than that of the output or objective. In particular, the simulation of adjoint equations arise in sensitivity analysis~\cite{Ca1981, CaLiPeSe2003}, adaptive mesh refinement~\cite{LiPe2003}, uncertainty quantification~\cite{WaDuAlIa2012}, automatic differentiation~\cite{Gr2003}, superconvergent functional recovery~\cite{PiGi2000}, optimal control~\cite{Ro2005}, optimal design~\cite{GiPi2000}, optimal estimation~\cite{NgGeBe2016}, and deep learning viewed as an optimal control problem~\cite{DeCeEhOwSc2019}.

The study of geometric aspects of adjoint systems arose from the observation that the combination of any system of differential equations and its adjoint equations are described by a formal Lagrangian~\cite{Ib2006, Ib2007}. This naturally leads to the question of when the formation of adjoints and discretization commutes~\cite{SoTz1997}, and prior work on this include the Ross--Fahroo lemma~\cite{RoFa2001}, and the observation by \citet{Sa2016} that the adjoints and discretization commute if and only if the discretization is symplectic.

\subsection{Symplectic and Presymplectic Geometry}\label{SymplecticGeometrySection}
Throughout the paper, we will assume that all manifolds and maps are smooth, unless otherwise stated. Let $(P,\Omega)$ be a (finite-dimensional) symplectic manifold, i.e., $\Omega$ is a closed nondegenerate two-form on $P$. Given a Hamiltonian $H: P \rightarrow \mathbb{R}$, the Hamiltonian system is defined by
$$ i_{X_H} \Omega = dH, $$
where the vector field $X_H$ is a section of the tangent bundle to $P$. By nondegeneracy, the vector field $X_H$ exists and is uniquely determined. For an open interval $I \subset \mathbb{R}$, we say that a curve $z: I \rightarrow P$ is a solution of Hamilton's equations if $z$ is an integral curve of $X_H$, i.e., $\dot{z}(t) = X_H(z(t))$ for all $t \in I$.

A particularly important example for our purposes is when the symplectic manifold is the cotangent bundle of a manifold, $P = T^*M$, equipped with the canonical symplectic form $\Omega = dq \wedge dp$ in natural coordinates $(q,p)$ on $T^*M$. A Hamiltonian system has the coordinate expression
\begin{align*}
\dot{q} &= \frac{\partial H(q,p)}{\partial p}, \\
\dot{p} &= - \frac{\partial H(q,p)}{\partial q}.
\end{align*}
By Darboux's theorem, any symplectic manifold is locally symplectomorphic to a cotangent bundle equipped with its canonical symplectic form. As such, any Hamiltonian system can be locally expressed in the above form (even when $P$ is not a cotangent bundle), using Darboux coordinates.

We now consider the generalization of Hamiltonian systems where we relax the condition that $\Omega$ is nondegenerate, i.e., presymplectic geometry. Let $(P,\Omega)$ be a presymplectic manifold, i.e., $\Omega$ is a closed two-form on $P$ with constant rank. As before, given a Hamiltonian $H: P \rightarrow \mathbb{R}$, we define the associated Hamiltonian system as
$$ i_{X_H}\Omega = dH. $$
Note that since $\Omega$ is now degenerate, $X_H$ is not guaranteed to exist and if it does, it need not be unique and in general is only partially defined on a submanifold of $P$. Again, we say a curve on $P$ is a solution to Hamilton's equations if it is an integral curve of $X_H$. Using Darboux coordinates $(q,p,r)$ adapted to $(P,\Omega)$, where $\Omega = dq \wedge dp$ and $\ker(\Omega) = \text{span}\{\partial/\partial r\}$, the local expression for Hamilton's equations is given by
\begin{align*}
\dot{q} &= \frac{\partial H(q,p,r)}{\partial p}, \\
\dot{p} &= -\frac{\partial H(q,p,r)}{\partial q}, \\
0 &= \frac{\partial H(q,p,r)}{\partial r}.
\end{align*}
The third equation above is interpreted as a constraint equation which any solution curve must satisfy. We will assume that the constraint defines a submanifold of $P$. It is clear that in order for a solution vector field $X_H$ to exist, it must be restricted to lie on this submanifold. However, in order for its flow to remain on the submanifold, it must be tangent to this submanifold, which further restricts where $X$ can be defined. Alternating restriction in order to satisfy these two constraints yields the presymplectic constraint algorithm  of \citet{GoNeHi1978}. The presymplectic constraint algorithm begins with the observation that for any $X$ satisfying the above system, so does $X+Z$, where $Z \in \text{ker}(\Omega)$. In order to obtain such a vector field $X$, one considers the subset $P_1$ of $P$ such that $Z_p(H) = 0$ for any $Z \in \text{ker}(\Omega), p \in P_1$. We will assume that the set $P_1$ is a submanifold of $P$. We refer to $P_1$ as the primary constraint manifold. In order for the flow of the resulting Hamiltonian vector field $X$ to remain on $P_1$, one further requires that $X$ is tangent to $P_1$. The set of points satisfying this property defines a subsequent secondary constraint submanifold $P_2$. Iterating this process, one obtains a sequence of submanifolds
$$ \dots \rightarrow P_k \rightarrow \dots \rightarrow P_1 \rightarrow P_0 \equiv P, $$
defined by
\begin{equation}\label{Presymplectic Constraint Algorithm Sequence}
P_{k+1} = \{ p \in P_k : Z_p(H_k) = 0 \text{ for all } Z \in \text{ker}(\Omega_k)\}, 
\end{equation}
where
\begin{align*}
\Omega_{k+1} &= \Omega_k|_{P_{k+1}}, \\
H_{k+1} &= H_k|_{P_{k+1}}.
\end{align*}
If there exists a nontrivial fixed point in this sequence, i.e., a submanifold $P_k$ of $P$ such that $P_{k} = P_{k+1}$, we refer to $P_{k}$ as the final constraint manifold. If such a fixed point exists, we denote by $\nu_P$ the minimum integer such that $P_{\nu_P} = P_{\nu_P+1}$, i.e., $\nu_P$ is the number of steps necessary for the presymplectic constraint algorithm to terminate. If such a final constraint manifold $P_{\nu_P}$ exists, there always exists a solution vector field $X$ defined on and tangent to $P_{\nu_P}$ such that $i_X \Omega_{\nu_P} = dH_{\nu_P}$ and $X$ is unique up to the kernel of $\Omega_{\nu_P}$. Furthermore, such a final constraint manifold is maximal in the sense that if there exists a submanifold $N$ of $P$ which admits a vector field $X$ defined on and tangent to $N$ such that $i_X\Omega|_N = dH|_N$, then $N \subset P_{\nu_P}$ (\citet{GoNe1979}).

\subsection{Main Contributions}
In this paper, we explore the geometric properties of adjoint systems associated with ordinary differential equations (ODEs) and differential-algebraic equations (DAEs). For a discussion of adjoint systems associated with ODEs and DAEs, see \citet{Sa2016} and \citet{CaLiPeSe2003}, respectively. In particular, we utilize the machinery of symplectic and presymplectic geometry as a basis for understanding such systems.

In Section \ref{AdjointVectorSpaceSection}, we review the notion of adjoint equations associated with ODEs over vector spaces. We show that the quadratic conservation law, which is the key to adjoint sensitivity analysis, arises from the symplecticity of the flow of the adjoint system. In Section \ref{AdjointManifoldSection}, we investigate the symplectic geometry of adjoint systems associated with ODEs on manifolds. We additionally discuss augmented adjoint systems, which are useful in the adjoint sensitivity of running cost functions. In Section \ref{AdjointDAESection}, we investigate the presymplectic geometry of adjoint systems associated with DAEs on manifolds. We investigate the relation between the index of the base DAE and the index of the associated adjoint system, using the notions of DAE reduction and the presymplectic constraint algorithm. We additionally consider augmented systems for such adjoint DAE systems. For the various adjoint systems that we consider, we derive various quadratic conservation laws which are useful in adjoint sensitivity analysis of terminal and running cost functions. We additionally discuss symmetry properties and present variational characterizations of such systems that provide a useful perspective for constructing geometric numerical methods for these systems.

In Section \ref{ApplicationsSection}, we discuss applications of the various adjoint systems to adjoint sensitivity and optimal control. In Section \ref{AdjointSensitivitySection}, we show how the quadratic conservation laws developed in Section \ref{AdjointSystems Main Section} can be used for adjoint sensitivity analysis of running and terminal cost functions, subject to ODE or DAE constraints. In Section \ref{DiscreteAdjointSystemsSection}, we construct structure-preserving discretizations of adjoint systems using the Galerkin Hamiltonian variational integrator construction  of \citet{LeZh2011}. For adjoint DAE systems, we introduce a presymplectic analogue of the Galerkin Hamiltonian variational integrator construction. We show that such discretizations admit discrete analogues of the aforementioned quadratic conservation laws and hence are suitable for the numerical computation of adjoint sensitivities. Furthermore, we show that such discretizations are natural when applied to DAE systems, in the sense that reduction, forming the adjoint system, and discretization all commute (for particular choices of these processes). As an application of this naturality, we derive a variational error analysis result for the resulting presymplectic variational integrator for adjoint DAE systems. Finally, in Section \ref{OCPSection}, we discuss adjoint systems in the context of optimal control problems, where we prove a similar naturality result, in that suitable choices of reduction, extremization, and discretization commute.

By developing a geometric theory for adjoint systems, the application areas that utilize such adjoint systems can benefit from the existing work on geometric and structure-preserving methods.

\subsection{Main Results}
In this paper, we prove that, starting with an index 1 DAE, appopriate choices of reduction, discretization, and forming the adjoint system commute. That is, the following diagram commutes.

\[\small\begin{tikzcd}[column sep=-4ex,row sep=10ex]
	\txt{Index 1 DAE} && \txt{ODE} \\
	& \txt{Discrete DAE} && \txt{Discrete ODE} \\
	\txt{Presymplectic Adjoint\\ DAE System} && \txt{Symplectic Adjoint\\ ODE System} \\
	& \txt{Presymplectic Galerkin\\ Hamiltonian Variational \\ Integrator} && \txt{Symplectic Galerkin \\ Hamiltonian Variational \\ Integrator}
	\arrow["{\text{Reduce}}", from=1-1, to=1-3]
	\arrow["{\text{Reduce}}"{pos=0.25}, from=3-1, to=3-3]
	\arrow["{\text{Adjoint}}"'{pos=0.6}, from=1-1, to=3-1]
	\arrow["{\text{Adjoint}}"{pos=0.6}, from=1-3, to=3-3]
	\arrow["{\text{Reduce}}"{pos=0.3}, from=2-2, to=2-4,crossing over]
	\arrow["{\text{Reduce}}"', from=4-2, to=4-4]
	\arrow["{\text{Adjoint}}"{pos=0.6}, from=2-2, to=4-2,crossing over]
	\arrow["{\text{Adjoint}}"{pos=0.6}, from=2-4, to=4-4]
	\arrow["{\text{Discretize}}"', from=1-1, to=2-2]
	\arrow["{\text{Discretize}}"{pos=0.6}, from=1-3, to=2-4]
	\arrow["{\text{Discretize}}"', from=3-1, to=4-2]
	\arrow["{\text{Discretize}}"', from=3-3, to=4-4]
\end{tikzcd}\]

In order to prove this result, we develop along the way the definitions of the various vertices and arrows in the above diagram. Roughly speaking, the four ``Adjoint" arrows are defined by forming the appropriate continuous or discrete action and enforcing the variational principle; the four ``Reduce" arrows are defined by solving the algebraic variables in terms of the kinematic variables through the continuous or discrete constraint equations; the two ``Discretize" arrows on the top face are given by a Runge--Kutta method, while the two ``Discretize" arrows on the bottom face are given by the associated symplectic partitioned Runge--Kutta method. The above commutative diagram can be understood as an extension of the result of \citet{Sa2016} (that discretization and forming the adjoint of an ODE commute when the discretization is a symplectic Runge--Kutta method) by adding the reduction operation. In order to appropriately define this reduction operation, we will show that the presymplectic adjoint DAE system has index 1 if the base DAE has index 1, so that the reduction of the presymplectic adjoint DAE system results in a symplectic adjoint ODE system; the tool for this will be the presymplectic constraint algorithm. 

In the process of defining the ingredients in the above diagram, we will additionally prove various properties of adjoint systems associated with ODEs and DAEs. The key properties that we will prove for such adjoint systems are the adjoint variational quadratic conservation laws, Propositions \ref{ManifoldQuadraticInvariantProp}, \ref{AugmentedQuadraticInvariantProp}, \ref{Quadratic Invariant Adjoint DAE Prop}, \ref{Quadratic Invariant Augmented DAE Adjoint Prop}. As we will show, these conservation laws can be used to compute adjoint sensitivities of running and terminal cost functions under the flow of an ODE or DAE. In order to prove these conservation laws, we will need to define the variational equations associated with an adjoint system. We will define them as the linearization of the base ODE or DAE; for the DAE case, we will show that the variational equations have the same index as the base DAE so that they have the same (local) solvability. 

\section{Adjoint Systems}\label{AdjointSystems Main Section}

\subsection{Adjoint Equations on Vector Spaces}\label{AdjointVectorSpaceSection}
In this section, we review the notion of adjoint equations on vector spaces and their properties, as preparation for adjoint systems on manifolds.

Let $Q$ be a finite-dimensional vector space and consider the ordinary differential equation on $Q$ given by
\begin{equation}\label{ODEvectorspace}
\dot{q} = f(q),
\end{equation}
where $f: Q \rightarrow Q$ is a differentiable vector field on $Q$. Let $Df(q)$ denote the linearization of $f$ at $q \in Q$, $Df(q) \in L(Q,Q)$. Denoting its adjoint by $[Df(q)]^* \in L(Q^*,Q^*)$, the adjoint equation associated with \eqref{ODEvectorspace} is given by
\begin{equation}\label{AdjointODEvectorspace}
\dot{p} = -[Df(q)]^* p,
\end{equation}
where $p$ is a curve on $Q^*$. 

Let $q^A$ be coordinates for $Q$ and let $p_A$ be the associated dual coordinates for $Q^*$, so that the duality pairing is given by $\langle p,q\rangle = p_Aq^A$. The linearization of $f$ at $q$ is given in coordinates by
$$ (Df(q))^A_B = \frac{\partial f^A(q)}{\partial q^B}, $$
where its action on $v \in Q$ in coordinates is
$$ (Df(q) v)^A = \frac{\partial f^A(q)}{\partial q^B} v^B. $$
Its adjoint then acts on $p \in Q^*$ by
$$ ([Df(q)]^* p)_A = \frac{\partial f^B(q)}{\partial q^A} p_B. $$
Thus, the ODE and its adjoint can be expressed in coordinates as
\begin{align*}
\dot{q}^A &= f^A(q), \\
\dot{p}_A &= - \frac{\partial f^B(q)}{\partial q^A} p_B.
\end{align*}

Next, we recall that the combined system \eqref{ODEvectorspace}-\eqref{AdjointODEvectorspace}, which we refer to as the adjoint system, arises from a variational principle. Letting $\langle \cdot,\cdot\rangle$ denote the duality pairing between $Q^*$ and $Q$, we define the Hamiltonian 
\begin{align*}
H: Q \times Q^* &\rightarrow \mathbb{R}, \\
(q,p) &\mapsto H(q,p) \equiv \langle p, f(q)\rangle.
\end{align*}
The associated action, defined on the space of curves on $Q \times Q^*$ covering some interval $(t_0,t_1)$, is given by
$$ S[q,p] = \int_{t_0}^{t_1} \left( \langle p,\dot{q}\rangle - H(q,p) \right) dt = \int_{t_0}^{t_1} \left(\langle p,\dot{q}\rangle - \langle p,f(q)\rangle \right) dt.  $$
\begin{prop}\label{VariationalPrincipleVectorSpaceCase}
The variational principle $\delta S = 0$, subject to variations $(\delta q,\delta p)$ which fix the endpoints $\delta q(t_0) = 0$, $\delta q(t_1) = 0$, yields the adjoint system \eqref{ODEvectorspace}-\eqref{AdjointODEvectorspace}.
\begin{proof}
Compute the variation of $S$ with respect to a compactly supported variation $(\delta q, \delta p)$,
\begin{align*}
\delta S[q,p] \cdot (\delta q, \delta p) &= \frac{d}{d\epsilon}\Big|_{\epsilon = 0} S[q + \epsilon \delta q, p + \epsilon \delta p] = \int_{t_0}^{t_1} \frac{d}{d\epsilon}\Big|_{\epsilon = 0} \langle p + \epsilon \delta p, \dot{q} + \epsilon \dot{\delta q} - f(q + \epsilon \delta q) \rangle dt \\
&= \int_{t_0}^{t_1} \Big( \langle \delta p, \dot{q} - f(q)\rangle + \langle p, \dot{\delta q} - Df(q) \delta q \rangle \Big) dt \\
&= \int_{t_0}^{t_1} \Big( \langle \delta p, \dot{q} - f(q)\rangle + \langle -\dot{p} - [Df(q)]^* p, \delta q\rangle \Big) dt.
\end{align*}
The fundamental lemma of the calculus of variations then yields \eqref{ODEvectorspace}-\eqref{AdjointODEvectorspace}.
\end{proof}
\end{prop}

\begin{remark}
Note that an analogous statement of Proposition \ref{VariationalPrincipleVectorSpaceCase} can also be stated using the Type II variational principle, where one instead considers the generating function
$$ H_+(q_0,p_1) = \ext \left[ p(t_1) q(t_1) - \int_{t_0}^{t_1} \langle p,\dot{q} - H(q,p)\rangle dt \right], $$
and one extremizes over $C^2$ curves from $[t_0,t_1]$ to $T^*Q$ such that $q(t_0) = q_0, p(t_1) = p_1$. The Type II variational principle again gives the above adjoint system, but with differing boundary conditions. These boundary conditions are typical in adjoint sensitivity analysis, where one fixes the initial position and the final momenta.
\end{remark}

The variational principle utilized above is formulated so that the stationarity condition $\delta S = 0$ is equivalent to Hamilton's equations, where we view $Q \times Q^* \cong T^*Q$ with the canonical symplectic form on the cotangent bundle $\Omega = dq \wedge dp$ and with the corresponding Hamiltonian $H: T^*Q \rightarrow \mathbb{R}$ given as above. It then follows that the flow of the adjoint system is symplectic. 

The symplecticity of the adjoint system is a key feature of the system. In fact, the symplecticity of the adjoint system implies that a certain quadratic invariant is preserved along the flow of the system. This quadratic invariant is the key ingredient to the use of adjoint equations for sensitivity analysis. To state the quadratic invariant, consider the variational equation associated with equation \eqref{ODEvectorspace},
\begin{equation}\label{VariationalODEvectorspace}
\frac{d}{dt}\delta q = Df(q)\delta q,
\end{equation}
which corresponds to the linearization of \eqref{ODEvectorspace} at $q \in Q$. For solution curves $p$ and $\delta q$ to \eqref{AdjointODEvectorspace} and \eqref{VariationalODEvectorspace}, respectively, over the same curve $q$, one has that the quantity $\langle p, \delta q\rangle$ is preserved along the flow of the system, since
\begin{align*}
\frac{d}{dt} \langle p,\delta q\rangle &= \langle \dot{p},\delta q\rangle + \langle p, \frac{d}{dt}\delta q\rangle = \langle -[Df(q)]^*p,\delta q\rangle + \langle p, Df(q)\delta q\rangle  \\
&= - \langle p, Df(q)\delta q\rangle + \langle p, Df(q)\delta q\rangle = 0.
\end{align*}
To see that symplecticity implies the preservation of this quantity, recall that symplecticity is the statement that, along a solution curve of the adjoint system \eqref{ODEvectorspace}-\eqref{AdjointODEvectorspace}, one has
$$ \frac{d}{dt}\Omega (V,W) = 0,$$
where $V$ and $W$ are first variations to the adjoint system (i.e., that the flow of $V$ and $W$ on solutions are again solutions). Infinitesimally, first variations $V$ and $W$ correspond to solutions of the linearization of the adjoint system \eqref{ODEvectorspace}-\eqref{AdjointODEvectorspace}. At a solution $(q,p)$ to the adjoint system, the linearization of the system is given by
\begin{align*}
\frac{d}{dt} \delta q &= Df(q)\delta q, \\
\frac{d}{dt} \delta p &= -[Df(q)]^* \delta p.
\end{align*}
Note that the first equation is just the variational equation \eqref{VariationalODEvectorspace} while the second equation is the adjoint equation \eqref{AdjointODEvectorspace}, with $p$ replaced by $\delta p$, since the adjoint equation is linear in $p$. The first variation vector field $V$ corresponding to a solution $(\delta q, \delta p)$ of this linearized system is 
$$ V = \delta q \frac{\partial}{\partial q} + \delta p \frac{\partial}{\partial p}. $$
Now, we make two choices for the first variations $V$ and $W$. For $W$, we take the solution $\delta q=0$, $\delta p = p$ of the linearized system, which gives $W = p\, \partial/\partial p$. For $V$, we take the solution $\delta q = \delta q$, $\delta p = 0$ of the linearized system, which gives $V = \delta q\, \partial/\partial q$. Inserting these into $\Omega$ gives
$$ \Omega(V,W) = p \frac{\partial}{\partial p} \lrcorner \left( \delta q \frac{\partial}{\partial q} \lrcorner ( dq \wedge dp ) \right) = \langle p,\delta q\rangle. $$
Thus, symplecticity $\frac{d}{dt}\Omega(V,W) = 0$ with this particular choice of first variations $V,W$ gives the preservation of the quadratic invariant $\langle p,\delta q\rangle$.

\subsection{Adjoint Systems on Manifolds}\label{AdjointManifoldSection}
We now extend the notion of the adjoint system to the case where the configuration space of the base ODE is a manifold. We will provide a symplectic characterization of these adjoint systems, prove the associated adjoint variational quadratic conservation laws, and additionally discuss symmetries and variational principles associated with these systems.

Let $M$ be a manifold and consider the ODE on $M$ given by
\begin{equation}\label{ODEmanifold}
\dot{q} = f(q),
\end{equation} 
where $f$ is a vector field on $M$. Letting $\pi: TM \rightarrow M$ denote the tangent bundle projection, we recall that a vector field $f$ is a map $f: M \rightarrow TM$ which satisfies $\pi \circ f = \mathbf{1}_M$, i.e., $f$ is a section of the tangent bundle.

Analogous to the adjoint system on vector spaces, we will define the adjoint system on a manifold as an ODE on the cotangent bundle $T^*M$ which covers \eqref{ODEmanifold}, such that the time evolution of the momenta in the fibers of $T^*M$ are given by an adjoint linearization of $f$.

To do this, in analogy with the vector space case, consider the Hamiltonian $H: T^*M \rightarrow \mathbb{R}$ given by $H(q,p) = \langle p, f(q) \rangle_q$ where $\langle\cdot,\cdot\rangle_q$ is the duality pairing of $T^*_qM$ with $T_qM$. When there is no possibility for confusion of the base point, we simply denote this duality pairing as $\langle\cdot,\cdot\rangle$. Recall that the cotangent bundle $T^*M$ possesses a canonical symplectic form $\Omega = -d\Theta$ where $\Theta$ is the tautological one-form on $T^*M$. With coordinates $(q,p) = (q^A, p_A)$ on $T^*M$, this symplectic form has the coordinate expression $\Omega = dq\wedge dp \equiv dq^A \wedge dp_A$. 

We define the adjoint system as the ODE on $T^*M$ given by Hamilton's equations, with the above choice of Hamiltonian $H$ and the canonical symplectic form. Thus, the adjoint system is given by the equation
$$ i_{X_H}\Omega = dH, $$
whose solution curves on $T^*M$ are the integral curves of the Hamiltonian vector field $X_H$. As is well-known, for the particular choice of Hamiltonian $H(q,p) = \langle p, f(q)\rangle$, the Hamiltonian vector field $X_H$ is given by the cotangent lift $\widehat{f}$ of $f$, which is a vector field on $T^*M$ that covers $f$ (for a discussion of the geometry of the cotangent bundle and lifts, see \citet{YaIs1973}; for a discussion of cotangent lifts in the context of optimal control, see \citet{BuLe2014}). With coordinates $z = (q,p)$ on $T^*M$, the adjoint system is the ODE on $T^*M$ given by
\begin{equation}\label{AdjointODEmanifold}
\dot{z} = \widehat{f}(z). 
\end{equation}
To be more explicit, recall that the cotangent lift of $f$ is constructed as follows. Let $\Phi_{\epsilon}: M \rightarrow M$ denote the one-parameter family of diffeomorphisms generated by $f$. Then, we consider the cotangent lifted diffeomorphisms given by $(\Phi_{-\epsilon})^*: T^*M \rightarrow T^*M$. This covers $\Phi_{\epsilon}$ in the sense that $\pi_{T^*M} \circ (\Phi_{-\epsilon})^* = \Phi_{\epsilon} \circ \pi_{T^*M} $ where $\pi_{T^*M}: T^*M \rightarrow M$ is the cotangent projection. The cotangent lift $\widehat{f}$ is then defined to be the infinitesimal generator of the cotangent lifted flow,
$$ \widehat{f}(z) = \frac{d}{d\epsilon}\Big|_{\epsilon=0} (\Phi_{-\epsilon})^* (z). $$
We can directly verify that $\widehat{f}$ is the Hamiltonian vector field for $H$, which follows from
$$ i_{\widehat{f}}\Omega = -i_{\widehat{f}}d\Theta = -\mathcal{L}_{\widehat{f}}\Theta + d( i_{\widehat{f}}\Theta ) = d( i_{\widehat{f}}\Theta) = dH, $$
where $\mathcal{L}_{\hat{f}}\Theta = 0$ follows from the fact that cotangent lifted flows preserve the tautological one-form and $H = i_{\widehat{f}}\Theta$ follows from a direct computation (where $i_{\widehat{f}}\Theta$ is interpreted as a function on the cotangent bundle which maps $(q,p)$ to $\langle \Theta(q,p), \widehat{f}(q,p)\rangle$)

The adjoint system \eqref{AdjointODEmanifold} covers \eqref{ODEmanifold} in the following sense.
\begin{prop}\label{LiftIntegralCurveProp}
Integral curves to the adjoint system \eqref{AdjointODEmanifold} lift integral curves to the system \eqref{ODEmanifold}.
\begin{proof}
Let $z = (q,p)$ be coordinates on $T^*M$. Let $(\dot{q},\dot{p}) \in T_{(q,p)}T^*M$. Then, $T\pi_{T^*M} (\dot{q},\dot{p}) = \dot{q}$ where $T\pi_{T^*M}$ is the pushforward of the cotangent projection. Furthermore,
\begin{align*}
T\pi_{T^*M} \hat{f}(q,p) &= T\pi_{T^*M} \frac{d}{d\epsilon}\Big|_{\epsilon = 0} (\Phi_{-\epsilon})^*(q,p) = \frac{d}{d\epsilon}\Big|_{\epsilon = 0} (\pi_{T^*M} \circ (\Phi_{-\epsilon})^*)(q,p) \\
&= \frac{d}{d\epsilon}\Big|_{\epsilon = 0} (\Phi_{\epsilon} \circ \pi_{T^*M})(q,p) = \frac{d}{d\epsilon}\Big|_{\epsilon=0} \Phi_{\epsilon}(q) = f(q). 
\end{align*}
Thus, the pushforward of the cotangent projection applied to \eqref{AdjointODEmanifold} gives \eqref{ODEmanifold}. It then follows that integral curves of \eqref{AdjointODEmanifold} lift integral curves of \eqref{ODEmanifold}.
\end{proof}
\end{prop}

\begin{remark}
This can also be seen explicitly in coordinates. Recalling that $i_{\widehat{f}}\Omega = dH$, one has
$$ dH = d(p_A f^A(q)) = f^A(q) dp_A + p_B \frac{\partial f^B(q)}{\partial q^A} dq^A,$$
and, on the other hand, denoting $\widehat{f}(q,p) = X^A(q,p) \partial/\partial{q^A} + Y_A(q,p) \partial/\partial{p_A}$,
$$ i_{\widehat{f}}\Omega = (X^A(q,p) \partial_{q^A} + Y_A(q,p) \partial_{p_A}) \lrcorner\, (dq^B \wedge dp_B) = X^A(q,p)dp_A - Y_A(q,p) dq^A. $$
Equating these two gives the coordinate expression for the cotangent lift $\widehat{f}$,
$$ \widehat{f}(q,p) = f^A(q) \frac{\partial}{\partial q^A} - p_B \frac{\partial f^B(q)}{\partial q^A} \frac{\partial}{\partial p_A}. $$
Thus, the system $\dot{z} = \widehat{f}(z)$ can be expressed in coordinates as
\begin{subequations}
\begin{align}
\dot{q}^A &= f^A(q), \label{AdjointODEa} \\
\dot{p}_A &=- p_B \frac{\partial f^B(q)}{\partial q^A}, \label{AdjointODEb}
\end{align}
\end{subequations}
which clearly covers the original ODE $\dot{q}^A = f^A(q)$. Also, note that this coordinate expression for the adjoint system recovers the coordinate expression for the adjoint system in the vector space case.
\end{remark} 

Analogous to the vector space case, the adjoint system possesses a quadratic invariant associated with the variational equations of \eqref{ODEmanifold}. The variational equation is given by considering the tangent lifted vector field on $TM$, $\widetilde{f}: TM \rightarrow TTM$, which is defined in terms of the flow $\Phi_{\epsilon}$ generated by $f$ by
$$ \widetilde{f}(q,\delta q) = \frac{d}{d\epsilon}\Big|_{\epsilon = 0} T\Phi_{\epsilon} (q,\delta q), $$
where $(q,\delta q)$ are coordinates on $TM$. That is, $\widetilde{f}$ is the infinitesimal generator of the tangent lifted flow. The variational equation associated with \eqref{ODEmanifold} is the ODE associated with the tangent lifted vector field. In coordinates, 
\begin{equation}\label{VariationalODEmanifold}
\frac{d}{dt}(q,\delta q) = \widetilde{f}(q,\delta q).
\end{equation}
\begin{prop}\label{ManifoldQuadraticInvariantProp}
For integral curves $(q,p)$ of \eqref{AdjointODEmanifold} and $(q,\delta q)$ of \eqref{VariationalODEmanifold}, which cover the same curve $q$,
\begin{equation}\label{ManifoldVariationalQuadratic}
\frac{d}{dt} \Big\langle (q(t),p(t)), (q(t),\delta q(t))\Big\rangle_{q(t)} = 0. \end{equation}
\begin{proof}
Note that $(q(t),p(t)) \in T^*_{q(t)}M$ and $(q(t),\delta q(t)) \in T_{q(t)}M$ so the duality pairing is well-defined. Then,
\begin{align*}
\Big\langle (q(t),p(t)), (q(t),\delta q(t))\Big\rangle_{q(t)} &= \Big\langle (\Phi_{-t})^* (q(0),p(0)), T\Phi_t (q(0),\delta q(0))\Big\rangle_{q(t)} \\
&= \Big\langle (q(0),p(0)), T\Phi_{-t} \circ T\Phi_t (q(0),\delta q(0))\Big\rangle_{q(0)} \\
&= \Big\langle (q(0),p(0)), T(\Phi_{-t} \circ \Phi_t) (q(0),\delta q(0))\Big\rangle_{q(0)} \\
&= \Big\langle (q(0),p(0)), (q(0),\delta q(0))\Big\rangle_{q(0)}, 
\end{align*}
so the pairing is constant. 
\end{proof}
\end{prop}
\begin{remark}\label{ManifoldQuadraticInvariantSymplecticityRemark}
In the vector space case, we saw that the preservation of the quadratic invariant is implied by symplecticity. The above result is analogously implied by symplecticity, noting that the flow of the adjoint system is symplectic since $\widehat{f}$ is a Hamiltonian vector field. 
\end{remark}

Another conserved quantity for the adjoint system \eqref{AdjointODEmanifold} is the Hamiltonian, since the adjoint system corresponds to a time-independent Hamiltonian flow, $\frac{d}{dt} H = \Omega(X_H,X_H) = 0.$

Additionally, conserved quantities for adjoint systems are generated, via cotangent lift, by symmetries of the original ODE \eqref{ODEmanifold}, where we say that a vector field $g$ is a symmetry of the ODE $\dot{x} = h(x)$ if $[g,h] = 0$.

\begin{prop}\label{AdjointSystemSymmetry}
Let $g$ be a symmetry of \eqref{ODEmanifold}, i.e., $[g,f] = 0$. Then, its cotangent lift $\widehat{g}$ is a symmetry of \eqref{AdjointODEmanifold} and additionally, the function
$$ \langle \Theta, \widehat{g}\rangle $$
on $T^*M$ is preserved along the flow of $\widehat{f}$, i.e., under the flow of the adjoint system \eqref{AdjointODEmanifold}.
\begin{proof}
We first show that $\widehat{g}$ is a symmetry of \eqref{AdjointODEmanifold}, i.e., that $[\widehat{g},\widehat{f}] = 0$. To see this, we recall that the cotangent lift of the Lie bracket of two vector fields equals the Lie bracket of their cotangent lifts,
$$ \widehat{[g,f]} = [\widehat{g},\widehat{f}]. $$
Then, since $[g,f]=0$ by assumption, $[\widehat{g},\widehat{f}] = \widehat{[g,f]} = \widehat{0} = 0$. 

To see that $\langle \Theta, \widehat{g}\rangle$ is preserved along the flow of $\widehat{f}$, we have
\begin{align*}
\mathcal{L}_{\widehat{f}}\langle \Theta,\widehat{g}\rangle &= \langle \mathcal{L}_{\widehat{f}} \Theta, \widehat{g}\rangle + \langle \Theta, \mathcal{L}_{\widehat{f}} \widehat{g}\rangle = \langle 0, \widehat{g}\rangle + \langle \Theta, [\widehat{f},\widehat{g}]\rangle = 0,
\end{align*}
where we used that $\mathcal{L}_{\widehat{f}}\Theta = 0$ since $\widehat{f}$ is a cotangent lifted vector field.
\end{proof}
\end{prop}

\begin{remark}
The above proposition states when $[f,g]=0$, the Hamiltonian for the adjoint system associated with $g$, $\langle \Theta,\widehat{g}\rangle$, is preserved along the Hamiltonian flow corresponding to the Hamiltonian for the adjoint system associated with $f$, $\langle \Theta,\widehat{f}\rangle$, and vice versa. Note, $\langle \Theta, \widehat{g}\rangle$ can be interpreted as the momentum map corresponding to the action on $T^*M$ given by the flow of $\widehat{g}$.

The above proposition shows that (at least some) symmetries of the adjoint system \eqref{AdjointODEmanifold} can be found by cotangent lifting symmetries of the original ODE \eqref{ODEmanifold}. Additionally, the above proposition states that such cotangent lifted symmetries give rise to conserved quantities.
\end{remark}

In light of the above proposition, it is natural to ask the following question. Given a symmetry $G$ of the adjoint system \eqref{AdjointODEmanifold} (i.e., $[G,\widehat{f}] = 0$), does it arise from a cotangent lifted symmetry in the sense of Proposition \ref{AdjointSystemSymmetry}? In general, the answer is no. However, for a projectable vector field $G$ which is a symmetry of the adjoint system, its projection by $T\pi_{T^*M}$ to a vector field on $M$ does satisfy the assumptions of Proposition \ref{AdjointSystemSymmetry}. This gives the following partial converse to the above proposition.

\begin{prop}
Let $G$ be a projectable vector field on the bundle $\pi_{T^*M}: T^*M \rightarrow M$ which is a symmetry of \eqref{AdjointODEmanifold}, i.e., $[G,\widehat{f}] = 0$. Then, the pushforward vector field $g = T\pi_{T^*M}(G)$ on $M$ satisfies the assumptions of Proposition \ref{AdjointSystemSymmetry} and $T\pi_{T^*M}\widehat{g} = T\pi_{T^*M}G$.
\begin{proof}
Since $G$ is a projectable vector field on the cotangent bundle, $g = T\pi_{T^*M}G$ defines a well-defined vector field on $M$. Thus,
$$ [g,f] = [T\pi_{T^*M}G, T\pi_{T^*M}\widehat{f}] = T\pi_{T^*M}[G,\widehat{f}] = T\pi_{T^*M} 0 = 0, $$
so $g$ is a symmetry of \eqref{ODEmanifold}. Furthermore, we also have
$$ T\pi_{T^*M}\widehat{g} = T\pi_{T^*M} \widehat{(T\pi_{T^*M} G)} = T\pi_{T^*M} G. $$
\end{proof}
\end{prop}

The preceding proposition shows that, for the class of projectable symmetries of the adjoint system \eqref{AdjointODEmanifold}, it is always possible to find an associated symmetry of the original ODE \eqref{ODEmanifold} which, by Proposition \ref{AdjointSystemSymmetry}, corresponds to a Hamiltonian symmetry. Note that this implies that we can associate a conserved quantity $\langle \Theta, \widehat{g}\rangle$ to $G$, where $g = T\pi_{T^*M}G$. Furthermore, since $T\pi_{T^*M}\widehat{g} = T\pi_{T^*M}G$ and the canonical form $\Theta$ is a horizontal one-form, this implies that $\langle \Theta, G\rangle$ equals $\langle \Theta, \widehat{g}\rangle$ and hence, is conserved.

These two propositions show that symmetries of an ODE can be identified with equivalence classes of projectable symmetries of the associated adjoint system, where two projectable symmetries are equivalent if their difference lies in the kernel of $T\pi_{T^*M}$.

We also recall that the adjoint system \eqref{AdjointODEmanifold} formally arises from a variational principle. To do so, let $\Theta$ be the tautological one-form on $T^*M$. The action is defined to be
\begin{equation}\label{AdjointAction}
S[\psi] = \int_I [\psi^* \Theta - (H \circ \psi) dt ],
\end{equation}
where $\psi(t) = (q(t),p(t))$ is a curve on $T^*M$ over the interval $I= (t_0,t_1)$. We consider the variational principle $\delta S[\psi] = 0$, subject to variations which fix the endpoints $q(t_0)$, $q(t_1)$.

\begin{prop}\label{manifold adjoint ODE variational princple}
Let $\psi$ be a curve on $T^*M$ over the interval $I$. Then, $\psi$ is a stationary point of $S$ with respect to variations which fix $q(t_0)$, $q(t_1)$ if and only if \eqref{AdjointODEmanifold} holds.
\end{prop}

The proof of the above proposition is standard in the literature, so we will omit it.

\begin{remark}
It should be noted that although the fixed endpoint conditions on $q(t_0)$ and $q(t_1)$ in the variational principle formally obtains the correct equations of motion for the adjoint system, these boundary conditions are incompatible with the adjoint system, since it covers an ODE on the base manifold. From a theoretical perspective, this is not an obstruction to Proposition \ref{manifold adjoint ODE variational princple} since the equations of motion are obtained after enforcing the variational principle. However, from a numerical perspective, a variational principle with Type II boundary conditions fixing $q(t_0)$ and $p(t_1)$ is preferable for constructing variational integrators for adjoint systems. In Appendix \ref{typeIIvp appendix}, we develop an intrinsic Type II variational principle to incorporate these boundary conditions.
\end{remark}

\begin{remark}
In coordinates, the above action \eqref{AdjointAction} takes the form
$$ S = \int_{t_0}^{t_1}(\langle p,\dot{q}\rangle - \langle p, f(q)\rangle)dt, $$
which is the same coordinate expression as the action in the vector space case. 
\end{remark}

\subsubsection{Adjoint Systems with Augmented Hamiltonians}\label{AugmentedAdjointODESection}
In this section, we consider a class of modified adjoint systems, where some function on the base manifold $M$ is added to the Hamiltonian of the adjoint system. More precisely, let $H: T^*M \rightarrow \mathbb{R}, H(q,p) = \langle p, f(q)\rangle$ be the Hamiltonian of the previous section, corresponding to the ODE $\dot{q} = f(q)$. Let $L: M \rightarrow \mathbb{R}$ be a function on $M$. We identify $L$ with its pullback through $\pi_{T^*M}: T^*M \rightarrow M$. Then, we define the augmented Hamiltonian
\begin{align*}
H_L \equiv H+L: T^*M &\rightarrow \mathbb{R} \\
(q,p) &\mapsto H(q,p) + L(q) = \langle p, f(q)\rangle + L(q).
\end{align*}
We define the augmented adjoint system as the Hamiltonian system associated with $H_L$ relative to the canonical symplectic form $\Omega$ on $T^*M$,
\begin{equation}\label{AugmentedAdjointSystem}
i_{X_{H_L}}\Omega = dH_L.
\end{equation}

\begin{remark}
The motivation for such systems arises from adjoint sensitivity analysis and optimal control. For adjoint sensitivity analysis of a running cost function, one is concerned with the sensitivity of some functional 
$$ \int_{0}^{t} L(q)dt $$
along the flow of the ODE $\dot{q} = f(q)$. In the setting of optimal control, the goal is to minimize such a functional, constrained to curves satisfying the ODE (see, for example, \citet{AgCaFo2021}). We will discuss such applications in more detail in Section \ref{ApplicationsSection}.
\end{remark}

In coordinates, the augmented adjoint system \eqref{AugmentedAdjointSystem} takes the form
\begin{subequations}
\begin{align}
\dot{q}^A &= \frac{\partial H}{\partial p_A} = f^A(q), \label{AugmentedAdjointCoordinate1} \\
\dot{p}_A &= - \frac{\partial H}{\partial q^A} = -p_B \frac{\partial f^B(q)}{\partial q^A} - \frac{\partial L(q)}{\partial q^A}. \label{AugmentedAdjointCoordinate2}
\end{align}
\end{subequations}

We now prove various properties of the augmented adjoint system, analogous to the previous section. To start, first note that we can decompose the Hamiltonian vector field $X_{H_L}$ as follows. Let $\widehat{f}$ be the cotangent lift of $f$. Let $X_L \equiv X_{H_L} - \widehat{f}$. Then, observe that
$$ i_{X_L}\Omega = i_{X_{H_L}}\Omega - i_{\widehat{f}}\Omega = dH_L - dH = dL. $$
Thus, we have the decomposition $X_{H_L} = \widehat{f} + X_L$, where $\widehat{f}$ and $X_L$ are the Hamiltonian vector fields for $H$ and $L$, respectively. In coordinates,
$$ X_L = - \frac{\partial L}{\partial q^A} \frac{\partial}{\partial p_A}. $$
From the coordinate expression, we see that $X_L$ is a vertical vector field over the bundle $T^*M \rightarrow M$. We can also see this intrinsically, since $dL$ is a horizontal one-form on $T^*M$, $X_L$ satisfies $i_{X_L}\Omega = dL$, and $\Omega$ restricts to an isomorphism from vertical vector fields on $T^*M$ to horizontal one-forms on $T^*M$. Thus, it is immediate to see intrinsically that an analogous statement to Proposition \ref{LiftIntegralCurveProp} holds, since the flow of $\widehat{f}$ lifts the flow of $f$, while the flow of $X_L$ is purely vertical. That is, since $T\pi_{T^*M}X_L = 0$,
$$ T\pi_{T^*M}X_{H_L} = T\pi_{T^*M}\widehat{f} = f. $$
We can of course also see that the augmented adjoint system lifts the original ODE from the coordinate expression for the augmented adjoint system, \eqref{AugmentedAdjointCoordinate1}-\eqref{AugmentedAdjointCoordinate2}.

We now prove analogous statements to Propositions \ref{ManifoldQuadraticInvariantProp} and \ref{AdjointSystemSymmetry}, modified appropriately for the presence of $L$ in the augmented Hamiltonian.

\begin{prop}\label{AugmentedQuadraticInvariantProp}
Let $(q,p)$ be an integral curve of the augmented adjoint system \eqref{AugmentedAdjointSystem} and let $(q,\delta q)$ be an integral curve of the variational equation \eqref{VariationalODEmanifold}, covering the same curve $q$. Then,
$$ \frac{d}{dt} \langle p,\delta q\rangle = - \langle dL,\delta q\rangle. $$
\begin{remark}
Note that the variational equation associated with the above system is the same as in the nonaugmented case, equation \eqref{VariationalODEmanifold}, since augmenting $L$ to the Hamiltonian system only shifts the Hamiltonian vector field in the vertical direction.
\end{remark}
\begin{proof}
We will prove this in coordinates. We have the equations
\begin{align*}
\dot{p}_A &= -p_B \frac{\partial f^B}{\partial q^A} - \frac{\partial L}{\partial q^A}, \\
\frac{d}{dt}\delta q^B &= \frac{\partial f^B}{\partial q^A} \delta q^A.
\end{align*}
Then,
\begin{align*}
\frac{d}{dt} \langle p,\delta q\rangle &= \frac{d}{dt} p_A\delta q^A = \dot{p}_A \delta q^A +  p_B \frac{d}{dt}\delta q^B \\
&= -p_B \frac{\partial f^B}{\partial q^A}\delta q^A - \frac{\partial L}{\partial q^A}\delta q^A + p_B \frac{\partial f^B}{\partial q^A}\delta q^A \\
&= - \frac{\partial L}{\partial q^A}\delta q^A = - \langle dL,\delta q\rangle.
\end{align*}
\end{proof}
\end{prop}
\begin{remark}\label{Augmented ODE Quadratic Invariant Remark}
Interestingly, the above proposition states that in the augmented case, $\langle p,\delta q\rangle$ is no longer preserved but rather, its change measures the change of $L$ with respect to the variation $\delta q$. This may at first seem contradictory since both the augmented and nonaugmented Hamiltonian vector fields, $X_{H_L}$ and $X_H$, preserve $\Omega$, and as we noted previously in Remark \ref{ManifoldQuadraticInvariantSymplecticityRemark}, the preservation of the quadratic invariant is implied by symplecticity. However, upon closer inspection, there is no contradiction because the two cases have different first variations, where recall a first variation is a symmetry vector field of the Hamiltonian system and symplecticity can be stated as
$$ \frac{d}{dt}\Omega(V,W) = 0, $$
for first variation vector fields $V$ and $W$. In the nonaugmented case, the equations satisfied by the first variation of the momenta $p$ can be identified with $p$ itself, since the adjoint equation for $p$ is linear in $p$. On the other hand, in the augmented case, the adjoint equation for $p$, \eqref{AugmentedAdjointCoordinate2}, is no longer linear in $p$, rather, it is affine in $p$. Furthermore, the failure of this equation to be linear in $p$ is given precisely by $-dL$. Thus, in the augmented case, first variations in $p$ can no longer be identified with $p$, and this leads to the additional term $-\langle dL,\delta q\rangle$ in the above proposition.
\end{remark}

To prove an analogous statement to Proposition \ref{AdjointSystemSymmetry}, we need the additional assumption that the symmetry vector field $g$ leaves $L$ invariant, $\mathcal{L}_gL = 0$.
\begin{prop}
Let $g$ be a symmetry of the ODE $\dot{q} = f(q)$, i.e., $[g,f] = 0$. Additionally, assume that $g$ is a symmetry of $L$, i.e., $\mathcal{L}_gL = 0$. Then, its cotangent lift $\widehat{g}$ is a symmetry of the augmented adjoint system, $[\widehat{g},X_{H_L}] = 0$ and additionally, the function
$$ \langle \Theta, \widehat{g}\rangle $$
on $T^*M$ is preserved along the flow of $X_{H_L}$.
\begin{proof}
To see that $[\widehat{g},X_{H_L}] = 0$, note that with the decomposition $X_{H_L} = \widehat{f} + X_L$, we have
$$ [\widehat{g}, X_{H_L}] = [\widehat{g},\widehat{f}] + [\widehat{g}, X_L] = [\widehat{g},X_L], $$
where we used that $[\widehat{g},\widehat{f}] = \widehat{[g,f]} = 0$. To see that $[\widehat{g},X_L] = 0$, we note that $[\widehat{g},X_L]$ can be expressed
$$ [\widehat{g},X_L] = \mathcal{L}_{\widehat{g}}X_L = \mathcal{L}_{\widehat{g}} (\Omega^{-1}(dL)), $$
where we interpret $\Omega: T(T^*M) \rightarrow T^*(T^*M)$. Then, note that $\widehat{g}$ preserves $\Omega$ since $\widehat{g}$ is a cotangent lift and it also preserves $L$ (where, since we identify $L$ with its pullback through $\pi_{T^*M}$, this is equivalent to $g$ preserving $L$). More precisely, since we are identifying $L$ with its pullback $(\pi_{T^*M})^*L$, we have
$$ \mathcal{L}_{\widehat{g}}((\pi_{T^*M})^*L) =\langle (\pi_{T^*M})^* dL, \widehat{g}\rangle = \langle dL, T\pi_{T^*M}\widehat{g}\rangle =\langle dL, g\rangle = \mathcal{L}_gL = 0. $$
Hence, $\mathcal{L}_{\widehat{g}} (\Omega^{-1}(dL)) = 0$. One can also verify this in coordinates, and a direct computation yields
$$ [\widehat{g}, X_L] = \frac{\partial}{\partial q^A}\left( g^B(q) \frac{\partial L}{\partial q^B} \right) \frac{\partial}{\partial p_A}, $$
which vanishes since $\mathcal{L}_gL = 0$.

Now, to show that $\langle \Theta, \widehat{g}\rangle$ is preserved along the flow of $X_{H_L}$, compute
\begin{align*}
\mathcal{L}_{X_{H_L}} \langle \Theta, \widehat{g} \rangle &= \mathcal{L}_{\widehat{f}}\langle \Theta, \widehat{g} \rangle + \mathcal{L}_{X_L} \langle \Theta, \widehat{g} \rangle = \mathcal{L}_{X_L}\langle \Theta, \widehat{g} \rangle, 
\end{align*}
where we used that $\mathcal{L}_{\widehat{f}}\langle\Theta,\widehat{g}\rangle = 0$ by Proposition \ref{AdjointSystemSymmetry}. Now, we have
\begin{align*}
\mathcal{L}_{X_{H_L}} \langle \Theta, \widehat{g} \rangle &= \mathcal{L}_{X_L}\langle \Theta,\widehat{g}\rangle = \langle \mathcal{L}_{X_L}\Theta, \widehat{g}\rangle + \langle \Theta, \mathcal{L}_{X_L} \widehat{g}\rangle = \langle \mathcal{L}_{X_L}\Theta, \widehat{g}\rangle + \langle \Theta, \underbrace{[X_L,\widehat{g}]}_{=0}\rangle \\
&= \langle i_{X_L}d\Theta + d(i_{X_L}\Theta), \widehat{g} \rangle = \langle -i_{X_L}\Omega, \widehat{g}\rangle + \langle d(i_{X_L}\Theta),\widehat{g}\rangle \\
&= -\langle dL, \widehat{g}\rangle + \langle d(i_{X_L}\Theta),\widehat{g}\rangle.
\end{align*}
The first term above vanishes since $\mathcal{L}_g L = 0$. 
Furthermore, $\langle d(i_{X_L}\Theta),\widehat{g}\rangle = 0$ since $X_L$ is a vertical vector field while $\Theta$ is a horizontal one-form. Hence, $\mathcal{L}_{X_{H_L}}\langle\Theta,\widehat{g}\rangle = 0$.
\end{proof}
\end{prop}

\subsection{Adjoint Systems for DAEs via Presymplectic Mechanics}\label{AdjointDAESection}
In this section, we generalize the notion of adjoint system to the case where the base equation is a (semi-explicit) DAE. We will prove analogous results to the ODE case. However, more care is needed than the ODE case, since the DAE constraint introduces issues with solvability. As we will see, the adjoint system associated with a DAE is a presymplectic system, so we will approach the solvability of such systems through the presymplectic constraint algorithm.

We consider the following setup for a differential-algebraic equation. Let $M_d$ and $M_a$ be two manifolds, where we regard $M_d$ as the configuration space of the ``dynamical" or ``differential" variables and $M_a$ as the configuration space of the ``algebraic" variables. Let $\pi_{\Phi}: \Phi \rightarrow M_d \times M_a$ be a vector bundle over $M_d \times M_a$. Furthermore, let $\pi_d: M_d \times M_a \rightarrow M_d$ be the projection onto the first factor and let $\pi_{\overline{TM}_d}: \overline{TM}_d \rightarrow M_d \times M_a$ be the pullback bundle of the tangent bundle $\pi_{TM_d}: TM_d \rightarrow M_d$ by $\pi_d$, i.e., $\overline{TM}_d = \pi_d^*(TM_d)$. Then, a (semi-explicit) DAE is specified by a section $f \in \Gamma(\overline{TM}_d)$ and a section $\phi \in \Gamma(\Phi)$, via the system
\begin{subequations}
\begin{align}
\dot{q} &= f(q,u), \label{DAEa} \\
0 &= \phi(q,u), \label{DAEb}
\end{align}
\end{subequations}
where $(q,u)$ are coordinates on $M_d \times M_a$. We refer to $\overline{TM}_d$ as the differential tangent bundle, with coordinates $(q,u,v)$ and to $\Phi$ as the constraint bundle. 

\begin{remark}\label{DAE local solvability}
For the local solvability of \eqref{DAEa}-\eqref{DAEb}, regard $\phi$ locally as a map $\mathbb{R}^{\dim(M_d)} \times \mathbb{R}^{\dim(M_a)} \rightarrow  \mathbb{R}^{\normalfont{\text{rank}}(\Phi)}$. If $\partial\phi/\partial u$ is an isomorphism at a point $(q_0,u_0)$ where $\Phi(q_0,u_0)=0$, then by the implicit function theorem, one can locally solve $u = u(q)$ about $(q_0,u_0)$ such that $\phi(q,u(q))=0$, and subsequently solve the unconstrained differential equation $\dot{q} = f(q,u(q))$ locally. This is the case for semi-explicit index $1$ DAEs.

In order for the $\normalfont{\text{rank}}(\Phi) \times \dim(M_a)$ matrix $\partial\phi/\partial u(q_0,u_0)$ to be an isomorphism, it is necessary that  $\normalfont{\text{rank}}(\Phi) = \dim(M_a)$. However, we will make no such assumption, so as to treat the theory in full generality, allowing for, e.g., nonunique solutions. We will, however, assume that the $D\phi$ has constant rank (this corresponds to a fixed index DAE, which is the case if $\partial\phi/\partial u$ is a pointwise isomorphism) since we utilize the results of presymplectic geometry for constant rank presymplectic manifolds, as discussed in Section \ref{SymplecticGeometrySection}.
\end{remark}

Now, let $\overline{T^*M}_d$ be the pullback bundle of the cotangent bundle $T^*M_d$ by $\pi_d$, with coordinates $(q,u,p)$, which we refer to as the differential cotangent bundle. Furthermore, let $\Phi^*$ be the dual vector bundle to $\Phi$, with coordinates $(q,u,\lambda)$. Let $\overline{T^*M}_d \oplus \Phi^*$ be the Whitney sum of these two vector bundles over $M_d \times M_a$ with coordinates $(q,u,p,\lambda)$, which we refer to as the generalized phase space bundle. We define a Hamiltonian on the generalized phase space,
\begin{align*}
&H: \overline{T^*M}_d \oplus \Phi^* \rightarrow \mathbb{R}, \\
&H(q,u,p,\lambda) = \langle p, f(q,u) \rangle + \langle \lambda, \phi(q,u)\rangle.
\end{align*}
Let $\Omega_d$ denote the canonical symplectic form on $T^*M_d$, with coordinate expression $\Omega_d = dq \wedge dp$. We define a presymplectic form $\Omega_0$ on $\overline{T^*M}_d \oplus \Phi^*$ as follows: the pullback bundle admits the map $\tilde{\pi}_d: \overline{T^*M}_d \rightarrow T^*M_d$ which covers $\pi_d$ and acts as the identity on fibers and furthermore, the generalized phase space bundle admits the projection $\Pi: \overline{TM}_d^* \oplus \Phi^* \rightarrow \overline{TM}_d^*$, since the Whitney sum has the structure of a double vector bundle. Hence, we can pullback $\Omega_d$ along the sequence of maps
$$ \overline{T^*M}_d \oplus \Phi^* \overset{\Pi}{\longrightarrow} \overline{T^*M}_d \overset{\tilde{\pi}_d}{\longrightarrow} T^*M_d, $$
which allows us to define a two-form $\Omega_0 \equiv \Pi^* \circ \tilde{\pi}_d^* (\Omega_d)$ on the generalized phase space bundle. Clearly, $\Omega_0$ is closed as the pullback of a closed form. In general, $\Omega_0$ will be degenerate except in the trivial case where $M_a$ is empty and the fibers of $\Phi$ are the zero vector space. Hence, $\Omega_0$ is a presymplectic form. Note that since $\Pi$ acts by projection and $\tilde{\pi}_d$ acts as the identity on fibers, the coordinate expression for $\Omega_0$ on $\overline{T^*M}_d \oplus \Phi^*$ with coordinates $(q,u,p,\lambda)$ is the same as the coordinate expression for $\Omega_d$, $\Omega_0 = dq \wedge dp$. The various spaces and their coordinates are summarized in the diagram below.

\adjustbox{scale=0.8,center}{%
\begin{tikzcd}[column sep=3ex,row sep=5ex]
	{(q,u,p,\lambda) \in \overline{T^*M}_d\oplus \Phi^*} && {(q,u,\lambda)\in\Phi^*} && \Phi \\
	\\
	{(q,u,p) \in \overline{T^*M}_d} &&& {M_d \times M_a \ni (q,u)} &&& {\overline{TM}_d \ni (q,u,v)} \\
	\\
	{(q,p) \in T^*M_d} &&& {M_d \ni q} &&& {TM_d \ni (q,v)}
	\arrow["{\pi_d}", from=3-4, to=5-4]
	\arrow[from=5-7, to=5-4]
	\arrow[from=5-1, to=5-4]
	\arrow[from=3-7, to=3-4]
	\arrow[from=3-7, to=5-7]
	\arrow[from=3-1, to=5-1]
	\arrow[from=3-1, to=3-4]
	\arrow[from=1-3, to=3-4]
	\arrow[from=1-5, to=3-4]
	\arrow[from=1-1, to=3-1]
	\arrow[from=1-1, to=1-3]
\end{tikzcd}
}

We now define the adjoint system associated with the DAE \eqref{DAEa}-\eqref{DAEb} as the Hamiltonian system
\begin{equation}\label{DAE Adjoint System Intrinsic}
i_X\Omega_0 = dH. 
\end{equation}
Given a (generally, partially defined) vector field $X$ on the generalized phase space satisfying \eqref{DAE Adjoint System Intrinsic}, we say a curve $(q(t),u(t),p(t),\lambda(t))$ is a solution curve of \eqref{DAE Adjoint System Intrinsic} if it is an integral curve of $X$. 

Let us find a coordinate expression for the above system. Expressing our coordinates with indices $(q^i, u^a, p_j, \lambda_A)$, the left hand side of \eqref{DAE Adjoint System Intrinsic} along a solution curve has the expression
\begin{align*}
i_X\Omega_0 &= \left( \dot{q}^i \frac{\partial}{\partial q^i} + \dot{u}^a \frac{\partial}{\partial u^a} + \dot{p}_j \frac{\partial}{\partial p_j} + \dot{\lambda}_A \frac{\partial}{\partial \lambda_A} \right) \lrcorner\, dq^k \wedge dp_k \\
&= \dot{q}^i dp_i - \dot{p}_j dq^j.
\end{align*}
On the other hand, the right hand side of \eqref{DAE Adjoint System Intrinsic} has the expression
\begin{align*}
dH &= d\Big(p_if^i(q,u) + \lambda_A \phi^A(q,u)\Big) \\
&= f^i(q,u) dp_i + \left( p_i \frac{\partial f^i}{\partial q^j} + \lambda_A \frac{\partial \phi^A}{\partial q^j} \right) dq^j + \phi^A(q,u) d\lambda_A + \left( p_i \frac{\partial f^i}{\partial u^a} + \lambda_A \frac{\partial \phi^A}{\partial u^a} \right) du^a.
\end{align*}
Equating these expressions gives the coordinate expression for the adjoint DAE system,
\begin{subequations}
\begin{align}
\dot{q}^i &= f^i(q,u), \label{DAE Adjoint System Coord 1} \\
\dot{p}_j &= -p_i \frac{\partial f^i}{\partial q^j} - \lambda_A \frac{\partial \phi^A}{\partial q^j}, \label{DAE Adjoint System Coord 2}\\
0 &= \phi^A(q,u),\label{DAE Adjoint System Coord 3} \\
0 &= p_i \frac{\partial f^i}{\partial u^a} + \lambda_A \frac{\partial \phi^A}{\partial u^a}.\label{DAE Adjoint System Coord 4}
\end{align}
\end{subequations}

\begin{remark}\label{Adjoint DAE local solvability}
As mentioned in Remark \ref{DAE local solvability}, in the index $1$ case, one can locally solve the original DAE \eqref{DAE Adjoint System Coord 1} and \eqref{DAE Adjoint System Coord 3}. Viewing such a solution $(q,u)$ as fixed, one can subsequently locally solve for $\lambda$ in equation \eqref{DAE Adjoint System Coord 4} as a function of $p$, since $\partial \phi/\partial u$ is locally invertible. Substituting this into \eqref{DAE Adjoint System Coord 2} gives an ODE solely in the variable $p$, which can be solved locally. 

Stated another way, if the original DAE \eqref{DAEa}-\eqref{DAEb} is an index $1$ system, then the adjoint DAE system \eqref{DAE Adjoint System Coord 1}-\eqref{DAE Adjoint System Coord 4} is an index $1$ system with dynamical variables $(q,p)$ and algebraic variables $(u,\lambda)$. To see this, if one denotes the constraints for the adjoint system \eqref{DAE Adjoint System Coord 3} and \eqref{DAE Adjoint System Coord 4} as 
$$ 0 = \tilde{\phi}(q,u,p,\lambda) \equiv \begin{pmatrix} \phi^A(q,u) \\ p_i \frac{\partial f^i}{\partial u^a} + \lambda_A \frac{\partial \phi^A}{\partial u^a} \end{pmatrix}, $$
then the matrix derivative of $\tilde{\phi}$ with respect to the algebraic variables $(u,\lambda)$ can be locally expressed in block form as
$$ \begin{pmatrix} \partial\phi/\partial u & A \\ 0 & \partial\phi/\partial u \end{pmatrix}, $$
where the block $A$ has components given by the derivative of the right hand side of \eqref{DAE Adjoint System Coord 4} with respect to $u$. It is clear from the block triangular form of this matrix that it is pointwise invertible if $\partial\phi/\partial u$ is.
\end{remark}

\begin{remark}
It is clear from the coordinate expression \eqref{DAE Adjoint System Coord 1}-\eqref{DAE Adjoint System Coord 4} that a solution curve of the adjoint DAE system, if it exists, covers a solution curve of the original DAE system. 
\end{remark}

We now prove several results regarding the structure of the adjoint DAE system.

First, we show that the constraint equations \eqref{DAE Adjoint System Coord 3}-\eqref{DAE Adjoint System Coord 4} can be interpreted as the statement that the Hamiltonian $H$ has the same time dependence as the ``dynamical" Hamiltonian,
\begin{align*}
&H_d: \overline{T^*M}_d \oplus \Phi^* \rightarrow \mathbb{R}, \\
&H_d(q,u,p,\lambda) = \langle p, f(q,u) \rangle,
\end{align*}
when evaluated along a solution curve. 

\begin{prop}\label{Dynamical Hamiltonian Time Dependence Prop}
For a solution curve $(q,u,p,\lambda)$ of \eqref{DAE Adjoint System Intrinsic},
$$ \frac{d}{dt}H(q(t),u(t),p(t),\lambda(t)) =  \frac{d}{dt}H_d(q(t),u(t),p(t),\lambda(t)). $$
\begin{proof}
For brevity, all functions below are appropriately evaluated along the solution curve. We have
\begin{align*}
\frac{d}{dt}H &= \frac{\partial H}{\partial q^i} \dot{q}^i + \frac{\partial H}{\partial p_j} \dot{p}_j + \frac{\partial H}{\partial u^a} \dot{u}^a + \frac{\partial H}{\partial \lambda_A} \dot{\lambda}_A \\
&= \frac{\partial H}{\partial q^i} \dot{q}^i + \frac{\partial H}{\partial p_j} \dot{p}_j + \left(p_i \frac{\partial f^i}{\partial u^a} + \lambda_A \frac{\partial \phi^A}{\partial u^a} \right) \dot{u}^a + \phi^A \dot{\lambda}_A \\
&= \frac{\partial H}{\partial q^i} \dot{q}^i + \frac{\partial H}{\partial p_j} \dot{p}_j \\
&= \frac{\partial H_d}{\partial q^i} \dot{q}^i + \frac{\partial H_d}{\partial p_j} \dot{p}_j = \frac{d}{dt}H_d,
\end{align*}
where in the third equality, we used \eqref{DAE Adjoint System Coord 3} and \eqref{DAE Adjoint System Coord 4}.
\end{proof}
\end{prop}

\begin{remark}
A more geometric way to view the above proposition is as follows: note that if a partially-defined vector field $X$ exists such that $i_X\Omega_0 = dH$, then the change of $H$ in a given direction $Y$, at any point where $X$ is defined, can be computed as $dH(Y) = \Omega_0(X,Y)$. Observe that the kernel of $\Omega_0$ is locally spanned by $\partial/\partial u$, $\partial/\partial \lambda$, i.e., it is spanned by the coordinate vectors in the algebraic coordinates. Hence, the change of $H$ in the algebraic coordinate directions is zero. This justifies referring to $(u,\lambda)$ as ``algebraic" variables.
\end{remark}

We now show that the adjoint system \eqref{DAE Adjoint System Coord 1}-\eqref{DAE Adjoint System Coord 4} formally arises from a variational principle. To do so, let $\Theta_0$ be the pullback of the tautological one-form $\Theta_d$ on the cotangent bundle $T^*M_d$ by the maps $\Pi$ and $\tilde{\pi}_d$, $\Theta_0 = \Pi^*\circ \tilde{\pi}_d^* (\Theta_d)$. Of course, one has $\Omega_0 = -d\Theta_0$. Consider the action $S$ defined by
$$ S[\psi] = \int_I [\psi^* \Theta_0 - (H \circ \psi) dt ], $$
where $\psi(t) = (q(t),u(t),p(t),\lambda(t))$ is a curve on the generalized phase space bundle over the interval $I = (t_0,t_1)$. We consider the variational principle $\delta S[\psi] = 0$, subject to variations which fix the endpoints $q(t_0)$, $q(t_1)$.

\begin{prop}\label{AdjointDAEVariationalPrincipleProp}
Let $\psi$ be a curve on the generalized phase space bundle over the interval $I$. Then, $\psi$ is a stationary point of $S$ with respect to variations which fix $q(t_0)$, $q(t_1)$ if and only if \eqref{DAE Adjoint System Coord 1}-\eqref{DAE Adjoint System Coord 4} hold.
\begin{proof}
In $\psi = (q,u,p,\lambda)$ coordinates, the action has the expression
$$ S[q,u,p,\lambda] = \int_I \left(p_i\dot{q}^i - p_if^i(q,u) - \lambda_A \phi^A(q,u) \right) dt = \int_I \left(p_i (\dot{q}^i - f^i(q,u)) - \lambda_A\phi^A(q,u) \right) dt. $$
The variation of the action reads
\begin{align*}
\delta &S[q,u,p,\lambda]\cdot (\delta q, \delta u, \delta p, \delta \lambda) \\
&= \int_I \left[ \delta p_i (\dot{q}^i - f^i) + p_j \frac{d}{dt} \delta q^j - p_i\frac{\partial f^i}{\partial q^j} \delta q^j - \lambda_A \frac{\partial\phi^A}{\partial q^j}\delta q^j - \delta\lambda_A \phi^A + \left(-p_i \frac{\partial f^i}{\partial u^a} - \lambda_A \frac{\partial \phi^A}{\partial u^a}\right)\delta u^a \right] dt\\ 
&=\int_I \left[ \delta p_i (\dot{q}^i - f^i) - \left(\dot{p}_j + p_i\frac{\partial f^i}{\partial q^j} + \lambda_A \frac{\partial \phi^A}{\partial q^j} \right) \delta q^j - \delta\lambda_A \phi^A + \left(-p_i \frac{\partial f^i}{\partial u^a} - \lambda_A \frac{\partial \phi^A}{\partial u^a}\right)\delta u^a \right] dt,
\end{align*}
where we used integration by parts and the vanishing of the variations at the endpoints to drop any boundary terms. Clearly, if \eqref{DAE Adjoint System Coord 1}-\eqref{DAE Adjoint System Coord 4} hold, then $\delta S = 0$ for all such variations. Conversely, by the fundamental lemma of the calculus of variations, if $\delta S = 0$ for all such variations, then \eqref{DAE Adjoint System Coord 1}-\eqref{DAE Adjoint System Coord 4} hold. 
\end{proof}
\end{prop}

\begin{remark}
We will use the variational structure associated with the adjoint DAE system to construct numerical integrators in Section \ref{DiscreteAdjointSystemsSection}.
\end{remark}

We now prove a result regarding the conservation of a quadratic invariant, analogous to the case of cotangent lifted adjoint systems in the ODE case. To do this, we define the variational equations as the linearization of the DAE \eqref{DAEa}-\eqref{DAEb}. The coordinate expressions for the variational equations are obtained by taking the variation of equations \eqref{DAEa}-\eqref{DAEb} with respect to variations $(\delta q, \delta u$),
\begin{subequations}
\begin{align}
\dot{q}^i &= f^i(q,u), \label{DAEVariationalEqn1}\\
0 &= \phi^A(q,u), \label{DAEVariationalEqn2} \\
\frac{d}{dt}\delta q^i &= \frac{\partial f^i(q,u)}{\partial q^j}\delta q^j + \frac{\partial f^i(q,u)}{\partial u^a}\delta u^a, \label{DAEVariationalEqn3}\\
0 &= \frac{\partial \phi^A(q,u)}{\partial q^j} \delta q^j + \frac{\partial\phi^A(q,u)}{\partial u^a} \delta u^a. \label{DAEVariationalEqn4}
\end{align}
\end{subequations}

\begin{prop}\label{Quadratic Invariant Adjoint DAE Prop}
For a solution $(q,u,p,\lambda)$ of the adjoint DAE system \eqref{DAE Adjoint System Coord 1}-\eqref{DAE Adjoint System Coord 4} and a solution $(q,u,\delta q,\delta u)$ of the variational equations \eqref{DAEVariationalEqn1}-\eqref{DAEVariationalEqn4}, covering the same curve $(q,u)$, one has
$$ \frac{d}{dt} \langle p(t),\delta q(t) \rangle = 0. $$
\begin{proof}
This follows from a direct computation,
\begin{align*}
\frac{d}{dt} \langle p, \delta q \rangle &= \frac{d}{dt} \left(p_i \delta q^i\right) = \dot{p}_j \delta q^j + p_i \frac{d}{dt}\delta q^i \\
&= -p_i \frac{\partial f^i}{\partial q^j} \delta q^j - \lambda_A \frac{\partial\phi^A}{\partial q^j}\delta q^j + p_i \frac{\partial f^i}{\partial q^j}\delta q^j + p_i\frac{\partial f^i}{\partial u^a}\delta u^a \\
&= - \lambda_A \frac{\partial \phi^A}{\partial q^j}\delta q^j + p_i \frac{\partial f^i}{\partial u^a}\delta u^a \\
&= \left(\lambda_A \frac{\partial \phi^A}{\partial u^a} + p_i\frac{\partial f^i}{\partial u^a} \right)\delta u^a = 0,
\end{align*}
where we used \eqref{DAE Adjoint System Coord 2}, \eqref{DAEVariationalEqn3}, \eqref{DAEVariationalEqn4}, and \eqref{DAE Adjoint System Coord 4}.
\end{proof}
\end{prop}

\begin{remark}\label{PresymplecticityQuadraticInvariantDAERemark}
Although we proved the previous proposition in coordinates, it can be understood intrinsically through the presymplecticity of the adjoint DAE flow. To see this, assume a partially-defined vector field $X$ exists such that $i_X\Omega_0 = dH$. Then, the flow of $X$ preserves $\Omega_0$, which follows from
$$ \mathcal{L}_X\Omega_0 = i_X d\Omega_0 + d(i_X\Omega_0) = d(i_X\Omega_0) = d^2H = 0. $$
The coordinate expression for the preservation of the presymplectic form $\Omega_0 = dq^i \wedge dp_i$, with the appropriate choice of first variations, gives the previous proposition, analogous to the argument that we made in the symplectic (unconstrained) case. 

Additionally, as we will see in Section \ref{AdjointSensitivitySection}, Proposition \ref{Quadratic Invariant Adjoint DAE Prop} will provide a method for computing adjoint sensitivities.

These two observations are interesting when constructing numerical methods to compute adjoint sensitivities, since if we can construct integrators that preserve the presymplectic form, then it will preserve the quadratic invariant and hence, be suitable for computing adjoint sensitivities efficiently. 
\end{remark}

\begin{remark}\label{DAEVariationalEqnExistenceRemark}
For an index 1 DAE \eqref{DAEa}-\eqref{DAEb}, since $\partial\phi/\partial u$ is (pointwise) invertible for a fixed curve $(q,u)$, one can solve for $\delta u$ as a function of $\delta q$ in the variational equation \eqref{DAEVariationalEqn4} and substitute this into \eqref{DAEVariationalEqn3} to obtain an explicit ODE for $\delta q$. Hence, in the index 1 case, given a solution $(q,u)$ of the DAE \eqref{DAEa}-\eqref{DAEb} and an initial condition $\delta q(0)$ in the tangent fiber over $q(0)$, there is a corresponding (at least local) unique solution of the variational equations.
\end{remark}

\subsubsection{DAE Index and the Presymplectic Constraint Algorithm}\label{DAEIndexPCASection}
In this section, we relate the index of the DAE \eqref{DAEa}-\eqref{DAEb} to the number of steps for convergence in the presymplectic constraint algorithm associated with the adjoint DAE system \eqref{DAE Adjoint System Intrinsic}. In particular, we show that for an index 1 DAE, the presymplectic constraint algorithm for the associated adjoint DAE system converges after $\nu_P = 1$ step. Subsequently, we discuss how one can formally handle the more general index $\nu$ DAE case.

We consider again the presymplectic system given by the adjoint DAE system, $P = \overline{T^*M}_d \oplus \Phi^*$ equipped with the presymplectic form $\Omega_0 = dq \wedge dp$ and Hamiltonian $H(q,u,p,\lambda) = \langle p,f(q,u)\rangle + \langle \lambda, \phi(q,u)\rangle$, as discussed in the previous section. Our goal is to bound the number of steps in the presymplectic constraint algorithm $\nu_P$ for this presymplectic system in terms of the index $\nu$ of the underlying DAE \eqref{DAEa}-\eqref{DAEb}. 

Recall the presymplectic constraint algorithm discussed in Section \ref{SymplecticGeometrySection}. We first determine the primary constraint manifold $P_1$. Observe that since $\Omega_0 = dq \wedge dp$, we have the local expression $\text{ker}(\Omega_0)|_{(q,u,p,\lambda)} = \text{span}\{\partial/\partial u, \partial/\partial \lambda\}$. Thus, we require that 
\begin{align*}
\frac{\partial H}{\partial u} &= 0, \\
\frac{\partial H}{\partial \lambda} &= 0,
\end{align*}
i.e., $P_1$ consists of the points $(q,u,p,\lambda)$ such that
\begin{align*}
0 &= \frac{\partial H(q,u,p,\lambda)}{\partial u^a} = p_i \frac{\partial f^i(q,u)}{\partial u^a} + \lambda_A \frac{\partial\phi^A(q,u)}{\partial u^a}, \\
0 &= \frac{\partial H(q,u,p,\lambda)}{\partial \lambda^A} = \phi^A(q,u).
\end{align*}
These are of course the constraint equations \eqref{DAE Adjoint System Coord 3}-\eqref{DAE Adjoint System Coord 4} of the adjoint DAE system.

We now consider first the case when the DAE system \eqref{DAEa}-\eqref{DAEb} has index $\nu=1$ and subsequently, consider the general case $\nu \geq 1$.

\textbf{The Presymplectic Constraint Algorithm for $\nu=1$.} For the case $\nu=1$, we will show that the presymplectic constraint algorithm terminates after $1$ step, i.e., $\nu_P = \nu = 1$.

Now, assume that the DAE system \eqref{DAEa}-\eqref{DAEb} has index $\nu=1$, i.e., for each $(q,u) \in M_d \times M_a$ such that $\phi(q,u) = 0$, the matrix with $A^{th}$ row and $a^{th}$ column entry
$$ \frac{\partial\phi^A(q,u)}{\partial u^a} $$
is invertible. Observe that the definition of the presymplectic constraint algorithm, equation \eqref{Presymplectic Constraint Algorithm Sequence}, is local and hence, we seek a local coordinate expression for $\Omega_1 \equiv \Omega_0|_{P_1}$ and its kernel. 

Let $(q,u,p,\lambda) \in P_1$. In particular, $\phi(q,u) = 0$. Since $\partial\phi(q,u)/\partial u$ is invertible, by the implicit function theorem, one can locally solve for $u$ as a function of $q$, which we denote $u = u(q)$, such that $\phi(q,u(q)) = 0$. Then, one can furthermore locally solve for $\lambda$ as a function of $q$ and $p$ from the second constraint equation,
$$ \lambda_A(q,p) = - \left[ \left( \frac{\partial\phi(q,u(q))}{\partial u} \right)^{-1}\right]^a_A p_i \frac{\partial f^i(q,u(q))}{\partial u^a}. $$
Thus, we can coordinatize $P_1$ via coordinates $(q',p')$, where the inclusion $i_1: P_1 \hookrightarrow P$ is given by the coordinate expression
$$ i_1: (q',p') \mapsto (q', u(q'), p', \lambda(q',p')). $$
Then, one obtains the local expression for $\Omega_1$,
$$ \Omega_1 = i_1^*\Omega_0 = i_1^*(dq) \wedge i_1^*(dp) = dq' \wedge dp'. $$
This is clearly nondegenerate, i.e., $Z_p = 0$ for any $Z \in \text{ker}(\Omega_1), p \in P_1$, so the presymplectic constraint algorithm terminates, $P_2 = P_1$. We conclude that $\nu_P = 1$. 

To conclude the discussion of the index $1$ case, we obtain coordinate expressions for the resulting nondegenerate Hamiltonian system. The Hamiltonian on $P_1$ can be expressed as
$$ H_1(q',p') = H(i_1(q',p')) = \langle p', f(q',u(q'))\rangle + \langle \lambda(q',p'), \phi(q',u(q'))\rangle = \langle p', f(q',u(q'))\rangle. $$
Thus, with the coordinate expression $X = \dot{q}'^i \partial/\partial q'^i + \dot{p}'_i \partial/\partial p'_i$, Hamilton's equations $i_X \Omega_1 = dH_1$ can be expressed as
\begin{align*}
\dot{q}'^i &= \frac{\partial H_1}{\partial p'_i} = f^i(q',u(q')), \\
\dot{p}'_i &= - \frac{\partial H_1}{\partial q'^i} = -p'_j \frac{\partial f^j(q',u(q'))}{\partial q^i} -  p'_j \frac{\partial f^j(q',u(q'))}{\partial u^a} \frac{\partial u^a(q')}{\partial q'^i}.
\end{align*}
We will now show explicitly that this Hamiltonian system solves \eqref{DAE Adjoint System Coord 1}-\eqref{DAE Adjoint System Coord 4} along the submanifold $P_1$. Clearly, the latter two equations \eqref{DAE Adjoint System Coord 3}-\eqref{DAE Adjoint System Coord 4} are satisfied, by definition of $P_1$. So, we want to show that the first two equations \eqref{DAE Adjoint System Coord 1}-\eqref{DAE Adjoint System Coord 2} are satisfied. Using the second constraint equation \eqref{DAE Adjoint System Coord 4}, we have 
$$- p'_j \frac{\partial f^j(q',u(q'))}{\partial u^a} = \lambda_A(q',p') \frac{\partial\phi^A(q',u(q'))}{\partial u^a}.$$
Substituting this into the equation for $\dot{p}'_i$ above gives
$$ \dot{p}'_i = -p'_j \frac{\partial f^j(q',u(q'))}{\partial q^i} + \lambda_A(q',p') \frac{\partial\phi^A(q',u(q'))}{\partial u^a} \frac{\partial u^a(q')}{\partial q'^i}. $$
By the implicit function theorem, one has
$$ \frac{\partial\phi^A(q',u(q'))}{\partial u^a} \frac{\partial u^a(q')}{\partial q'^i} = - \frac{\partial \phi^A(q',u(q'))}{\partial q^i}. $$
Hence, the Hamiltonian system on $P_1$ can be equivalently expressed as
\begin{align*}
\dot{q}'^i &= f^i(q',u(q')), \\
\dot{p}'_i &= -p'_j \frac{\partial f^j(q',u(q'))}{\partial q^i} - \lambda_A(q',p') \frac{\partial \phi^A(q',u(q'))}{\partial q^i}.
\end{align*}
Thus, we have explicitly verified that \eqref{DAE Adjoint System Coord 1}-\eqref{DAE Adjoint System Coord 4} are satisfied along $P_1$. Note that since the presymplectic constraint algorithm terminates at $\nu_P = 1$, $X$ is guaranteed to be tangent to $P_1$. One can also verify this explicitly by computing the pushforward $Ti_1(X)$ and verifying that it annihilates the constraint functions whose zero level set defines $P_1$,
\begin{align*}
(q,u,p,\lambda) &\mapsto \phi^A(q,u), \\
(q,u,p,\lambda) &\mapsto p_i \frac{\partial f^i(q,u)}{\partial u^a} + \lambda_A \frac{\partial\phi^A(q,u)}{\partial u^a}.
\end{align*}

\begin{remark}\label{Index1ReductionPCA}
It is interesting to note that the Hamiltonian system $i_X\Omega_1 = dH_1$, which we obtained by forming the adjoint system of the underlying index 1 DAE and subsequently, reducing the index of the adjoint DAE system through the presymplectic constraint algorithm, can be equivalently obtained (at least locally) by first reducing the index of the underlying DAE and then forming the adjoint system.

More precisely, if one locally solves $\phi(q,u) = 0$ for $u = u(q)$, then the index 1 DAE can be reduced to an ODE,
$$ \dot{q} = f(q,u(q)). $$
Subsequently, we can form the adjoint system to this ODE, as discussed in Section \ref{AdjointManifoldSection}. The corresponding Hamiltonian is $H(q,p) = \langle p, f(q,u(q)) \rangle$, which is the same as $H_1$.

Thus, for the index 1 case, the process of forming the adjoint system and reducing the index commute. 
\end{remark}

\begin{remark}
In the language of the presymplectic constraint algorithm, Proposition~\ref{Dynamical Hamiltonian Time Dependence Prop} can be restated as the statement that the Hamiltonian $H$ and its first derivatives, restricted to the primary constraint manifold, agrees with the dynamical Hamiltonian $H_1$ and its first derivatives.
\end{remark}

\begin{remark}
An alternative view of the solution theory of the presymplectic adjoint DAE system \eqref{DAE Adjoint System Coord 1}-\eqref{DAE Adjoint System Coord 4} is through singular perturbation theory (see, for example, \citet{Be2007} and \citet{ChTr2021}). We proceed by writing \eqref{DAE Adjoint System Coord 1}-\eqref{DAE Adjoint System Coord 4} as
\begin{align*}
\dot{q} &= \frac{\partial H}{\partial p} = f(q,u), \\
\dot{p} &= -\frac{\partial H}{\partial q} = - [D_qf(q,u)]^*p - [D_q\phi(q,u)]^*\lambda, \\
0 &= \frac{\partial H}{\partial \lambda} = \phi(q,u), \\
0 &= - \frac{\partial H}{\partial u} = - [D_uf(q,u)]^*p - [D_u\phi(q,u)]^*\lambda.
\end{align*}
Applying a singular perturbation to the constraint equations yields the system 
\begin{align*}
\dot{q} &= \frac{\partial H}{\partial p}, \\
\dot{p} &= -\frac{\partial H}{\partial q}, \\
\epsilon \dot{u} &= \frac{\partial H}{\partial \lambda},  \\
\epsilon \dot{\lambda} &= - \frac{\partial H}{\partial u},
\end{align*}
where $\epsilon > 0$. Observe that this is a nondegenerate Hamiltonian system with $H(q,u,p,\lambda)$ as previously defined but with the modified symplectic form $\Omega_\epsilon = dq \wedge dp + \epsilon\, du \wedge d\lambda$.  Then, the above system can be expressed $i_{X_H}\Omega_\epsilon = dH$. In the language of perturbation theory, the primary constraint manifold for the presymplectic system is precisely the slow manifold of the singularly perturbed system. One can utilize techniques from singular perturbation theory to develop a solution theory for this system, using Tihonov's theorem, whose assumptions for this particular system depend on the eigenvalues of the algebraic Hessian $D_{u,\lambda}^2H$ (see, \citet{Be2007}). Although we will not elaborate on this here, this could be an interesting approach for the existence, stability, and approximation theory of such systems. In particular, the slow manifold integrators introduced in \citet{BuKl2020} may be relevant to their discretization. It is also interesting to note that for a solution $(q_\epsilon, p_\epsilon, u_\epsilon, \lambda_\epsilon)$ of the singularly perturbed system and a solution $(\delta q_\epsilon, \delta u_\epsilon)$ of the variational equations,
\begin{align*}
\frac{d}{dt} \delta q_\epsilon &= D_qf(q_\epsilon, u_\epsilon) \delta q_\epsilon + D_uf(q_\epsilon,u_\epsilon)\delta u_\epsilon, \\
\epsilon \frac{d}{dt} \delta u_\epsilon &= D_q\phi(q_\epsilon,u_\epsilon)\delta q_\epsilon + D_u\phi(q_\epsilon, u_\epsilon) \delta u_\epsilon,
\end{align*}
one has the perturbed adjoint variational quadratic conservation law
$$ \frac{d}{dt} \Big( \langle p_\epsilon, \delta q_\epsilon \rangle + \epsilon \langle \lambda_\epsilon, \delta u_\epsilon\rangle \Big) = 0, $$
which follows immediately from the preservation of $\Omega_\epsilon$ under the symplectic flow.
\end{remark}

\textbf{The Presymplectic Constraint Algorithm for General $\nu \geq 1$.} Note that for the general case, we assume that the index of the DAE is finite, $1 \leq \nu < \infty$.

In this case, there are two possible approaches to reduce the adjoint system: either form the adjoint system associated with the index $\nu$ DAE and then successively apply the presymplectic constraint algorithm or, alternatively, reduce the index of the DAE, form the adjoint system, and then apply the presymplectic constraint algorithm as necessary. 


Since we have already worked out the presymplectic constraint algorithm for the index 1 case, we will take the latter approach. Namely, we reduce an index $\nu$ DAE to an index $1$ DAE, and subsequently, apply the presymplectic constraint algorithm to the reduced index 1 DAE. Given an index $\nu$ DAE, it is generally possible to reduce the DAE to an index 1 DAE using the algorithm introduced in \citet{MaSo1993}. The process of index reduction is given by differentiating the equations of the DAE to reveal hidden constraints. Geometrically, the process of index reduction can be understood as the successive jet prolongation of the DAE and subsequent projection back onto the first jet (see, \citet{ReLiWi2001}).

Thus, given an index $\nu$ DAE $\dot{x} = \tilde{f}(x,y)$, $\tilde{\phi}(x,y) = 0$, we can, after $\nu-1$ reduction steps, transform it into an index 1 DAE of the form $\dot{q} = f(q,u)$, $\phi(q,u) = 0$. Subsequently, we can form the adjoint DAE system and apply one iteration of the presymplectic constraint algorithm to obtain the underlying nondegenerate dynamical system. If we let the $\nu_{R,P}$ denote the minimum number of DAE index reduction steps plus presymplectic constraint algorithm iterations necessary to take an index $\nu$ DAE and obtain the underlying nondegenerate Hamiltonian system associated with the adjoint, we have $\nu_{R,P} \leq \nu$.

\begin{remark}
Note that we could have reduced the index $\nu$ DAE to an explicit ODE after $\nu$ reduction steps, and subsequently, formed the adjoint. While this is formally equivalent to the above procedure by Remark \ref{Index1ReductionPCA}, we prefer to keep the DAE in index 1 form. This is especially preferable from the viewpoint of numerics: if one reduces an index 1 DAE to an ODE and attempts to apply a numerical integrator, it is generically the case that the discrete flow drifts off the constraint manifold. For this reason, it is preferable to develop numerical integrators for the index 1 adjoint DAE system directly to prevent constraint violation.
\end{remark}

\begin{example}[Hessenberg Index 2 DAE]
Consider a Hessenberg index 2 DAE, i.e., a DAE of the form
\begin{align*}
\dot{q} &= f(q,u), \\
0 &= g(q),
\end{align*}
where $(q,u) \in \mathbb{R}^n \times \mathbb{R}^m$, $f: \mathbb{R}^n \times \mathbb{R}^m \rightarrow \mathbb{R}^n$, $g: \mathbb{R}^n \rightarrow \mathbb{R}^m$, and $\frac{\partial g}{\partial q} \frac{\partial f}{\partial u}$ is pointwise invertible. We reduce this to an index 1 DAE \eqref{DAEa}-\eqref{DAEb} as follows. Let $M_d = g^{-1}(\{0\})$ be the dynamical configuration space which we will assume is a submanifold of $\mathbb{R}^n$. For example, this is true if $g$ is a constant rank map. Furthermore, let $M_a = \mathbb{R}^m$ be the algebraic configuration space. To reduce the index, we differentiate the constraint $g(q) = 0$ with respect to time. This is equivalent to enforcing that the dynamics are tangent to $M_d$. This gives
$$ 0 = \frac{\partial g^A(q)}{\partial q^i}{\dot{q}^i} = \frac{\partial g^A(q)}{\partial q^i} f^i(q,u) \equiv \phi^A(q,u). $$
Hence, we can form the semi-explicit index 1 system on $M_d \times M_a$ given by
\begin{align*}
\dot{q} &= f(q,u), \\
0 &= \phi(q,u).
\end{align*}
The above system is an index 1 DAE since $\frac{\partial \phi}{\partial u} = \frac{\partial g}{\partial q}\frac{\partial f}{\partial u}$ is pointwise invertible. 

We now form the adjoint DAE system associated with this index 1 DAE, \eqref{DAE Adjoint System Coord 1}-\eqref{DAE Adjoint System Coord 4}. Expressing the constraint in terms of $g$ and $f$, instead of $\phi$, gives
\begin{align*}
\dot{q}^i &= f^i(q,u), \\
\dot{p}_j &= -p_i \frac{\partial f^i(q,u)}{\partial q^j} - \lambda_A \left( \frac{\partial^2 g^A(q)}{\partial q^j \partial q^i} f^i(q,u) + \frac{\partial g^A(q)}{\partial q^i} \frac{\partial f^i(q,u)}{\partial q^j} \right), \\
0 &= \frac{\partial g^A(q)}{\partial q^i} f^i(q,u),\\
0 &= p_i \frac{\partial f^i(q,u)}{\partial u^a} + \lambda_A \left( \frac{\partial g^A(q)}{\partial q^i} \frac{\partial f^i(q,u)}{\partial u^a} \right).
\end{align*}
We can then apply one iteration of the presymplectic constraint algorithm, as discussed above in the index $\nu=1$ case, to obtain the underlying nondegenerate Hamiltonian dynamics. Restricting to the primary constraint manifold, using the first constraint equation to solve for $u=u(q)$ by the implicit function theorem and subsequently, using the second constraint equation to solve for $\lambda = \lambda(q,p)$ by inverting $\left(\frac{\partial g}{\partial q} \frac{\partial f}{\partial u}\right)^T$, gives the Hamiltonian system
\begin{align*}
\dot{q}'^i &= f^i(q',u(q')), \\
\dot{p}'_j &= -p'_i \frac{\partial f^i(q',u(q'))}{\partial q^j} - \lambda_A(q',p') \left( \frac{\partial^2 g^A(q')}{\partial q^j \partial q^i} f^i(q',u(q')) + \frac{\partial g^A(q')}{\partial q^i} \frac{\partial f^i(q',u(q'))}{\partial q^j} \right).
\end{align*}
\end{example}

\subsubsection{Adjoint Systems for DAEs with Augmented Hamiltonians}\label{AugmentedAdjointDAESection}
In Section \ref{AugmentedAdjointODESection}, we augmented the adjoint ODE Hamiltonian by some function $L$. In this section, we do analogously for the adjoint DAE system.

To begin, let $H(q,u,p,\lambda) = \langle p,f(q,u)\rangle + \langle \lambda,\phi(q,u)\rangle$ be the Hamiltonian on the generalized phase space bundle corresponding to the DAE $\dot{q}=f(q,u)$, $0 = \phi(q,u)$, and let $L: M_d \times M_a \rightarrow \mathbb{R}$ be the function that we would like to augment. We identify $L$ with its pullback through $\overline{T^*M}_d \oplus \Phi^* \rightarrow M_d \times M_a$. Then, we define the augmented Hamiltonian
\begin{align*}
H_L \equiv H+L: \overline{T^*M}_d \oplus \Phi^* &\rightarrow \mathbb{R} \\
(q,u,p,\lambda) &\mapsto H(q,u,p,\lambda) + L(q,u).
\end{align*}
We define the augmented adjoint DAE system as the presymplectic system
\begin{equation}\label{AugmentedAdjointDAEIntrinsic}
i_{X_{H_L}}\Omega_0 = dH_L.
\end{equation}
A direct calculation yields the coordinate expression, along an integral curve of such a (generally, partially-defined) vector field $X_{H_L}$,
\begin{subequations}
\begin{align}
\dot{q}^i &= f^i(q,u), \label{Augmented Adjoint DAE System Coord 1} \\
\dot{p}_j &= -p_i \frac{\partial f^i}{\partial q^j} - \lambda_A \frac{\partial \phi^A}{\partial q^j} - \frac{\partial L}{\partial q^j}, \label{Augmented Adjoint DAE System Coord 2}\\
0 &= \phi^A(q,u),\label{Augmented Adjoint DAE System Coord 3} \\
0 &= p_i \frac{\partial f^i}{\partial u^a} + \lambda_A \frac{\partial \phi^A}{\partial u^a} + \frac{\partial L}{\partial u^a}.\label{Augmented Adjoint DAE System Coord 4}
\end{align}
\end{subequations}
\begin{remark}
Observe that if the base DAE \eqref{DAEa}-\eqref{DAEb} has index 1, then the above system has index 1 by the exact same argument given in the nonaugmented case. After reduction by applying the presymplectic constraint algorithm and solving for $u$ as a function of $q$ and $\lambda$ as a function of $(q,p)$, the underlying nondegenerate Hamiltonian system on the primary (final) constraint manifold corresponds to the Hamiltonian 
$$(H_L)_1(q',p') = \langle p',f(q',u(q'))\rangle + L(q',u(q')),$$
which is the adjoint Hamiltonian for the ODE $\dot{q}' = f(q',u(q'))$, augmented by $L(q',u(q'))$.

However, as we will discuss in Section \ref{OCPSection}, it is not uncommon in optimal control problems for $\partial\phi/\partial u$ to be singular, but the presence of $\int L\, dt$ in the minimization objective may uniquely specify the singular degrees of freedom.
\end{remark}

We now prove an analogous proposition to Proposition \ref{Quadratic Invariant Adjoint DAE Prop}, modified by the presence of $L$ in the Hamiltonian. We again consider the variational equations \eqref{DAEVariationalEqn1}-\eqref{DAEVariationalEqn4} associated with the base DAE \eqref{DAEa}-\eqref{DAEb}, which for simplicity we express in matrix derivative notation as
\begin{subequations}
\begin{align}
\dot{q} &= f(q,u), \label{Full q,u DAE variation 1} \\
0 &= \phi(q,u), \label{Full q,u DAE variation 2} \\
\frac{d}{dt} \delta q &= D_qf(q,u)\delta q + D_uf(q,u)\delta u, \label{Full q,u DAE variation 3} \\
0 &= D_q\phi(q,u)\delta q + D_u\phi(q,u)\delta u. \label{Full q,u DAE variation 4}
\end{align}
\end{subequations}

\begin{prop}\label{Quadratic Invariant Augmented DAE Adjoint Prop}
For a solution $(q,u,p,\lambda)$ of the augmented adjoint DAE system \eqref{Augmented Adjoint DAE System Coord 1}-\eqref{Augmented Adjoint DAE System Coord 4} and a solution $(q,u,\delta q, \delta u)$ of the variational equations \eqref{Full q,u DAE variation 1}-\eqref{Full q,u DAE variation 4}, covering the same solution $(q,u)$ of the base DAE \eqref{DAEa}-\eqref{DAEb},
\begin{equation}\label{Quadratic Invariant Augmented DAE Adjoint Eqn}
\frac{d}{dt} \langle p,\delta q \rangle = -\langle \nabla_qL, \delta q\rangle - \langle \nabla_uL,\delta u\rangle.
\end{equation}
\begin{proof}
This follows from a direct computation:
\begin{align*}
\frac{d}{dt} \langle p,\delta q\rangle &= \langle \dot{p},\delta q\rangle + \langle p, \frac{d}{dt}\delta q\rangle \\
&= - \langle [D_qf]^*p, \delta q \rangle - \langle [D_q\phi]^*\lambda, \delta q\rangle - \langle \nabla_qL, \delta q\rangle + \langle p, D_qf \delta q\rangle + \langle p, D_uf \delta u\rangle \\
&= - \langle \lambda, D_q\phi\delta q\rangle - \langle \nabla_qL, \delta q\rangle + \langle p, D_uf \delta u\rangle \\
&= \langle \lambda, D_u\phi \delta u\rangle - \langle \nabla_qL, \delta q\rangle + \langle p, D_uf \delta u\rangle \\
&= - \langle \nabla_qL, \delta q\rangle + \langle [D_u\phi]^*\lambda + [D_uf]^*p, \delta u\rangle \\
&= -\langle \nabla_qL, \delta q\rangle - \langle \nabla_uL,\delta u\rangle,
\end{align*}
where in the fourth equality above we used \eqref{Full q,u DAE variation 4} and in the sixth equality above we used \eqref{Augmented Adjoint DAE System Coord 4}.
\end{proof}
\end{prop}
\begin{remark}
Analogous to the ODE case discussed in Remark \ref{Augmented ODE Quadratic Invariant Remark}, we remark that for the nonaugmented adjoint DAE system \eqref{DAE Adjoint System Coord 1}-\eqref{DAE Adjoint System Coord 4}, we have preservation of $\langle p, \delta q\rangle$ by virtue of presymplecticity. On the other hand, for the augmented adjoint DAE system, despite preserving the same presymplectic form, the change of $\langle p,\delta q\rangle$ now measures the change in $L$ with respect to variations in $q$ and $u$. This can be understood from the fact that the adjoint equations for $(p,\lambda)$ in the nonaugmented case, \eqref{DAE Adjoint System Coord 2} and \eqref{DAE Adjoint System Coord 4}, are linear in $(p,\lambda)$, so that one can identify first variations in $(p,\lambda)$ with $(p,\lambda)$; whereas, in the augmented case, equations \eqref{Augmented Adjoint DAE System Coord 2} and \eqref{Augmented Adjoint DAE System Coord 4} are affine in $(p,\lambda)$, so such an identification cannot be made. Furthermore, the failure of \eqref{Augmented Adjoint DAE System Coord 2} and \eqref{Augmented Adjoint DAE System Coord 4} to be linear in $(p,\lambda)$ are given precisely by $\nabla_qL$ and $\nabla_uL$, respectively. Thus, in the augmented case, this leads to the additional terms $-\langle \nabla_uL,\delta q\rangle - \langle \nabla_q L,\delta u\rangle$ in equation \eqref{Quadratic Invariant Augmented DAE Adjoint Eqn}.
\end{remark}

\section{Applications}\label{ApplicationsSection}

\subsection{Adjoint Sensitivity Analysis for Semi-explicit Index 1 DAEs}\label{AdjointSensitivitySection}
In this section, we discuss how one can utilize adjoint systems to compute sensitivities. We will split this into four cases; namely, we want to compute sensitivities for ODEs or DAEs (we will focus on index 1 DAEs), and whether we are computing the sensitivity of a terminal cost or the sensitivity of a running cost.

The relevant adjoint system used to compute sensitivities in all four cases can be summarized:
\renewcommand{\arraystretch}{1.3}
\begin{center}
\begin{tabular}{ |c|c|c| } 
 \hline
  & Terminal Cost & Running Cost \\ 
  \hline
 ODE & Adjoint ODE System \eqref{AdjointODEa}-\eqref{AdjointODEb} & Augmented Adjoint ODE System \eqref{AugmentedAdjointCoordinate1}-\eqref{AugmentedAdjointCoordinate2} \\ 
 DAE & Adjoint DAE System \eqref{DAE Adjoint System Coord 1}-\eqref{DAE Adjoint System Coord 4} & Augmented Adjoint DAE System \eqref{Augmented Adjoint DAE System Coord 1}-\eqref{Augmented Adjoint DAE System Coord 4} \\ 
 \hline
\end{tabular}
\end{center}
Note that in our calculations below, the top row (the ODE case) can be formally obtained from the bottom row (the DAE case) simply by ignoring the algebraic variables $(u,\lambda)$ and letting the constraint function $\phi$ be identically zero. Thus, we will focus on the bottom row, i.e., computing sensitivities of a terminal cost function and of a running cost function, subject to a DAE constraint. In both cases, we will first show how the adjoint sensitivity can be derived using a traditional variational argument. Subsequently, we will show how the adjoint sensitivity can be derived more simply by using Propositions \ref{Quadratic Invariant Adjoint DAE Prop} and \ref{Quadratic Invariant Augmented DAE Adjoint Prop}.

\textbf{Adjoint Sensitivity of a Terminal Cost.}
Consider the DAE $\dot{q} = f(q,u)$, $0 = \phi(q,u)$ as in Section \ref{AdjointDAESection}. We will assume that $M_d$ is a vector space and additionally, that the DAE has index 1. We would like to extract the gradient of a terminal cost function $C(q(t_f))$ with respect to the initial condition $q(0) = \alpha$, i.e., we want to extract the sensitivity of $C(q(t_f))$ with respect to an infinitesimal perturbation in the initial condition, given by $\nabla_\alpha C(q(t_f))$. Consider the functional $J$ defined by
$$ J = C(q(t_f)) - \langle p_0, q(0) - \alpha\rangle - \int_0^{t_f} [\langle p, \dot{q}-f(q,u)\rangle - \langle \lambda, \phi(q,u)\rangle ]dt. $$
Observe that for $(q,u)$ satisfying the given DAE with initial condition $q(0) = \alpha$, $J$ coincides with $C(q(t_f))$. We think of $p_0$ as a free parameter. For simplicity, we will use matrix derivative notation instead of indices. Computing the variation of $J$ yields
\begin{align*}
\delta J &= \langle \nabla_qC(q(t_f)), \delta q(t_f)\rangle - \langle p_0, \delta q(0) - \delta \alpha\rangle \\ 
&\qquad - \int_0^{t_f} \Big[ \langle p, \frac{d}{dt} \delta q - D_qf(q,u)\delta q\rangle - \langle p, D_uf(q,u)\delta u\rangle - \langle \lambda, D_q\phi(q,u)\delta q + D_u\phi(q,u)\delta u\rangle \Big]dt.
\end{align*}
Integrating by parts in the term containing $\frac{d}{dt}\delta q$ and restricting to a solution $(q,u,p,\lambda)$ of the adjoint DAE system \eqref{DAE Adjoint System Coord 1}-\eqref{DAE Adjoint System Coord 4} yields
\begin{align*}
\delta J &= \langle \nabla_qC(q(t_f)) - p(t_f), \delta q(t_f)\rangle - \langle p_0, \delta \alpha\rangle + \langle p(0) - p_0, \delta q(0)\rangle .
\end{align*}
We enforce the endpoint condition $p(t_f) = \nabla_qC(q(t_f))$ and choose $p_0 = p(0)$, which yields
$$ \delta J = \langle p(0), \delta\alpha\rangle. $$
Hence, the sensitivity of $C(q(t_f))$ is given by
$$ p(0) = \nabla_\alpha J = \nabla_\alpha C(q(t_f)), $$
with initial condition $q(0) = \alpha$ and terminal condition $p(t_f) = \nabla_qC(q(t_f))$. Thus, the adjoint sensitivity can be computed by setting the terminal condition on $p(t_f)$ above and subsequently, solving for the momenta $p$ at time $0$. In order for this to be well-defined, we have to verify that the given initial and terminal conditions lie on the primary constraint manifold $P_1$. However, as discussed in Section \ref{DAEIndexPCASection}, since the DAE has index 1, we can always solve for the algebraic variables $u = u(q)$ and $\lambda = \lambda(q,p)$ and thus, we are free to choose the initial and terminal values of $q$ and $p$, respectively. For higher index DAEs, one has to ensure that these conditions are compatible with the final constraint manifold. For example, this is done in \cite{CaLiPeSe2003} in the case of Hessenberg index 2 DAEs. Alternatively, at least theoretically, for higher index DAEs, one can reduce the DAE to an index 1 DAE and then the above discussion applies, however, this reduction may fail in practice due to numerical cancellation.

Note that the above adjoint sensitivity result is also a consequence of the preservation of the quadratic invariant $\langle p,v\rangle$ as in Proposition \ref{Quadratic Invariant Adjoint DAE Prop}. From this proposition, one has that
$$ \langle p(t_f), \delta q(t_f) \rangle = \langle p(0), \delta q(0)\rangle, $$
where $\delta q$ satisfies the variational equations. Setting $p(t_f) = \nabla_q C(q(t_f))$ and $\delta q(0) = \delta\alpha$ gives the same result. As mentioned in Remark \ref{PresymplecticityQuadraticInvariantDAERemark}, this quadratic invariant arises from the presymplecticity of the adjoint DAE system. Thus, a numerical integrator which preserves the presymplectic structure is desirable for computing adjoint sensitivities, as it exactly preserves the quadratic invariant that allows the adjoint sensitivities to be accurately and efficiently computed. We will discuss this in more detail in Section \ref{DiscreteAdjointSystemsSection}.

\textbf{Adjoint Sensitivity of a Running Cost.} Again, consider an index 1 DAE $\dot{q} = f(q,u)$, $0 = \phi(q,u)$. We would like to extract the sensitivity of a running cost function 
$$ \int_0^{t_f} L(q,u) dt,$$
where $L: M_d \times M_a \rightarrow \mathbb{R}$, with respect to an infinitesimal perturbation in the initial condition $q(0) = \alpha$. Consider the functional $J$ defined by
$$ J = -\langle p_0, q(0)-\alpha\rangle + \int_{0}^{t_f}[L(q,u) + \langle p,f(q,u) - \dot{q}\rangle + \langle \lambda, \phi(q,u)\rangle]dt. $$
Observe that when the DAE is satisfied with initial condition $q(0)=\alpha$, $J = \int_0^{t_f}L\, dt$. Now, we would to compute the implicit change in $\int_0^{t_f}L\,dt$ with respect to a perturbation $\delta\alpha$ in the initial condition. Taking the variation in $J$ yields
\begin{align*}
\delta J &= -\langle p_0, \delta q(0)-\delta \alpha\rangle\\
&\qquad + \int_0^{t_f} \Big[\langle \nabla_qL, \delta q\rangle + \langle \nabla_uL, \delta u\rangle + \langle p, D_qf \delta q - \frac{d}{dt}\delta q \rangle  + \langle p, D_uf \delta u\rangle + \langle \lambda, D_q\phi \delta q + D_u\phi \delta u \rangle \Big]dt \\
&= -\langle p_0, \delta q(0)-\delta \alpha\rangle - \langle p(t_f), \delta q(t_f)\rangle + \langle p(0), \delta q(0)\rangle \\ 
&\qquad + \int_0^{t_f} \Big[ \langle \nabla_qL + [D_qf]^*p + [D_q\phi]^*\lambda + \dot{p}, \delta q\rangle + \langle \nabla_uL + [D_uf]^*p + [D_u\phi]^*\lambda, \delta u\rangle \Big] dt.
\end{align*}
Restricting to a solution $(q,u,p,\lambda)$ of the augmented adjoint DAE system \eqref{Augmented Adjoint DAE System Coord 1}-\eqref{Augmented Adjoint DAE System Coord 4}, setting the terminal condition $p(t_f) = 0$, and choosing $p_0 = p(0)$ gives $ \delta J = \langle p(0), \delta \alpha\rangle.$ Hence, the implicit sensitivity of $\int_{0}^{t_f} L\, dt$ with respect to a change $\delta\alpha$ in the initial condition is given by 
$$ p(0) = \delta_\alpha J = \delta_\alpha \int_0^{t_f}L(q,u)dt.$$
Thus, the adjoint sensitivity of a running cost functional with respect to a perturbation in the initial condition can be computed by using the augmented adjoint DAE system \eqref{Augmented Adjoint DAE System Coord 1}-\eqref{Augmented Adjoint DAE System Coord 4} with terminal condition $p(t_f) = 0$ to solve for the momenta $p$ at time 0. 

Note that the above adjoint sensitivity result can be obtained from Proposition \ref{Quadratic Invariant Augmented DAE Adjoint Prop} as follows. We write equation \eqref{Quadratic Invariant Augmented DAE Adjoint Eqn} as
$$ \frac{d}{dt} \langle p,\delta q\rangle = -\langle dL, (\delta q, \delta u)\rangle, $$
to highlight that the right hand side measures the total induced variation of $L$. Now, we integrate this equation from $0$ to $t_f$, which gives
$$ \langle p(t_f), \delta q(t_f)\rangle - \langle p(0),\delta q(0)\rangle = - \int_0^{t_f} \langle dL, (\delta q, \delta u)\rangle dt. $$ 
Since we want to determine the change in the running cost functional with respect to a perturbation in the initial condition, we set $p(t_f) = 0$ which yields
$$ \langle p(0),\delta q(0)\rangle = \int_0^{t_f} \langle dL, (\delta q, \delta u)\rangle dt. $$
The right hand side is the total change induced on the running cost functional, whereas the left hand side tells us how this change is implicitly induced from a perturbation $\delta q(0)$ in the initial condition. Note that a perturbation in the initial condition $\delta q(0)$ will generally induce perturbations in both $q$ and $u$, according to the variational equations. Such a curve $(\delta q, \delta u)$ satisfying the variational equations exists in the index 1 case as noted in Remark \ref{DAEVariationalEqnExistenceRemark}. Thus, we arrive at the same conclusion as the variational argument: $p(0)$ is the desired adjoint sensitivity.

To summarize, adjoint sensitivities for terminal and running costs can be computed using the properties of adjoint systems, such as the various aforementioned propositions regarding $\frac{d}{dt} \langle p, \delta q\rangle$, which is zero in the nonaugmented case and measures the variation of $L$ in the augmented case. In the case of a terminal cost, one sets an inhomogeneous terminal condition $p(t_f) = \nabla_qC(q(t_f))$ and backpropagates the momenta through the nonaugmented adjoint DAE system \eqref{DAE Adjoint System Coord 1}-\eqref{DAE Adjoint System Coord 4} to obtain the sensitivity $p(0)$. On the other hand, in the case of a running cost, one sets a homogeneous terminal condition $p(t_f) = 0$ and backpropagates the momenta through the augmented adjoint DAE system \eqref{Augmented Adjoint DAE System Coord 1}-\eqref{Augmented Adjoint DAE System Coord 4} to obtain the sensitivity $p(0)$.

The various propositions used to derive the above adjoint sensitivity results are summarized below. We also include the ODE case, since it follows similarly.
\renewcommand{\arraystretch}{1.7}
\begin{center}
\begin{tabular}{ |c|c|c| } 
 \hline
  & Terminal Cost & Running Cost \\ 
  \hline ODE & Proposition \ref{ManifoldQuadraticInvariantProp}, $\frac{d}{dt}\langle p,\delta q\rangle = 0$ & Proposition \ref{AugmentedQuadraticInvariantProp}, $\frac{d}{dt} \langle p,\delta q\rangle = - \langle dL, \delta q\rangle$ \\ 
 DAE & Proposition \ref{Quadratic Invariant Adjoint DAE Prop}, $\frac{d}{dt}\langle p,\delta q\rangle = 0$  & Proposition \ref{Quadratic Invariant Augmented DAE Adjoint Prop}, $\frac{d}{dt}\langle p,\delta q\rangle = - \langle dL, (\delta q,\delta u)\rangle$ \\ 
 \hline
\end{tabular}
\end{center}
In Section \ref{DiscreteAdjointSystemsSection}, we will construct integrators that admit discrete analogues of the above propositions, and hence, are suitable for computing discrete adjoint sensitivities. 

\subsection{Structure-Preserving Discretizations of Adjoint Systems}\label{DiscreteAdjointSystemsSection}
In this section, we utilize the Galerkin Hamiltonian variational integrators of \citet{LeZh2011} to construct structure-preserving integrators which admit discrete analogues of Propositions \ref{ManifoldQuadraticInvariantProp}, \ref{AugmentedQuadraticInvariantProp}, \ref{Quadratic Invariant Adjoint DAE Prop}, and \ref{Quadratic Invariant Augmented DAE Adjoint Prop}, and are therefore suitable for numerical adjoint sensitivity analysis. For brevity, the proofs of these discrete analogues can be found in Appendix \ref{QuadraticConservationProofs}.

We start by recalling the construction of Galerkin Hamiltonian variational integrators as introduced in \citet{LeZh2011}. We assume that the base manifold $Q$ is a vector space and thus, we have the identification $T^*Q \cong Q \times Q^*$. To construct a variational integrator for a Hamiltonian system on $T^*Q$, one starts with the exact Type II generating function
$$ H^+_{d,\text{exact}}(q_0,p_1) = \ext\left[\langle p_1,q_1\rangle - \int_0^{\Delta t} [\langle p,\dot{q}\rangle - H(q,p)]dt\right], $$
where one extremizes over $C^2$ curves on the cotangent bundle satisfying $q(0) = q_0, p(\Delta t) = p_1$. This is a Type II generating function in the sense that it defines a symplectic map $(q_0,p_1) \mapsto (q_1, p_0)$ by $q_1 = D_2H^+_{d,\text{exact}}(q_0,p_1)$, $p_0 = D_1H^+_{d,\text{exact}}(q_0,p_1)$.

To approximate this generating function, one approximates the integral above using a quadrature rule and extremizes the resulting expression over a finite-dimensional subspace satisfying the prescribed boundary conditions. This yields the Galerkin discrete Hamiltonian
$$ H_d^+(q_0,p_1) = \ext \left[\langle p_1, q_1\rangle - \Delta t \sum_i b_i \Big( \langle P^i, V^i\rangle - H(Q^i,P^i) \Big) \right], $$
where $\Delta t > 0$ is the timestep, $q_0, q_1, p_0, p_1$ are numerical approximations to $q(0), q(\Delta t), p(0), p(\Delta t)$, respectively, $b_i > 0$ are quadrature weights corresponding to quadrature nodes $c_i \in [0,1]$, $Q^i$ and $P^i$ are internal stages representing $q(c_i\Delta t), p(c_i\Delta t)$, respectively, and $V$ is related to $Q$ by $Q^i = q_0 + \Delta t \sum_j a_{ij}V^j$, where the coefficients $a_{ij}$ arise from the choice of function space. The expression above is extremized over the internal stages $Q^i, P^i$ and subsequently, one applies the discrete right Hamilton's equations 
\begin{align*}
q_1 &= D_2H_d^+(q_0,p_1), \\
p_0 &= D_1H_d^+(q_0,p_1),
\end{align*}
to obtain a Galerkin Hamiltonian variational integrator. The extremization conditions and the discrete right Hamilton's equations can be expressed as
\begin{subequations}
\begin{align}
q_1 &= q_0 + \Delta t \sum_i b_i D_pH(Q^i,P^i), \label{GHVI1} \\
Q^i &= q_0 + \Delta t \sum_j a_{ij} D_pH(Q^j,P^j), \label{GHVI2} \\
p_1 &= p_0 - \Delta t \sum_i b_i D_qH(Q^i,P^i), \label{GHVI3} \\
P^i &= p_0 - \Delta t \sum_j \tilde{a}_{ij} D_qH(Q^i,P^i), \label{GHVI4}
\end{align}
\end{subequations}
where we interpret $a_{ij}$ as Runge--Kutta coefficients and $\tilde{a}_{ij} = (b_ib_j - b_ja_{ji})/b_i$ as the symplectic adjoint of the $a_{ij}$ coefficients. Thus, \eqref{GHVI1}-\eqref{GHVI4} can be viewed as a symplectic partitioned Runge--Kutta method.

We will consider such methods in four cases: adjoint systems corresponding to a base ODE or DAE, and whether or not the corresponding system is augmented. Note that in the DAE case, we will have to modify the above construction because the system is presymplectic. Furthermore, we will assume that all of the relevant configuration spaces are vector spaces. 

\textbf{Nonaugmented Adjoint ODE System.} The simplest case to consider is the nonaugmented adjoint ODE system \eqref{AdjointODEa}-\eqref{AdjointODEb}. Since the quadratic conservation law in Proposition \ref{ManifoldQuadraticInvariantProp}, 
$$\frac{d}{dt} \langle p,\delta q\rangle = 0,$$
arises from symplecticity, a structure-preserving discretization can be obtained by applying a symplectic integrator. This case is already discussed in \citet{Sa2016}, so we will only outline it briefly.

Applying the Galerkin Hamiltonian variational integrator \eqref{GHVI1}-\eqref{GHVI4} to the Hamiltonian for the adjoint ODE system, $H(q,p) = \langle p, f(q)\rangle, $
yields
\begin{subequations}
\begin{align}
q_1 &= q_0 + \Delta t \sum_i b_i f(Q^i), \label{SPRKAdjointODE1} \\
Q^i &= q_0 + \Delta t \sum_j a_{ij} f(Q^j), \label{SPRKAdjointODE2} \\
p_1 &= p_0 - \Delta t \sum_i b_i [Df(Q^i)]^*P^i, \label{SPRKAdjointODE3} \\
P^i &= p_0 - \Delta t \sum_j \tilde{a}_{ij} [Df(Q^j)]^*P^j. \label{SPRKAdjointODE4}
\end{align}
\end{subequations}
In the setting of adjoint sensitivity analysis of a terminal cost function, the appropriate boundary condition to prescribe on the momenta is $p_1 = \nabla_qC(q(t_f))$, as discussed in Section \ref{AdjointSensitivitySection}.

Since the above integrator is symplectic, we have the symplectic conservation law,
$$ dq_1 \wedge dp_1 = dq_0 \wedge dp_0, $$
when evaluated on discrete first variations of \eqref{SPRKAdjointODE1}-\eqref{SPRKAdjointODE4}. In this setting, a discrete first variation can be identified with solutions of the linearization of \eqref{SPRKAdjointODE1}-\eqref{SPRKAdjointODE4}. For the linearization of the equations in the position variables, \eqref{SPRKAdjointODE1}-\eqref{SPRKAdjointODE2}, we have
\begin{subequations}
\begin{align}
\delta q_1 &= \delta q_0 + \Delta t \sum_i b_i Df(Q^i)\delta Q^i, \label{DiscreteVariationalEquations1} \\
\delta Q^i &= \delta q_0 + \Delta t \sum_j a_{ij} Df(Q^j) \delta Q^j. \label{DiscreteVariationalEquations2}
\end{align}
\end{subequations}
As observed in \citet{Sa2016}, while we obtained this by linearizing the discrete equations, one could also obtain this by first linearizing \eqref{ODEvectorspace} and subsequently, applying the Runge--Kutta scheme to the linearization. For the linearization of the equations for the adjoint variables, \eqref{SPRKAdjointODE3}-\eqref{SPRKAdjointODE4}, observe that they are already linear in the adjoint variables, so we can identify the linearization with itself. Thus, we can choose for first variations vector fields $V$ as the first variation corresponding to the solution of the linearized position equation and $W$ as the first variation corresponding to the solution of the adjoint equation itself. With these choices, the above symplectic conservation law yields
$$ 0 = dq_1 \wedge dp_1(V,W)|_{(q_1,p_1)} - dq_0 \wedge dp_0 (V,W)|_{(q_0,p_0)} = \langle p_1, \delta q_1\rangle - \langle p_0, \delta q_0\rangle. $$
This is of course a discrete analogue of Proposition \ref{ManifoldQuadraticInvariantProp}. Note that one can derive the conservation law $\langle p_1,\delta q_1 \rangle = \langle p_0,\delta q_0\rangle$ directly by starting with the expression $\langle p_1,\delta q_1\rangle$ and substituting the discrete equations where appropriate. We will do this in the more general augmented case below. 

\textbf{Augmented Adjoint ODE System.} We now consider the case of the augmented adjoint ODE system \eqref{AugmentedAdjointCoordinate1}-\eqref{AugmentedAdjointCoordinate2}. In the continuous setting, we have from Proposition \ref{AugmentedQuadraticInvariantProp},
$$ \frac{d}{dt}\langle p,\delta q\rangle = -\langle dL, \delta q\rangle. $$
We would like to construct an integrator which admits a discrete analogue of this equation. To do this, we apply the Galerkin Hamiltonian variational integrator, equations \eqref{GHVI1}-\eqref{GHVI4}, to the augmented Hamiltonian $H_L(q,p) = \langle p,f(q)\rangle + L(q)$. This gives
\begin{subequations}
\begin{align}
q_1 &= q_0 + \Delta t \sum_i b_i f(Q^i), \label{SPRKAugmentedAdjointODE1} \\
Q^i &= q_0 + \Delta t \sum_j a_{ij} f(Q^j), \label{SPRKAugmentedAdjointODE2} \\
p_1 &= p_0 - \Delta t \sum_i b_i ([Df(Q^i)]^*P^i + dL(Q^i)) , \label{SPRKAugmentedAdjointODE3} \\
P^i &= p_0 - \Delta t \sum_j \tilde{a}_{ij} ([Df(Q^j)]^*P^j + dL(Q^j)). \label{SPRKAugmentedAdjointODE4}
\end{align}
\end{subequations}

We now prove a discrete analogue of Proposition \ref{AugmentedQuadraticInvariantProp}. To do this, we again consider the discrete variational equations for the position variables, \eqref{DiscreteVariationalEquations1}-\eqref{DiscreteVariationalEquations2}.
\begin{prop}\label{DiscreteAugmentedQuadraticInvariantProp}
With the above notation, the above integrator satisfies
\begin{equation}\label{DiscreteAugmentedQuadraticInvariant}
\langle p_1, \delta q_1\rangle = \langle p_0, \delta q_0\rangle - \Delta t \sum_i b_i \langle dL(Q^i), \delta Q^i\rangle.
\end{equation}
\begin{proof}
See Appendix \ref{QuadraticConservationProofs}.
\end{proof}
\end{prop}
\begin{remark}
To see that this is a discrete analogue of $\frac{d}{dt} \langle p,\delta q\rangle = -\langle dL,\delta q\rangle$, we write it in integral form as
$$ \langle p_1, \delta q_1\rangle = \langle p_0,\delta q_0\rangle - \int_0^{\Delta t} \langle dL(q),\delta q\rangle dt. $$
Then, applying the quadrature rule on $[0,\Delta t]$ given by quadrature weights $b_i\Delta t$ and quadrature nodes $c_i\Delta t$, the above integral is approximated by
$$ \int_0^{\Delta t} \langle dL(q),\delta q\rangle dt \approx \Delta t\sum_i b_i \langle dL(q(c_i\Delta t)), \delta q(c_i\Delta t) \rangle = \Delta t \sum_i b_i \langle dL(Q^i), \delta Q^i\rangle,$$
which yields equation \eqref{DiscreteAugmentedQuadraticInvariant}. The discrete analogue is natural in the sense that the quadrature rule for which the discrete equation \eqref{DiscreteAugmentedQuadraticInvariant} approximates the continuous equation is the same as the quadrature rule used to approximate the exact discrete generating function. This occurs more generally for such Hamiltonian variational integrators, as noted in \citet{TrLe2022} for the more general setting of multisymplectic Hamiltonian variational integrators.
\end{remark}

For adjoint sensitivity analysis of a running cost $\int L \, dt$, the appropriate boundary condition to prescribe on the momenta is $p_1 = 0$, as discussed in Section \ref{AdjointSensitivitySection}. With such a boundary condition, equation \eqref{DiscreteAugmentedQuadraticInvariant} reduces to
$$ \langle p_0, \delta q_0\rangle = \Delta t\sum_i b_i\langle dL(Q^i), \delta Q^i\rangle. $$
Thus, $p_0$ gives the discrete sensitivity, i.e., the change in the quadrature approximation of $\int L\, dt$ induced by a change in the initial condition along a discrete solution trajectory. 
One can compute this quantity directly via the direct method, where one needs to integrate the discrete variational equations for every desired search direction $\delta q_0$. On the other hand, by the above proposition, one can compute this quantity using the adjoint method: one integrates the adjoint equation with $p_1 = 0$ once to compute $p_0$ and subsequently, pair $p_0$ with any search direction $\delta q_0$ to obtain the sensitivity in that direction. By the above proposition, both methods give the same sensitivities. However, assuming the search space has dimension $n>1$, the adjoint method is more efficient since it only requires $\mathcal{O}(1)$ integrations and $\mathcal{O}(n)$ vector-vector products, whereas the direct method requires $\mathcal{O}(n)$ integrations and $\mathcal{O}(ns)$ vector-vector products where $s \geq 1$ is the number of Runge--Kutta stages, since, in the direct method, one has to compute $\langle dL(Q^i), \delta Q^i\rangle$ for each $i$ and for each choice of $\delta q_0$.

\textbf{Nonaugmented Adjoint DAE System.}
We will now construct discrete Hamiltonian variational integrators for the adjoint DAE system \eqref{DAE Adjoint System Coord 1}-\eqref{DAE Adjoint System Coord 4}, where we assume that the base DAE has index 1. To construct such a method, we have to modify the Galerkin Hamiltonian variational integrator \eqref{GHVI1}-\eqref{GHVI4}, so that it is applicable to the presymplectic adjoint DAE system. 

First, consider a general presymplectic system $i_X\Omega' = dH$. Note that, locally, any presymplectic system can be transformed to the canonical form~(see, \citet{CaIbGoRo1987}),
\begin{align*}
\dot{q} &= D_pH(q,p,r), \\
\dot{p} &= -D_qH(q,p,r), \\
0 &= D_rH(q,p,r),
\end{align*}
where, in these coordinates, $\Omega' = dq \wedge dp$, so that $\text{ker}(\Omega') = \text{span}\{\partial/\partial r\}.$ The action for this system is given by $\int_0^{\Delta t} (\langle p, \dot{q} \rangle - H(q,p,r) )dt$. We approximate this integral by quadrature, introduce internal stages for $q,p$ as before, and additionally introduce internal stages $R^i = r(c_ih)$. This gives the discrete generating function
$$ H_d^+(q_0,p_1) = \text{ext}\Big[ \langle p_1, q_1\rangle - \Delta t \sum_i b_i \left( \langle P^i,V^i\rangle - H(Q^i,P^i,R^i) \right) \Big], $$
where again $V$ is related to the internal stages of $Q$ by $Q^i = q_0 + \Delta t \sum_j a_{ij}V^j$ and the above expression is extremized over the internal stages $Q^i, P^i, R^i$. The discrete right Hamilton's equations are again given by
$$ q_1 = H_d^+(q_0,p_1),\ p_0 = H_d^+(q_0,p_1), $$
which we interpret as the evolution equations of the system. There are no evolution equations for $r$ due to the presymplectic structure and the absence of derivatives of $r$ in the action. This gives the integrator
\begin{subequations}
\begin{align}
q_1 &= q_0 + \Delta t \sum_i b_i D_pH(Q^i, P^i, R^i), \label{PresymplecticGHVI1}\\
Q^i &= q_0 + \Delta t \sum_j a_{ij} D_pH(Q^i,P^i,R^i), \label{PresymplecticGHVI2} \\
p_1 &= p_0 - \Delta t \sum_i b_i D_qH(Q^i, P^i, R^i), \label{PresymplecticGHVI3} \\
P^i &= p_0 - \Delta t \sum_j \tilde{a}_{ij} D_qH(Q^i, P^i, R^i), \label{PresymplecticGHVI4} \\
0 &= D_rH(Q^i,P^i,R^i), \label{PresymplecticGHVI5}
\end{align}
\end{subequations}
where \eqref{PresymplecticGHVI2}, \eqref{PresymplecticGHVI4}, \eqref{PresymplecticGHVI5} arise from extremizing with respect to $P^i, Q^i, R^i$, respectively, while \eqref{PresymplecticGHVI1} and \eqref{PresymplecticGHVI3} arise from the discrete right Hamilton's equations. This integrator is presymplectic, in the sense that
$$ dq_1 \wedge dp_1 = dq_0 \wedge dp_0, $$
when evaluated on discrete first variations. The proof is formally identical to the symplectic case. For this reason, we refer to \eqref{PresymplecticGHVI1}-\eqref{PresymplecticGHVI5} as a presymplectic Galerkin Hamiltonian variational integrator.

\begin{remark}
In general, the system \eqref{PresymplecticGHVI1}-\eqref{PresymplecticGHVI5} evolves on the primary constraint manifold given implicitly by the zero level set of $D_rH$, however, it may not evolve on the final constraint manifold. This is not an issue for us since we are dealing with adjoint DAE systems for index 1 DAEs, for which we know the primary constraint manifold and the final constraint manifold coincide. For the general case, one may need to additionally differentiate the constraint equation $D_rH = 0$ to obtain hidden constraints.

Thus, the method \eqref{PresymplecticGHVI1}-\eqref{PresymplecticGHVI5} is generally only applicable to index 1 presymplectic systems, unless we add in further hidden constraints. In order for the continuous presymplectic system to have index 1, it is sufficient that the Hessian of $H$ with respect to the algebraic variables, $D_r^2H$, is (pointwise) invertible on the primary constraint manifold. This is the case for the adjoint DAE system corresponding to an index 1 DAE.
\end{remark}

We now specialize to the adjoint DAE system \eqref{DAE Adjoint System Coord 1}-\eqref{DAE Adjoint System Coord 4}, corresponding to an index 1 DAE, which is already in the above canonical form with $r = (u,\lambda)$ and $H(q,u,p,\lambda) = \langle p,f(q,u)\rangle + \langle \lambda, \phi(q,u)\rangle$. Note that we reordered the argument of $H$, $(q,p,r) = (q,p,u,\lambda) \rightarrow (q,u,p,\lambda)$, in order to be consistent with the previous notation used throughout. We label the internal stages for the algebraic variables as $R^i = (U^i, \Lambda^i)$. Applying the presymplectic Galerkin Hamiltonian variational integrator to this particular system yields
\begin{subequations}
\begin{align} \label{SPRKAdjointDAE1}
q_1 &= q_0 + \Delta t \sum_i b_i f(Q^i,U^i), \\
Q^i &= q_0 + \Delta t \sum_j a_{ij} f(Q^j, U^j), \label{SPRKAdjointDAE2}\\
p_1 &= p_0 - \Delta t \sum_i b_i \left( [D_qf(Q^i,U^i)]^*P^i + [D_q\phi(Q^i,U^i)]^*\Lambda^i \right), \label{SPRKAdjointDAE3} \\
P^i &= p_0 - \Delta t \sum_j \tilde{a}_{ij} \left( [D_qf(Q^j,U^j)]^*P^j + [D_q\phi(Q^j,U^j)]^*\Lambda^j \right), \label{SPRKAdjointDAE4} \\
0 &= \phi(Q^i,U^i), \label{SPRKAdjointDAE5} \\
0 &= [D_uf(Q^i,U^i)]^*P^i + [D_u\phi(Q^i,U^i)]^*\Lambda^i, \label{SPRKAdjointDAE6}
\end{align}
\end{subequations}
where \eqref{SPRKAdjointDAE2}, \eqref{SPRKAdjointDAE4}, \eqref{SPRKAdjointDAE5}, \eqref{SPRKAdjointDAE6} arise from extremizing over $P^i, Q^i, \Lambda^i, U^i$, respectively, while \eqref{SPRKAdjointDAE1}, \eqref{SPRKAdjointDAE3} arise from the discrete right Hamilton's equations. 

\begin{remark}
In order for $q_1$ to appropriately satisfy the constraint, we should take the final quadrature point to be $c_s = 1$ (for an $s$-stage method), so that $\phi(q_1, U^s) = \phi(Q^s,U^s) = 0$. In this case, equation \eqref{SPRKAdjointDAE1} and equation \eqref{SPRKAdjointDAE2} with $i=s$ are redundant. Note that with the choice $c_s=1$, they are still consistent (i.e., are the same equation), since in the Galerkin construction, the coefficients $a_{ij}$ and $b_i$ are defined as
$$ a_{ij} = \int_0^{c_i} \phi_j(\tau)d\tau,\ b_i = \int_0^1 \phi_j(\tau)d\tau, $$
where $\phi_j$ are functions on $[0,1]$ which interpolate the nodes $c_j$ (see, \citet{LeZh2011}). Hence, $a_{sj} = b_j$, so that the two equations are consistent. However, we will write the system as above for conceptual clarity. Furthermore, even in the case where one does not take $c_s = 1$, the proposition that we prove below still holds, despite the possibility of constraint violations. 

A similar remark holds for the adjoint variable $p$ and the associated constraint \eqref{SPRKAdjointDAE6}, except we think of $p_0$ as the unknown, instead of $p_1$.
\end{remark}

Note that \eqref{SPRKAdjointDAE1}, \eqref{SPRKAdjointDAE2}, \eqref{SPRKAdjointDAE5} is a standard Runge--Kutta discretization of an index 1 DAE $\dot{q} = f(q,u)$, $0 = \phi(q,u)$, where again, usually $c_s = 1$. Associated with these equations are the variational equations given by their linearization,
\begin{subequations}
\begin{align}\label{DiscreteDAEVariationalEquations1}
\delta q_1 &= \delta q_0 + \Delta t \sum_i b_i(D_qf(Q^i,U^i)\delta Q^i + D_uf(Q^i,U^i)\delta U^i), \\
\delta Q^i &= \delta q_0 + \Delta t \sum_j a_{ij}(D_qf(Q^j,U^j)\delta Q^j + D_uf(Q^j,U^j)\delta U^j), \label{DiscreteDAEVariationalEquations2}\\
0 &= D_q\phi(Q^i,U^i)\delta Q^i + D_u\phi(Q^i,U^i) \delta U^i,\label{DiscreteDAEVariationalEquations3}
\end{align}
\end{subequations}
which is the Runge--Kutta discretization of the continuous variational equations \eqref{DAEVariationalEqn3} - \eqref{DAEVariationalEqn4}.
\begin{prop}\label{DiscreteDAEQuadraticInvariantProp}
With the above notation, the above integrator satisfies
$$ \langle p_1, \delta q_1\rangle = \langle p_0, \delta q_0\rangle. $$
\begin{proof}
See Appendix \ref{QuadraticConservationProofs}.
\end{proof}
\end{prop}
Thus, the above integrator admits a discrete analogue of Proposition \ref{Quadratic Invariant Adjoint DAE Prop} for the nonaugmented adjoint DAE system. By setting $p_1 = \nabla_q C(q(t_f))$, one can use this integrator to compute the sensitivity $p_0$ of a terminal cost function with respect to a perturbation in the initial condition. As discussed before, this only requires $\mathcal{O}(1)$ integrations instead of $\mathcal{O}(n)$ integrations via the direct method (for a dimension $n$ search space). Furthermore, the adjoint method requires only $\mathcal{O}(1)$ numerical solves of the constraints, while the direct method requires $\mathcal{O}(n)$ numerical solves.

\begin{remark}
Since we are assuming the DAE has index 1, it is always possible to prescribe an arbitrary initial condition $q_0$ (and $\delta q_0$) and terminal condition $p_1$, since the corresponding algebraic variables can always formally be solved for using the corresponding constraints. In practice, one generally has to solve the constraints to some tolerance, e.g., through an iterative scheme. If the constraints are only satisfied to a tolerance $\mathcal{O}(\epsilon)$, then the above proposition holds to $\mathcal{O}(s\epsilon)$, where $s$ is the number of Runge--Kutta stages.
\end{remark}

\begin{remark}
The above method \eqref{SPRKAdjointDAE1}-\eqref{SPRKAdjointDAE6} is presymplectic, since it is a special case of the more general presymplectic Galerkin Hamiltonian variational integrator \eqref{PresymplecticGHVI1}-\eqref{PresymplecticGHVI5}. Although we proved it directly, the above proposition could also have been proven from presymplecticity, with the appropriate choices of first variations. 
\end{remark}

\textbf{Augmented Adjoint DAE System.} Finally, we construct a discrete Hamiltonian variational integrator for the augmented adjoint DAE system \eqref{Augmented Adjoint DAE System Coord 1}-\eqref{Augmented Adjoint DAE System Coord 4} associated with an index 1 DAE. To do this, we apply the presymplectic Galerkin Hamiltonian variational integrator \eqref{PresymplecticGHVI1}-\eqref{PresymplecticGHVI5} with $r = (u,\lambda)$ and with Hamiltonian given by the augmented adjoint DAE Hamiltonian, 
$$ H_L(q,u,p,\lambda) = \langle p,f(q,u)\rangle + \langle \lambda,\phi(q,u)\rangle + L(q,u). $$
The presymplectic integrator is then
\begin{subequations}
\begin{align} \label{SPRKAugmentAdjointDAE1}
q_1 &= q_0 + \Delta t \sum_i b_i f(Q^i,U^i), \\
Q^i &= q_0 + \Delta t \sum_j a_{ij} f(Q^j, U^j), \label{SPRKAugmentAdjointDAE2}\\
p_1 &= p_0 - \Delta t \sum_i b_i \left( [D_qf(Q^i,U^i)]^*P^i + [D_q\phi(Q^i,U^i)]^*\Lambda^i + D_qL(Q^i,U^i) \right), \label{SPRKAugmentAdjointDAE3} \\
P^i &= p_0 - \Delta t \sum_j \tilde{a}_{ij} \left( [D_qf(Q^j,U^j)]^*P^j + [D_q\phi(Q^j,U^j)]^*\Lambda^j + D_qL(Q^i,U^i) \right), \label{SPRKAugmentAdjointDAE4} \\
0 &= \phi(Q^i,U^i), \label{SPRKAugmentAdjointDAE5} \\
0 &= [D_uf(Q^i,U^i)]^*P^i + [D_u\phi(Q^i,U^i)]^*\Lambda^i + D_uL(Q^i,U^i). \label{SPRKAugmentAdjointDAE6}
\end{align}
\end{subequations}
The associated variational equations are again \eqref{DiscreteDAEVariationalEquations1}-\eqref{DiscreteDAEVariationalEquations3}.
Remarks analogous to the nonaugmented case regarding setting the quadrature node $c_s=1$ and solvability of these systems under the index 1 assumption can be made.

\begin{prop}\label{DiscreteAugmentedAdjointDAEQuadraticProp}
With the above notation, the above integrator satisfies
$$ \langle p_1, \delta q_1\rangle = \langle p_0, \delta q_0\rangle - \Delta t \sum_i b_i \langle dL(Q^i,U^i), (\delta Q^i, \delta U^i)\rangle. $$
\begin{proof}
See Appendix \ref{QuadraticConservationProofs}.
\end{proof}
\end{prop}

\begin{remark}
Analogous to the remark in the augmented adjoint ODE case, the above proposition is a discrete analogue of Proposition \ref{Quadratic Invariant Augmented DAE Adjoint Prop}, in integral form,
$$ \langle p_1, \delta q_1\rangle - \langle p_0,\delta q_0\rangle = - \int_0^{\Delta t}\langle dL(q,u), (\delta q, \delta u)\rangle dt. $$
The discrete analogue is natural in the sense that it is just quadrature applied to the right hand side of this equation, with the same quadrature rule used to discretize the generating function.
\end{remark}

\begin{remark}
As with the augmented adjoint ODE case, the above proposition allows one to compute numerical sensitivities of a running cost function by solving for $p_0$ with $p_1 = 0$, which is more efficient than the direct method.
\end{remark}

To summarize, we have utilized Galerkin Hamiltonian variational integrators to construct methods which admit natural discrete analogues of the various propositions used for sensitivity analysis. We summarize the results below.

\begin{center}
\begin{tabular}{ |c|c|c| } 
 \hline
  & Terminal Cost & Running Cost \\ 
  \hline
 ODE & $\langle p_1, \delta q_1\rangle = \langle p_0, \delta q_0\rangle$ & $\langle p_1, \delta q_1\rangle = \langle p_0, \delta q_0\rangle - \Delta t \sum_i b_i \langle dL(Q^i), \delta Q^i\rangle$ \\ 
 DAE & $\langle p_1, \delta q_1\rangle = \langle p_0, \delta q_0\rangle$ & $\langle p_1,\delta q_1\rangle = \langle p_0, \delta q_0\rangle - \Delta t\sum_i b_i \langle dL(Q^i,U^i), (\delta Q^i,\delta U^i)\rangle$ \\ 
 \hline
\end{tabular}
\end{center}

\subsubsection{Naturality of the Adjoint DAE System Discretization}\label{NaturalityDiscretizationSection}
To conclude our discussion of discretizing adjoint systems, we prove a discrete extension of the fact that, for an index 1 DAE, the process of index reduction and forming the adjoint system commute, as discussed in Section \ref{DAEIndexPCASection}. Namely, we will show that, starting from an index 1 DAE \eqref{DAEa}-\eqref{DAEb}, the processes of reduction, forming the adjoint system, and discretization all commute, for particular choices of these processes which we will define and choose below. This can be summarized in the following commutative diagram.

\[\small\begin{tikzcd}[column sep=-4ex,row sep=10ex]
	\txt{Index 1 DAE} && \txt{ODE} \\
	& \txt{Discrete DAE} && \txt{Discrete ODE} \\
	\txt{Presymplectic Adjoint\\ DAE System} && \txt{Symplectic Adjoint\\ ODE System} \\
	& \txt{Presymplectic Galerkin\\ Hamiltonian Variational \\ Integrator} && \txt{Symplectic Galerkin \\ Hamiltonian Variational \\ Integrator}
	\arrow["{\text{Reduce}}", from=1-1, to=1-3]
	\arrow["{\text{Reduce}}"{pos=0.25}, from=3-1, to=3-3]
	\arrow["{\text{Adjoint}}"'{pos=0.6}, from=1-1, to=3-1]
	\arrow["{\text{Adjoint}}"{pos=0.6}, from=1-3, to=3-3]
	\arrow["{\text{Reduce}}"{pos=0.3}, from=2-2, to=2-4,crossing over]
	\arrow["{\text{Reduce}}"', from=4-2, to=4-4]
	\arrow["{\text{Adjoint}}"{pos=0.6}, from=2-2, to=4-2,crossing over]
	\arrow["{\text{Adjoint}}"{pos=0.6}, from=2-4, to=4-4]
	\arrow["{\text{Discretize}}"', from=1-1, to=2-2]
	\arrow["{\text{Discretize}}"{pos=0.6}, from=1-3, to=2-4]
	\arrow["{\text{Discretize}}"', from=3-1, to=4-2]
	\arrow["{\text{Discretize}}"', from=3-3, to=4-4]
\end{tikzcd}\]

In the above diagram, we will use the convention that the ``Discretize" arrows point forward, the ``Adjoint" arrows point downward, and the ``Reduce" arrows point to the right. For the ``Discretize" arrows on the top face, we take the discretization to be a Runge--Kutta discretization (of a DAE on the left and of an ODE on the right, with the same Runge--Kutta coefficients in both cases). For the ``Discretize" arrows on the bottom face, we take the discretization to be the symplectic partitioned Runge--Kutta discretization induced by the discretization of the base DAE or ODE, i.e., the momenta expansion coefficients $\tilde{a}_{ij}$ are the symplectic adjoint of the coefficients $a_{ij}$ used on the top face. We have already defined the ``Adjoint" arrows on the back face, as discussed in Section \ref{AdjointSystems Main Section}. For the ``Adjoint" arrows on the front face, we define them as forming the discrete adjoint system corresponding to a discrete (and generally nonlinear) system of equations and we will review this notion where needed in the proof. We have already defined the ``Reduce" arrows on the back face, as discussed in Section \ref{DAEIndexPCASection}. For the ``Reduce" arrows on the front face, we define this as solving for the discrete algebraic variables in terms of the discrete kinematic variables through the discrete constraint equations. With these choices, the above diagram commutes, as we will show. To prove this, it suffices to prove that the diagrams on each of the six faces commutes. To keep the exposition concise, we provide the proof in Appendix \ref{AppendixNaturalityProof} and move on to discuss the implications of this result.

The previous discussion shows that the presymplectic Galerkin Hamiltonian variational integrator construction is natural for discretizing adjoint (index 1) DAE systems, in the sense that the integrator is equivalent to the integrator produced from applying a symplectic Galerkin Hamiltonian variational integrator to the underlying nondegenerate Hamiltonian system. Of course, in practice, one cannot generally determine the function $u = u(q)$ needed to reduce the DAE to an ODE. Therefore, one generally works with the presymplectic Galerkin Hamiltonian variational integrator instead, where one iteratively solves the constraint equations. However, although reduction then symplectic integration is often impractical, one can utilize this naturality to derive properties of the presymplectic integrator. For example, we will use this naturality to prove a variational error analysis result. 

The basic idea for the variational error analysis result goes as follows: one utilizes the naturality to relate the presymplectic variational integrator to a symplectic variational integrator of the underlying nondegenerate Hamiltonian system and subsequently, applies the variational error analysis result in the symplectic case (\citet{ScLe2017}). Recall the discrete generating function for the previously constructed presymplectic variational integrator,
$$ H_d^+(q_0,p_1; \Delta t) = \text{ext}\Big[ \langle p_1, q_1\rangle - \Delta t \sum_i b_i \left( \langle P^i,V^i\rangle - H(Q^i,U^i,P^i,\Lambda^i) \right) \Big], $$
where we have now explicitly included the timestep dependence in $H_d^+$ and $H$ is the Hamiltonian for the adjoint DAE system (augmented or nonaugmented), corresponding to an index 1 DAE. 
\begin{prop}
Suppose the discrete generating function $H_d^+(q_0,p_1; \Delta t)$ for the presymplectic variational integrator approximates the exact discrete   generating function $H_d^{+,E}(q_0,p_1; \Delta t)$ to order $r$, i.e.,
$$ H_d^+(q_0,p_1; \Delta t) = H_d^{+,E}(q_0,p_1;\Delta t) + \mathcal{O}(\Delta t^{r+1}), $$
and the Hamiltonian $H$ is continuously differentiable, then the Type II map $(q_0,p_1) \mapsto (q_1, p_0)$ and the evolution map $(q_0,p_0) \mapsto (q_1, p_1)$ are order-$r$ accurate. 
\begin{proof}
The proof follows from two simple steps. First, observe that the discrete generating function $H_d^+(q_0,p_1; \Delta t)$ for the presymplectic integrator is also the discrete generating function for the symplectic integrator for the underlying nondegenerate Hamiltonian system. This follows since in the definition of $H_d^+$, one extremizes over the algebraic variables $U^i,\Lambda^i$ which enforces the constraints and hence, determines $U^i,\Lambda^i$ as functions of the kinematic variables $Q^i,P^i$. Thus, the discrete (or continuous) Type II map determined by $H_d^+$ (or $H_d^{+,E}$, respectively), $(q_0,p_1) \mapsto (q_1,p_0)$, is the same as the Type II map for the underlying nondegenerate Hamiltonian system, which is just another consequence of the aforementioned naturality. One then applies the variational error analysis result in \citet{ScLe2017}.
\end{proof}
\end{prop}
\begin{remark}
Another way to view this result is that the order of an implicit (partitioned) Runge--Kutta scheme for index 1 DAEs is the same as the order of an implicit (partitioned) Runge--Kutta scheme for ODEs (\citet{Ro1989}), since the aforementioned discretization generates a partitioned Runge--Kutta scheme. To be complete, we should determine the order for the full presymplectic flow, i.e., including also the algebraic variables. As discussed in \citet{Ro1989}, as long as $a_{si} = b_i$ for each $i$, which, as we have discussed, is a natural choice and holds as long as $c_s=1$, there is no order reduction arising from the algebraic variables. Thus, with this assumption, the presymplectic variational integrator in the previous proposition approximates the presymplectic flow, in both the kinematic and algebraic variables, to order $r$.
\end{remark}
\begin{remark}
In the above proposition, we considered both the Type II map $(q_0, p_1) \mapsto (q_1, p_0)$ and the evolution map $(q_0,p_0) \mapsto (q_1,p_1)$.  The latter is of course the traditional way to view the map corresponding to a numerical method, but the former is the form of the map used in adjoint sensitivity analysis.
\end{remark}

Furthermore, in light of this naturality, we can view Propositions \ref{DiscreteDAEQuadraticInvariantProp} and \ref{DiscreteAugmentedAdjointDAEQuadraticProp} as following from the analogous propositions for symplectic Galerkin Hamiltonian variational integrators, applied to the underlying nondegenerate Hamiltonian system.

\subsubsection{Numerical Example}
For our numerical example, we consider the planar pendulum. Although one can formulate this system as an ODE in the angular variable $\theta$, we instead work with this system in Cartesian coordinates $xy$ where this system is formulated as a DAE, as an academic example of the theory presented in this paper. We will derive the adjoint DAE system associated to the planar pendulum DAE, and subsequently, perform a numerical test demonstrating the presymplecticity of a presymplectic Galerkin Hamiltonian variational integrator applied to this system.

Consider a pendulum of mass $m > 0$ and length $L > 0$ confined to the $xy$ plane, where gravity acts in the vertical $y$ direction, with acceleration $-g < 0$. This is described by the system
\begin{align*}
m \ddot{x} &= \rho x, \\
m \ddot{y} &= \rho y - mg, \\
x^2 + y^2 &= L^2.
\end{align*}
This system can be derived from the Lagrangian
$$ L = \frac{1}{2}m (\dot{x}^2 + \dot{y}^2) - mg(y - L) + \frac{1}{2}\rho(x^2+y^2-L^2), $$
where the first term is the kinetic energy, the second term is (minus) the potential energy, and the third term enforces the constraint $x^2+y^2 = L^2$ where $\rho$ is interpreted as a Lagrange multiplier.

If we restrict to the region $y<0$, the above system can be expressed as a semi-explicit index 1 DAE of the form
\begin{subequations}
\begin{align}
\dot{x} &= v_x, \\
\dot{v}_x &= \rho x / m, \\
0 &= x^2+y^2 - L^2, \\
0 &= v_x x + v_y y, \\
0 &= m(v_x^2+v_y^2) - mgy + L^2\rho.
\end{align}
\end{subequations}

In terms of the notation of Section \ref{AdjointDAESection}, we have $(x,v_x) \in M_d = (-1,1) \times \mathbb{R}$ and $(y,v_y,\rho) \in M_a = \mathbb{R}_{-} \times \mathbb{R} \times \mathbb{R}$. Letting $q = (x,v_x)$ denote the coordinates for the dynamical variables and $u = (y,v_y,\rho)$ denote the coordinates for the algebraic variables, this system can be expressed in the form \eqref{DAEa}-\eqref{DAEb}, where
\begin{align*}
f(q,u) &= \begin{pmatrix} v_x \\ \rho x/m \end{pmatrix}, \\
\phi(q,u) &= \begin{pmatrix} x^2+y^2-L^2 \\ v_x x + v_y y \\ m(v_x^2+v_y^2) - mgy + L^2\rho \end{pmatrix}.
\end{align*}
We regard $\phi$ as a section of the constraint bundle $\Phi$ given by the trivial vector bundle $(M_d \times M_a) \times \mathbb{R}^3 \rightarrow M_d \times M_a$. Coordinatize $\overline{T^*M_d}$ by $(q,u,p)$ where $p = (p_x, p_{v_x})$ are the momenta dual to $q = (x, v_x)$ and coordinatize $\Phi^*$ by $(q,u,\lambda)$ where $\lambda = (\lambda_1, \lambda_2, \lambda_3)$ are the coordinates of the fibers dual to the constraint bundle fibers. The Hamiltonian $H: \overline{T^*M_d} \oplus \Phi^* \rightarrow \mathbb{R}$ is then given by
\begin{align*}
H(q,u,p,\lambda) &= \langle p, f(q,u)\rangle + \langle \lambda, \phi(q,u)\rangle \\
&= ( p_x \ p_{v_x} ) \begin{pmatrix} v_x \\ \rho x \end{pmatrix} + (\lambda_1 \ \lambda_2\ \lambda_3 ) \begin{pmatrix} x^2+y^2-L^2 \\ v_x x + v_y y \\ m(v_x^2+v_y^2) - mgy + L^2\rho \end{pmatrix}.
\end{align*}
The presymplectic form $\Omega_0$ on $\overline{T^*M_d} \oplus \Phi^*$ is given by
$$ \Omega_0 = dq \wedge dp = dx \wedge dp_x + dv_x \wedge dp_{v_x}. $$
To obtain an expression for the adjoint DAE system \eqref{DAE Adjoint System Coord 1}-\eqref{DAE Adjoint System Coord 4}, we compute the derivative matrices of $f$ and $\phi$.
\begin{align*}
D_qf(q,u) &= \begin{pmatrix} 0 & 1 \\ \rho/m & 0 \end{pmatrix}, \\
D_uf(q,u) &= \begin{pmatrix} 0 & 0 & 0 \\ 0 & 0 & x/m \end{pmatrix}, \\
D_q\phi(q,u) &= \begin{pmatrix} 2x & 0 \\ v_x & x \\ 0 & 2mv_x \end{pmatrix}, \\
D_u\phi(q,u) &= \begin{pmatrix} 2y & 0 & 0 \\ v_y & y & 0 \\ -mg & 2mv_y & L^2 \end{pmatrix}.
\end{align*}
Note that $\det(D_u\phi(q,u)) = 2L^2y^2 \neq 0$ for $(q,u) \in M_d \times M_a$ and hence, the system is an index 1 DAE as previously claimed.

The adjoint DAE system \eqref{DAE Adjoint System Coord 1}-\eqref{DAE Adjoint System Coord 4} for the planar pendulum is then given by
\begin{subequations}
\begin{align}
\frac{d}{dt} \begin{pmatrix} x \\ v_x \end{pmatrix} &= \begin{pmatrix} v_x \\ \rho x/m \end{pmatrix}, \label{AdjointDAEpendulum1} \\
\frac{d}{dt} \begin{pmatrix} p_x \\ p_{v_x} \end{pmatrix} &= - \begin{pmatrix} 0 & 1 \\ \rho/m & 0 \end{pmatrix}^T \begin{pmatrix} p_x \\ p_{v_x} \end{pmatrix} - \begin{pmatrix} 2x & 0 \\ v_x & x \\ 0 & 2mv_x \end{pmatrix}^T \begin{pmatrix} \lambda_1 \\ \lambda_2 \\ \lambda_3 \end{pmatrix}, \\
0 &= \begin{pmatrix} x^2+y^2-L^2 \\ v_x x + v_y y \\ m(v_x^2+v_y^2) - mgy + L^2\rho \end{pmatrix}, \\
0 &= \begin{pmatrix} 0 & 0 & 0 \\ 0 & 0 & x/m \end{pmatrix}^T \begin{pmatrix} p_x \\ p_{v_x} \end{pmatrix} + \begin{pmatrix} 2y & 0 & 0 \\ v_y & y & 0 \\ -mg & 2mv_y & L^2 \end{pmatrix}^T \begin{pmatrix} \lambda_1 \\ \lambda_2 \\ \lambda_3 \end{pmatrix}. \label{AdjointDAEpendulum4}
\end{align}
\end{subequations}

We will apply a presymplectic Galerkin Hamiltonian variational integrator \eqref{SPRKAdjointDAE1}-\eqref{SPRKAdjointDAE6} to the above system. We choose a first-order Runge--Kutta method, with Runge--Kutta coefficients $a=1, b=1, c=1$ and hence, $\tilde{a} = 0$. Thus, the internal stages for the position and momenta are given by $Q = q_1$ and $P = p_0$. With these choices, the presymplectic Galerkin Hamiltonian variational integrator can be expressed as
\begin{align*} 
q_1 &= q_0 + \Delta t f(q_1,U), \\
p_1 &= p_0 - \Delta t \left( [D_qf(q_1,U)]^* p_0 + [D_q\phi(q_1,U)]^*\Lambda \right),  \\
0 &= \phi(q_1,U),  \\
0 &= [D_uf(q_1,U)]^*p_0 + [D_u\phi(q_1,U)]^*\Lambda.
\end{align*}
For our example, we set $m=g=L=1$. Letting $U = (Y, V_y, \mathcal{P})$ and $\Lambda = (\Lambda_1, \Lambda_2, \Lambda_3)$ denote the internal stages corresponding to $u = (y, v_y, \rho)$ and $\lambda = (\lambda_1, \lambda_2, \lambda_3)$, respectively, the above integrator applied to the adjoint DAE system for the planar pendulum \eqref{AdjointDAEpendulum1}-\eqref{AdjointDAEpendulum4}, with $m=g=L=1$, can be expressed as
\begin{align*}
\begin{pmatrix} x_1 \\ (v_x)_1 \end{pmatrix} &= \begin{pmatrix} x_0 \\ (v_x)_0 \end{pmatrix} + \Delta t \begin{pmatrix} (v_x)_1 \\ \mathcal{P} x_1 \end{pmatrix}, \\
\begin{pmatrix} (p_x)_1 \\ (p_{v_x})_1 \end{pmatrix} &= \begin{pmatrix} (p_x)_0 \\ (p_{v_x})_0 \end{pmatrix} - \Delta t \left( \begin{pmatrix} 0 & 1 \\ \mathcal{P} & 0 \end{pmatrix}^T \begin{pmatrix} (p_x)_0 \\ (p_{v_x})_0 \end{pmatrix} + \begin{pmatrix} 2x_1 & 0 \\ (v_x)_1 & x_1 \\ 0 & 2(v_x)_1 \end{pmatrix}^T \begin{pmatrix} \Lambda_1 \\ \Lambda_2 \\ \Lambda_3 \end{pmatrix} \right), \\
0 &= \begin{pmatrix} x_1^2+Y^2-1 \\ (v_x)_1 x_1 + V_y Y \\ (v_x)_1^2+V_y^2 - Y + \mathcal{P} \end{pmatrix}, \\
0 &= \begin{pmatrix} 0 & 0 & 0 \\ 0 & 0 & x_1 \end{pmatrix}^T \begin{pmatrix} (p_x)_0 \\ (p_{v_x})_0 \end{pmatrix} + \begin{pmatrix} 2Y & 0 & 0 \\ V_y & Y & 0 \\ -1 & 2V_y & 1 \end{pmatrix}^T \begin{pmatrix} \Lambda_1 \\ \Lambda_2 \\ \Lambda_3 \end{pmatrix}.
\end{align*}
We refer to this method as PGHVI--1. We will compare this to the first-order method where the Runge--Kutta coefficients are the same for both $q$ and $p$, i.e., $a = 1 = \tilde{a}$. This method, which we refer to as BE--1, is given by applying the backward Euler method in both the $q$ and $p$ variables, i.e.,
\begin{align*}
\begin{pmatrix} x_1 \\ (v_x)_1 \end{pmatrix} &= \begin{pmatrix} x_0 \\ (v_x)_0 \end{pmatrix} + \Delta t \begin{pmatrix} (v_x)_1 \\ \mathcal{P} x_1 \end{pmatrix}, \\
\begin{pmatrix} (p_x)_1 \\ (p_{v_x})_1 \end{pmatrix} &= \begin{pmatrix} (p_x)_0 \\ (p_{v_x})_0 \end{pmatrix} - \Delta t \left( \begin{pmatrix} 0 & 1 \\ \mathcal{P} & 0 \end{pmatrix}^T \begin{pmatrix} (p_x)_1 \\ (p_{v_x})_1 \end{pmatrix} + \begin{pmatrix} 2x_1 & 0 \\ (v_x)_1 & x_1 \\ 0 & 2(v_x)_1 \end{pmatrix}^T \begin{pmatrix} \Lambda_1 \\ \Lambda_2 \\ \Lambda_3 \end{pmatrix} \right), \\
0 &= \begin{pmatrix} x_1^2+Y^2-1 \\ (v_x)_1 x_1 + V_y Y \\ (v_x)_1^2+V_y^2 - Y + \mathcal{P} \end{pmatrix}, \\
0 &= \begin{pmatrix} 0 & 0 & 0 \\ 0 & 0 & x_1 \end{pmatrix}^T \begin{pmatrix} (p_x)_1 \\ (p_{v_x})_1 \end{pmatrix} + \begin{pmatrix} 2Y & 0 & 0 \\ V_y & Y & 0 \\ -1 & 2V_y & 1 \end{pmatrix}^T \begin{pmatrix} \Lambda_1 \\ \Lambda_2 \\ \Lambda_3 \end{pmatrix}.
\end{align*}

For our numerical test, we will qualitatively compare the preservation of the presymplectic form $\Omega_0 = dx \wedge dp_x + dv_x \wedge dp_{v_x}$ between the two methods. Since Type II boundary conditions arise in adjoint sensitivity analysis, we place Type II boundary conditions, i.e., by specifying $q_0 = (x_0, (v_x)_0)$ and $p_1 = ( (p_x)_1, (p_{v_x})_1 )$, and subsequently, numerically solve the resulting system for $q_1, p_0, U, \Lambda$. We use various nearby values for the initial position $q_0 = (x_0, (v_x)_0)$ and various nearby values for the final momenta $p_1 = ( (p_x)_1, (p_{v_x})_1 )$. For a presymplectic integrator applied to a presymplectic system with presymplectic form $dx \wedge dp_x + dv_x \wedge dp_{v_x}$, we expect that the area occupied by the distribution of points $(x_0, (p_x)_0)$ is the same as the area occupied by the distribution of points $(x_1, (p_x)_1)$; similarly, we expect that the area occupied by the distribution of points $((v_x)_0, (p_{v_x})_0)$ is the same as the area occupied by the distribution of points $((v_x)_1, (p_{v_x})_1)$. Since we choose to only solve the system for one timestep, we take a large timestep to highlight the difference between the two methods, $\Delta t = 2$, which corresponds to roughly one-third of the period of the pendulum.

\begin{figure}[H]
\begin{center}
\includegraphics[width=130mm]{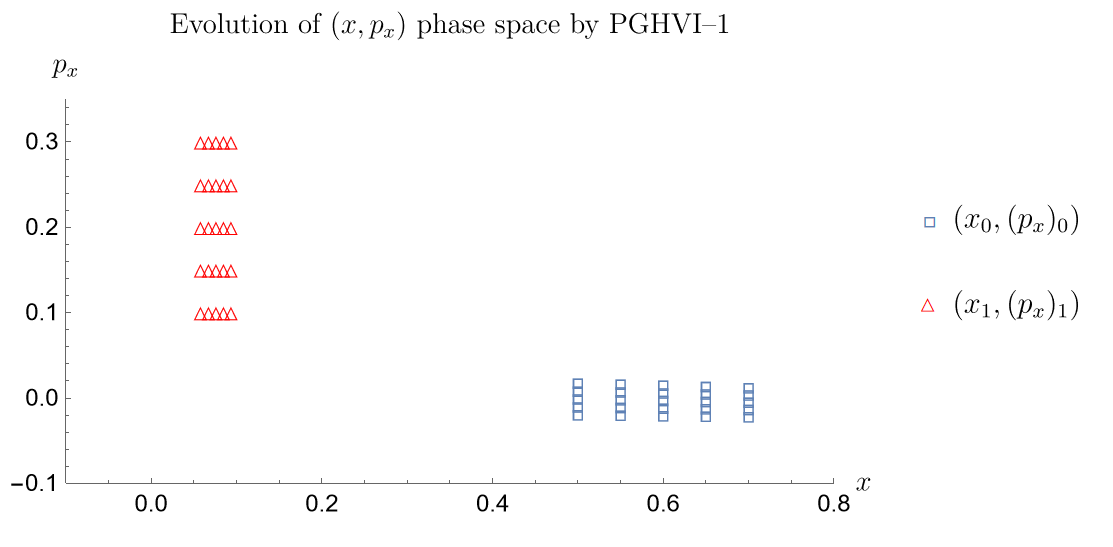}
\caption{$(x,p_x)$ phase space cross-section of PGHVI--1 applied to a distribution of initial conditions $q_0$ and final momenta $p_1$}
\label{pghvi_xp}
\end{center}
\end{figure}

\begin{figure}[H]
\begin{center}
\includegraphics[width=130mm]{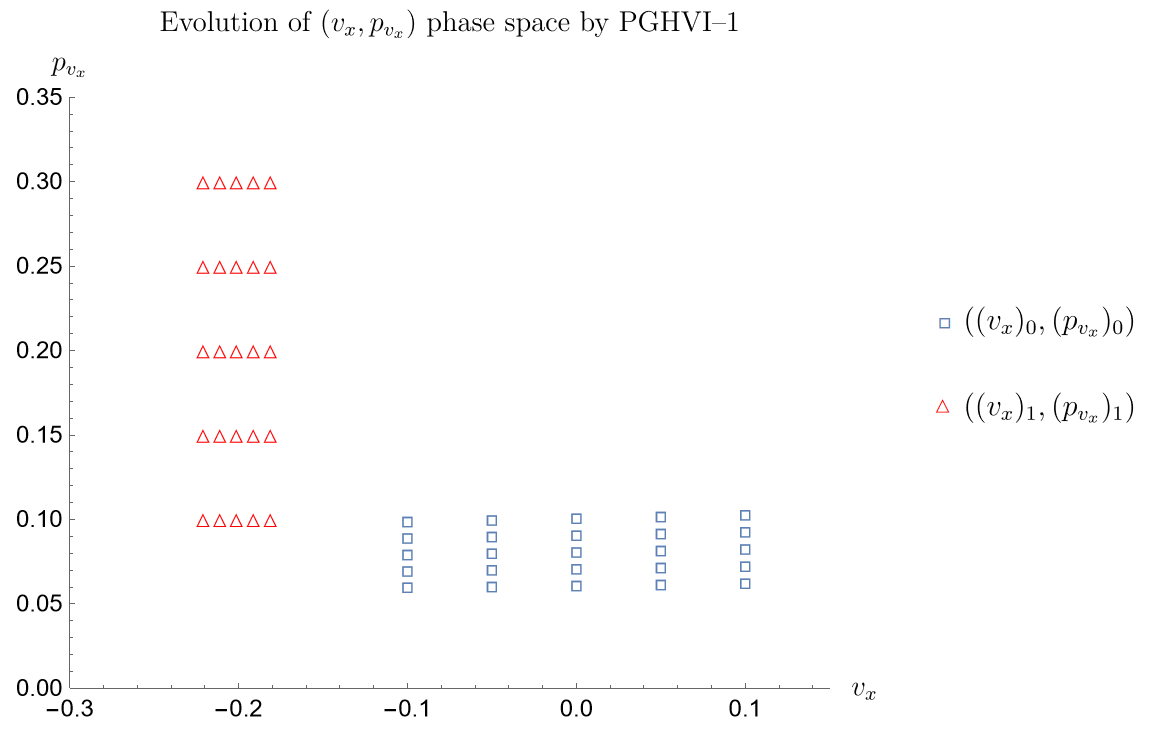}
\caption{$(v_x,p_{v_x})$ phase space cross-section of PGHVI--1 applied to a distribution of initial conditions $q_0$ and final momenta $p_1$}
\label{pghvi_vxpvx}
\end{center}
\end{figure}

\begin{figure}[H]
\begin{center}
\includegraphics[width=130mm]{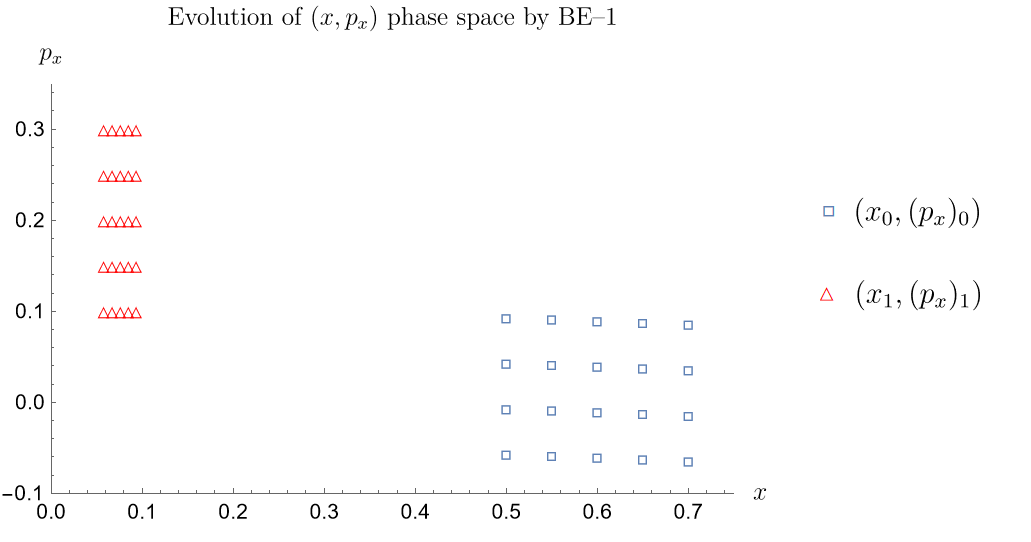}
\caption{$(x,p_x)$ phase space cross-section of BE--1 applied to a distribution of initial conditions $q_0$ and final momenta $p_1$}
\label{be_xp}
\end{center}
\end{figure}

\begin{figure}[H]
\begin{center}
\includegraphics[width=120mm]{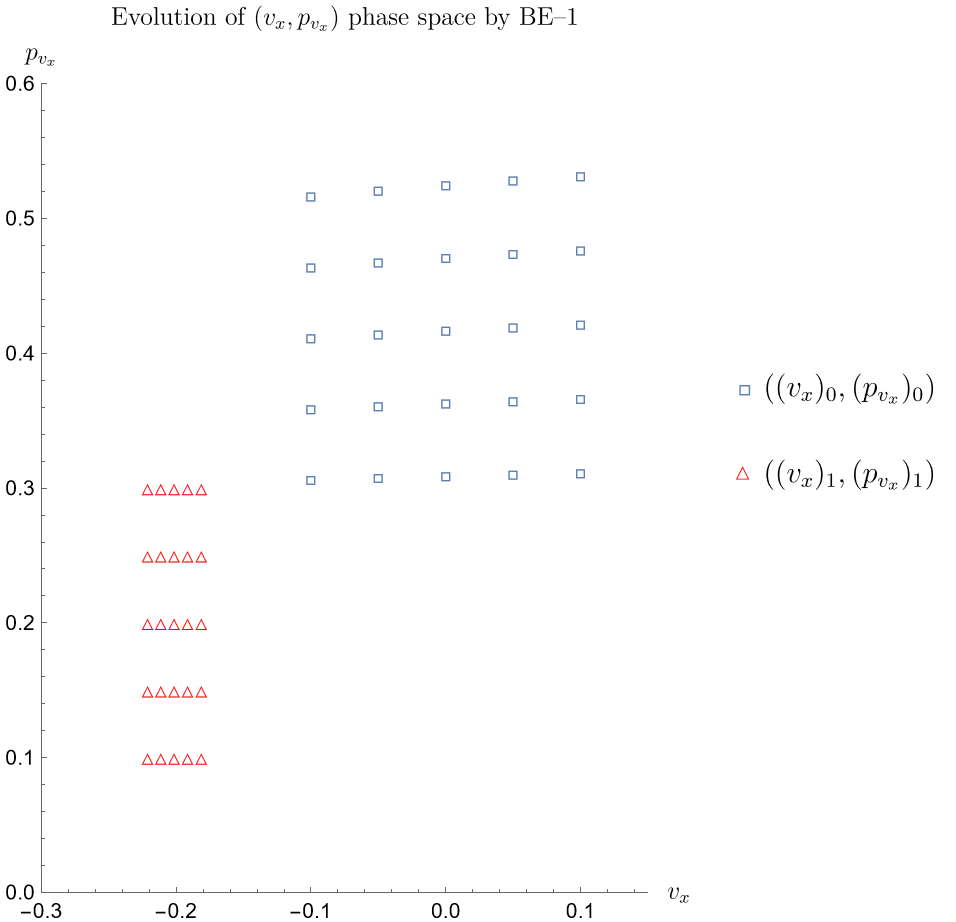}
\caption{$(v_x,p_{v_x})$ phase space cross-section of BE--1 applied to a distribution of initial conditions $q_0$ and final momenta $p_1$}
\label{be_vxpvx}
\end{center}
\end{figure}

Note that, with Type II boundary conditions, both methods give a map $(q_0,p_1) \mapsto (q_1,p_0)$ which implicitly determines an evolution map $(q_0,p_0) \mapsto (q_1,p_1)$; below, we plot the phase space cross-sections of these implicit evolution maps. The evolution of the $(x,p_x)$ and $(v_x, p_{v_x})$ distributions by PGHVI--1 is shown in Figure \ref{pghvi_xp} and Figure \ref{pghvi_vxpvx}, respectively. The evolution of the $(x,p_x)$ and $(v_x, p_{v_x})$ distributions by BE--1 is shown in Figure \ref{be_xp} and Figure \ref{be_vxpvx}, respectively. As can be qualitatively seen from Figures \ref{pghvi_xp}, \ref{pghvi_vxpvx}, \ref{be_xp}, \ref{be_vxpvx}, the PGHVI--1 method preserves the phase space area in both the $(x,p_x)$ and $(v_x, p_{v_x})$ cross-sections, whereas the BE--1 method does not.


\subsection{Optimal Control of DAE Systems}\label{OCPSection}
In this section, we derive the optimality conditions for an optimal control problem (OCP) subject to a semi-explicit DAE constraint. It is known that the optimality conditions can be described as a presymplectic system on the generalized phase space bundle (\citet{DeIb2003}, \citet{EcMaMuRo2003}). For a discussion of the presymplectic geometry of optimal control systems and in particular, symmetries of such systems, see \citet{dLCdDM2004}. We will subsequently consider a variational discretization of such OCPs and discuss the naturality of such discretizations.

Consider the following optimal control problem in Bolza form, subject to a DAE constraint, which we refer to as (OCP-DAE),
\begin{align*}
&\text{min } C(q(t_f)) + \int_0^{t_f}L(q,u)dt \\
&\text{subject to} \\
& \qquad \dot{q} = f(q,u), \\
& \qquad 0 = \phi(q,u), \\
& \qquad q_0 = q(0), \\
& \qquad 0 = \phi_f(q(t_f)),
\end{align*}
where the DAE system $\dot{q} = f(q,u)$, $0 = \phi(q,u)$ is over $M_d \times M_a$ as described in Section \ref{AdjointDAESection}, $C: M_d \rightarrow \mathbb{R}$ is the terminal cost, $L: M_d \times M_a \rightarrow \mathbb{R}$ is the running cost, the initial condition $q(0) = q_0$ is prescribed, and for generality, a terminal constraint $\phi_f(q(t_f)) = 0$ is also imposed, where $\phi_f$ is a map from $M_d$ into some vector space $V$.

We assume a local optimum to (OCP-DAE). We then adjoin the constraints to $J$ using adjoint variables, which gives the adjoined functional
$$ \mathcal{J} = C(q(t_f)) + \langle \lambda_f, \phi_f(q(t_f))\rangle + \int_{0}^{t_f} \left[ L(q,u) + \langle p, f(q,u) - \dot{q}\rangle + \langle \lambda, \phi(q,u)\rangle \right]dt. $$
The optimality conditions are given by the condition that $\mathcal{J}$ is stationary about the local optimum, $\delta \mathcal{J} = 0$ (\citet{Bi2010}). For simplicity in the notation, we will use matrix derivative instead of indices. Note also that we will implicitly leave out the variation of the adjoint variables, since those terms pair with the DAE constraints, which vanish at the local optimum. The optimality condition $\delta \mathcal{J} = 0$ is then
\begin{align*}
0 &= \delta \mathcal{J} = \langle \nabla_q C(q(t_f)), \delta q(t_f) \rangle + \langle \lambda_f, D_q\phi_f(q(t_f)) \delta q(t_f) \rangle \\ 
&\qquad \qquad \quad + \int_0^{t_f} \Big[ \langle \nabla_q L(q,u), \delta q\rangle + \langle \nabla_u L(q,u), \delta u \rangle + \langle p, D_qf(q,u) \delta q \rangle + \langle p, D_uf(q,u)\delta u\rangle  \\
& \qquad \qquad \qquad \qquad - \langle p, \frac{d}{dt} \delta q\rangle + \langle \lambda, D_q\phi(q,u)\delta q \rangle + \langle \lambda, D_u\phi(q,u)\delta u \rangle \Big]dt \\
&= \langle \nabla_q C(q(t_f)) + [D_q\phi_f(q(t_f))]^* \lambda_f - p(t_f), \delta q(t_f)\rangle \\
& \qquad + \int_0^{t_f} \Big[ \langle \nabla_qL(q,u) + [D_qf(q,u)]^*p + \dot{p} + [D_q\phi(q,u)]^*\lambda, \delta q \rangle  \\
& \qquad \quad\quad \quad + \langle \nabla_u L(q,u) + [D_uf(q,u)]^*p + [D_u\phi(q,u)]^*\lambda, \delta u\rangle \Big] dt,
\end{align*}
where we integrated by parts on the term $\langle p, \frac{d}{dt} \delta q\rangle$ and used $\delta q(0) = 0$ since the initial condition is fixed. Enforcing stationarity for all such variations gives the optimality conditions,
\begin{subequations}
\begin{align}
\dot{q} &= f(q,u), \label{OCPDAEoptimality1} \\
\dot{p} &= -[D_qf(q,u)]^*p - [D_q\phi(q,u)]^*\lambda - \nabla_qL(q,u), \label{OCPDAEoptimality2} \\
0 &= \phi(q,u), \label{OCPDAEoptimality3} \\
0 &= \nabla_u L(q,u) + [D_uf(q,u)]^*p + [D_u\phi(q,u)]^*\lambda, \label{OCPDAEoptimality4} \\
0 &= \phi_f(q(t_f)), \label{OCPDAEoptimality5} \\
p(t_f) &= \nabla_q C(q(t_f)) + [D_q\phi_f(q(t_f))]^* \lambda_f. \label{OCPDAEoptimality6}
\end{align}
\end{subequations}
The first four optimality conditions \eqref{OCPDAEoptimality1}-\eqref{OCPDAEoptimality4} are precisely the augmented adjoint DAE equations, \eqref{Augmented Adjoint DAE System Coord 1}-\eqref{Augmented Adjoint DAE System Coord 4}. The last two optimality conditions \eqref{OCPDAEoptimality5}, \eqref{OCPDAEoptimality6} are the terminal constraint and the associated transversality condition, respectively. Note that these conditions are only sufficient for a trajectory $(q,u,p,\lambda)$ to be an extremum of the optimal control problem; whether or not the trajectory is optimal depends on the properties of the DAE constraint and cost function, e.g., convexity of $L$.

\textbf{Regular Index 1 Optimal Control.} In the literature, the problem (OCP-DAE) is usually formulated by making a distinction between algebraic variables and control variables, $(q,y,u)$, instead of $(q,u)$ (see, for example, \citet{Bi2010} and \citet{AgCaFo2021}). This does not change any of the previous discussion of the optimality conditions, except that \eqref{OCPDAEoptimality4} splits into two equations for $y$ and $u$. That is, the distinction is not formally important for the previous discussion. It is of course important when actually solving such an optimal control problem. For example, the constraint function $\phi(q,y,u)$ may have a singular matrix derivative with respect to $(y,u)$ but may have a nonsingular matrix derivative with respect to $y$. In such a case, one interprets $y$ as the algebraic variable, in that it can locally be solved in terms of $(q,u)$ via the constraint, and the control variable $u$ as ``free" to optimize over. We now briefly elaborate on this case.

We take the configuration manifold for the algebraic variables to be $M_a = Y_a \times U \ni (y,u)$, where $y$ is interpreted as the algebraic constraint variable and $u$ is interpreted as the control variable. We will assume that the control space $U$ is compact. The constraint has the form $\phi(q,y,u) = 0$, and we assume that $\partial\phi/\partial y$ is pointwise invertible. We consider the following optimal control problem,
\begin{align*}
&\text{min } \int_0^{t_f}L(q,y,u)dt \\
&\text{subject to} \\
& \qquad \dot{q} = f(q,y,u), \\
& \qquad 0 = \phi(q,y,u), \\
& \qquad q_0 = q(0). 
\end{align*}
We perform an analogous argument to before, except that, in this case, since $U$ may have a boundary, the optimality for the control variable $u$ will either require $u$ to lie on $\partial U$ or will require the stationarity of the adjoined functional with respect to variations in $u$. In any case, the necessary conditions for optimality can be expressed as
\begin{subequations}
\begin{align}
\dot{q} &= f(q,y,u), \label{SecondOCPDAEoptimality1} \\
\dot{p} &= -[D_qf(q,y,u)]^*p - [D_q\phi(q,y,u)]^*\lambda - \nabla_qL(q,y,u), \label{SecondOCPDAEoptimality2} \\
0 &= \phi(q,y,u), \label{SecondOCPDAEoptimality3} \\
0 &= \nabla_y L(q,y,u) + [D_yf(q,y,u)]^*p + [D_y\phi(q,y,u)]^*\lambda, \label{SecondOCPDAEoptimality4} \\
u &= \argmin_{u' \in U} H_L(q,y,u'), \label{SecondOCPDAEoptimality5} \\
0 &= p(t_f), \label{SecondOCPDAEoptimality6}
\end{align}
\end{subequations}
where $H_L$ is the augmented Hamiltonian $H_L(q,y,u) = L(q,y,u) + \langle p,f(q,y,u)\rangle + \langle \lambda,\phi(q,y,u)\rangle$. Assuming that $u$ lies in the interior of $U$, \eqref{SecondOCPDAEoptimality5} can be expressed as
$$ 0 = \nabla_u L(q,y,u) + [D_uf(q,y,u)]^*p + [D_u\phi(q,y,u)]^*\lambda, $$
or $D_u H_L(q,y,u) = 0.$ We say that an optimal control problem with a DAE constraint forms a regular index 1 system if both $\partial\phi/\partial y$ and the Hessian $D_u^2 H_L$ are pointwise invertible. In this case, whenever $u$ lies on the interior of $U$, $(y,u,\lambda)$ can be locally solved as functions of $(q,p)$. Thus, in principle, the resulting Hamiltonian ODE for $(q,p)$ can be integrated to yield extremal trajectories for the optimal control problem. As mentioned before, without additional assumptions on the DAE and cost function, such a trajectory will only generally be an extremum but not necessarily optimal.

Of course, in practice, one cannot generally analytically integrate the resulting ODE nor determine the functions which give $(y,u,\lambda)$ in terms of $(q,p)$. Thus, the only practical option is to discretize the presymplectic system above to compute approximate extremal trajectories. To  integrate such a presymplectic system, one can again use the presymplectic Galerkin Hamiltonian variational integrator construction discussed in Section \ref{DiscreteAdjointSystemsSection}. Such an integrator would be natural in the following sense. First, as discussed in Section \ref{DiscreteAdjointSystemsSection}, a presymplectic Galerkin Hamiltonian variational integrator applied to the augmented adjoint DAE system is equivalent to applying a symplectic Galerkin Hamiltonian variational integrator to the underlying Hamiltonian ODE, with the same Runge--Kutta expansions for $q_1, Q^i$ in both methods. Furthermore, as shown in \citet{Sa2016}, utilizing a symplectic integrator to discretize the extremality conditions is equivalent to first discretizing the ODE constraint by a Runge--Kutta method and then enforcing the associated discrete extremality conditions. This also holds in the DAE case.

More precisely, beginning with a regular index 1 optimal control problem, the processes of reduction, extremization, and discretization commute, for suitable choices of these processes, analogous to those used in the naturality result discussed in Section \ref{NaturalityDiscretizationSection}. The proof is similar to the naturality result discussed in Section \ref{NaturalityDiscretizationSection}, where the arrow given by forming the adjoint is replaced by extremization. In essence, these are the same, since the extremization condition is given by the adjoint system, so we will just elaborate briefly. We already know how to extremize the continuous optimal control problem, with either a DAE constraint or an ODE constraint after reduction, which results in an adjoint system. We also already know how to discretize the resulting adjoint system after discretization, using a (pre)symplectic partitioned Runge--Kutta method. Furthermore, at any step, reduction is just defined to be solving the continuous or discrete constraints for $y$ in terms of $(q,u)$. Thus, the only major difference compared to the previous naturality result is defining the discretization of the optimal control problem and subsequently, how to extremize the discrete optimal control problem. For the regular index 1 optimal control problem, 
\begin{align*}
&\text{min } \int_0^{t_f}L(q,y,u)dt \\
&\text{subject to} \\
& \qquad \dot{q} = f(q,y,u), \\
& \qquad 0 = \phi(q,y,u), \\
& \qquad q_0 = q(0),
\end{align*}
its discretization is obtained by replacing the constraints with a Runge--Kutta discretization and replacing the cost function with its quadrature approximation, using the same quadrature weights as those in the Runge--Kutta discretization. This can be written as
\begin{align*}
&\text{min } \Delta t \sum_i b_i L(Q^i,Y^i,U^i) \\
&\text{subject to} \\
& \qquad V^i = f(Q^i,Y^i,U^i), \\
& \qquad 0 = \phi(Q^i,Y^i,U^i),
\end{align*}
where $Q^i = q_0 + \Delta t\sum_j a_{ij}V^j$, which implicitly encodes $q(0)=q_0$. One can then extremize this discrete system, which is given by the discrete Euler--Lagrange equations for the discrete action
$$ \mathbb{S} = \Delta t \sum_i b_i \Big( \langle P^i,V^i-f(Q^i,Y^i,U^i) \rangle - \langle \Lambda^i,\phi(Q^i,Y^i,U^i)\rangle - L(Q^i,Y^i,U^i) \Big). $$
That is, we enforce the discrete constraints by adding to the discrete Lagrangian the appropriate Lagrange multiplier terms paired with the constraints, where we weighted the Lagrange multipliers $P^i,\Lambda^i$ by $\Delta t b_i$ just as convention, in order to interpret them as the appropriate variables, as discussed in Appendix \ref{AppendixNaturalityProof}. Enforcing extremality of this action recovers a partitioned Runge--Kutta method applied to the adjoint system corresponding to extremizing the continuous optimal control problem, as discussed in Appendix \ref{AppendixNaturalityProof}, where the Runge--Kutta coefficients for the momenta are the symplectic adjoint of the original Runge--Kutta coefficients. Alternatively, starting from the original continuous optimal control problem, one could first reduce the DAE constraint to an ODE constraint using the invertibility of $D_y\phi$ to give 
\begin{align*}
&\text{min } \int_0^{t_f}L(q,y(q,u),u)dt \\
&\text{subject to} \\
& \qquad \dot{q} = f(q,y(q,u),u), \\
& \qquad q_0 = q(0).
\end{align*}
One can then discretize this using the same Runge--Kutta method as before, where the cost function is replaced with a quadrature approximation, and then extremize using Lagrange multipliers. Alternatively, one can extremize the continuous problem to yield an adjoint system and then apply a partitioned Runge--Kutta method to that system, where the momenta Runge--Kutta coefficients are again the symplectic adjoint of the original Runge--Kutta coefficients. Having defined all of these processes, a direct computation yields that all of the processes commute, analogous to the computation in Appendix \ref{AppendixNaturalityProof}.

\section{Conclusion and Future Research Directions}
In this paper, we utilized symplectic and presymplectic geometry to study the properties of adjoint systems associated with ODEs and DAEs, respectively. The (pre)symplectic structure of these adjoint systems led us to a geometric characterization of the adjoint variational quadratic conservation law used in adjoint sensitivity analysis. As an application of this geometric characterization, we constructed structure-preserving discretizations of adjoint systems by utilizing (pre)symplectic integrators, which led to natural discrete analogues of the quadratic conservation laws. 

A natural research direction is to extend the current framework to adjoint systems for differential equations with nonholonomic constraints, in order to more generally allow for constraints between configuration variables and their derivatives. In this setting, it is reasonable to expect that the geometry of the associated adjoint systems can be described using Dirac structures (see, for example, \citet{YoMa2006a, YoMa2006b}), which generalize the symplectic and presymplectic structures of adjoint ODE and DAE systems, respectively. Structure-preserving discretizations of such systems could then be studied through the lens of discrete Dirac structures (\citet{LeOh2008}). These discrete Dirac structures make use of the notion of a retraction (\citet{AbMaSe2008}). The tangent and cotangent lifts of a retraction also provide a useful framework for constructing geometric integrators (\citet{BaMa2021}). It would be interesting to synthesize the notion of tangent and cotangent lifts of retraction maps with discrete Dirac structures in order to construct discrete Dirac integrators for adjoint systems with nonholonomic constraints which generalize the presymplectic integrators constructed in \citet{BaMa2022}.

Another natural research direction is to extend the current framework to evolutionary partial differential equations (PDEs). There are two possible approaches in this direction. The first is to consider evolutionary PDEs as ODEs evolving on infinite-dimensional spaces, such as Banach or Hilbert manifolds. One can then investigate the geometry of the infinite-dimensional symplectic structure associated with the corresponding adjoint system. In practice, adjoint systems for evolutionary PDEs are often formed after semi-discretization, leading to an ODE on a finite-dimensional space. Understanding the reduction of the infinite-dimensional symplectic structure of the adjoint system to a finite-dimensional symplectic structure under semi-discretization could provide useful insights into structure-preservation. The second approach would be to explore the multisymplectic structure of the adjoint system associated with a PDE. This approach would be insightful for several reasons. First, an adjoint variational quadratic conservation law arising from multisymplecticity would be adapted to spacetime instead of just time. With appropriate spacetime splitting and boundary conditions, such a quadratic conservation law would induce either a temporal or spatial conservation law. As such, one could use the multisymplectic conservation law to determine adjoint sensitivities for a PDE with respect to spatial or temporal directions, which could be useful in practice \citep{LiPe2004}. Furthermore, the multisymplectic framework would apply equally as well to nonevolutionary (elliptic) PDEs, where there is no interpretation of a PDE as an infinite-dimensional evolutionary ODE. Additionally, adjoint systems for PDEs with constraints could be investigated with multi-Dirac structures (\citet{VaYoLeMa2012}). In future work, we aim to explore both approaches, relate them once a spacetime splitting has been chosen, and investigate structure-preserving discretizations of such systems by utilizing the multisymplectic variational integrators constructed in \citet{TrLe2022}.

\section*{Acknowledgements}

BT was supported by the NSF Graduate Research Fellowship DGE-2038238, and by NSF under grants DMS-1813635. ML was supported by NSF under grants DMS-1345013, DMS-1813635, and by AFOSR under grant FA9550-18-1-0288. 

\section*{Data Availability Statement}
Data sharing is not applicable to this article as no datasets were generated or analyzed during the current study.

\appendix

\section{Proofs of Discrete Adjoint Variational Quadratic Conservation Laws}\label{QuadraticConservationProofs}

\textbf{Proof of Proposition \ref{DiscreteAugmentedQuadraticInvariantProp}.}
\begin{proof}
We begin by substituting \eqref{SPRKAugmentedAdjointODE3} and \eqref{DiscreteVariationalEquations1} into the left hand side of \eqref{DiscreteAugmentedQuadraticInvariant},
\begin{align*}
\langle p_1, \delta q_1\rangle &= \langle p_0, \delta q_0\rangle + \Delta t \sum_i b_i\langle p_0, Df(Q^i)\delta Q^i\rangle - \Delta t \sum_i b_i \langle [Df(Q^i)]^*P^i, \delta q_0\rangle - \Delta t\sum_i b_i\langle dL(Q^i), \delta q_0\rangle \\
& \qquad - \Delta t^2 \sum_{ij} b_ib_j \langle [Df(Q^i)]^*P^i, Df(Q^j)\delta Q^j\rangle - \Delta t^2 \sum_{ij}b_ib_j \langle dL(Q^i), Df(Q^j)\delta Q^j\rangle \\
&= \langle p_0, \delta q_0\rangle + \Delta t \sum_i b_i\left\langle P^i + \Delta t \sum_j\tilde{a}_{ij}([Df(Q^j)]^*P^j + dL(Q^j)), Df(Q^i)\delta Q^i\right\rangle \\
& \qquad - \Delta t \sum_i b_i \left\langle [Df(Q^i)]^*P^i, \delta Q^i - \Delta t \sum_j a_{ij} Df(Q^j)\delta Q^j \right\rangle - \Delta t\sum_ib_i \langle dL(Q^i), \delta q_0\rangle \\
& \qquad - \Delta t^2 \sum_{ij} b_ib_j \langle [Df(Q^i)]^*P^i, Df(Q^j)\delta Q^j\rangle - \Delta t^2 \sum_{ij}b_ib_j \langle dL(Q^i), Df(Q^j)\delta Q^j\rangle,
\end{align*}
where, in the last equality, we substituted \eqref{SPRKAugmentedAdjointODE4} and \eqref{DiscreteVariationalEquations2}. We now group and simplify the above expression,
\begin{align*}
\langle p_1, \delta q_1\rangle &= \langle p_0, \delta q_0\rangle + \Delta t \sum_i b_i\left\langle P^i + \Delta t \sum_j\tilde{a}_{ij}[Df(Q^j)]^*P^j, Df(Q^i) \delta Q^i\right\rangle \\
& \qquad - \Delta t \sum_i b_i \left\langle [Df(Q^i)]^*P^i, \delta Q^i - \Delta t \sum_j a_{ij} Df(Q^j)\delta Q^j \right\rangle \\
& \qquad - \Delta t^2 \sum_{ij} b_ib_j \langle [Df(Q^i)]^*P^i, Df(Q^j)\delta Q^j\rangle \\
& \qquad + \Delta t \sum_i b_i\left\langle \Delta t \sum_j\tilde{a}_{ij}dL(Q^j), Df(Q^i)\delta Q^i\right\rangle - \Delta t\sum_ib_i \langle dL(Q^i), \delta q_0\rangle \\
& \qquad - \Delta t^2 \sum_{ij}b_ib_j \langle dL(Q^i), Df(Q^j)\delta Q^j\rangle \\
&= \langle p_0, \delta q_0\rangle + \Delta t^2 \sum_{ij} \underbrace{(b_j\tilde{a}_{ji} + b_ia_{ij} - b_ib_j)}_{=0} \langle [Df(Q^i)]^*P^i, Df(Q^j)\delta Q^j\rangle \\
& \qquad + \Delta t \sum_i b_i\left\langle \Delta t \sum_j\tilde{a}_{ij}dL(Q^j), Df(Q^i)\delta Q^i\right\rangle - \Delta t\sum_i b_i\langle dL(Q^i), \delta q_0\rangle \\
& \qquad - \Delta t^2 \sum_{ij}b_ib_j \langle dL(Q^i), Df(Q^j)\delta Q^j\rangle \\
&= \langle p_0, \delta q_0\rangle - \Delta t \sum_ib_i\langle dL(Q^i), \delta q_0\rangle \\
&\qquad - \Delta t^2 \sum_{ij} \underbrace{(b_ib_j - b_j\tilde{a}_{ji})}_{=b_ia_{ij}} \langle dL(Q^i), Df(Q^j)\delta Q^j\rangle \\
&= \langle p_0, \delta q_0\rangle - \Delta t^2 \sum_i b_i \left\langle dL(Q^i), \frac{\delta q_0}{\Delta t} + \sum_j a_{ij}Df(Q^j)\delta Q^j\right\rangle \\
&= \langle p_0, \delta q_0\rangle - \Delta t\sum_i b_i \langle dL(Q^i), \delta Q^i\rangle,
\end{align*}
where, in the last equality, we used \eqref{DiscreteVariationalEquations2}.
\end{proof}

\textbf{Proof of Proposition \ref{DiscreteDAEQuadraticInvariantProp}.}
\begin{proof}
For brevity, we denote
\begin{align*}
D_qf_i &\equiv D_qf(Q^i,U^i), \\
D_uf_i &\equiv D_uf(Q^i,U^i), \\
D_q\phi_i &\equiv D_q\phi(Q^i,U^i), \\
D_u\phi_i &\equiv D_u\phi(Q^i,U^i).
\end{align*}
Starting from $\langle p_1,\delta q_1\rangle$, we substitute the evolution equations \eqref{SPRKAdjointDAE3}, \eqref{SPRKAdjointDAE4}, \eqref{DiscreteDAEVariationalEquations1}, \eqref{DiscreteDAEVariationalEquations2},
\begin{align*}
\langle p_1, \delta q_1 \rangle &= \langle p_0, \delta q_0\rangle - \Delta t \sum_i b_i \langle [D_qf_i]^*P^i + [D_q\phi_i]^*\Lambda^i, \delta q_0\rangle 
+ \Delta t \sum_i b_i \langle p_0, D_qf_i \delta Q^i + D_uf_i \delta U^i\rangle \\
&\qquad - \Delta t^2 \sum_{ij}b_ib_j \langle [D_qf_i]^*P^i + [D_q\phi_i]^*\Lambda^i, D_qf_j \delta Q^j + D_uf_j \delta U^j)\rangle \\
&= \langle p_0, \delta q_0\rangle - \Delta t\sum_ib_i \left\langle [D_qf_i]^*P^i + [D_q\phi_i]^*\Lambda^i, \delta Q^i - \Delta t \sum_j a_{ij}(D_qf_j \delta Q^j + D_uf_j \delta U^j)\right\rangle  \\
& \qquad + \Delta t \sum_i b_i \left\langle P^i + \Delta t \sum_j \tilde{a}_{ij}( [D_qf_j]^*P^j + [D_q\phi_j]^*\Lambda^j), D_qf_i \delta Q^i + D_uf_i \delta U^i\right\rangle \\
&\qquad - \Delta t^2 \sum_{ij}b_ib_j \langle [D_qf_i]^*P^i + [D_q\phi_i]^*\Lambda^i, D_qf_j \delta Q^j + D_uf_j \delta U^j\rangle \\
&= \langle p_0, \delta q_0\rangle - \Delta t \sum_i b_i \langle [D_qf_i]^*P^i + [D_q\phi_i]^*\Lambda^i, \delta Q^i\rangle \\
&\qquad + \Delta t \sum_i b_i \langle P^i, D_qf_i \delta Q^i + D_uf_i \delta U^i\rangle \\
&\qquad + \Delta t^2 \sum_{ij}\underbrace{(b_j\tilde{a}_{ji} + b_ia_{ij} - b_ib_j)}_{=0} \langle [D_qf_i]^*P^i + [D_q\phi_i]^*\Lambda^i, D_qf_j \delta Q^j + D_uf_j \delta U^j\rangle \\
&= \langle p_0, \delta q_0\rangle + \Delta t \sum_i b_i \Big(\langle [D_uf_i]^*P^i, \delta U^i\rangle - \langle [D_q\phi_i]^*\Lambda^i,\delta Q^i\rangle \Big) \\
&= \langle p_0, \delta q_0\rangle + \Delta t \sum_i b_i \Big( - \langle [D_u\phi_i]^*\Lambda^i, \delta U^i\rangle - \langle [D_q\phi_i]^*\Lambda^i,\delta Q^i\rangle \Big) \\
&= \langle p_0, \delta q_0\rangle - \Delta t \sum_i b_i \langle \Lambda^i, D_u\phi_i \delta U^i + D_q\phi_i \delta Q^i\rangle \\ 
&= \langle p_0, \delta q_0\rangle,
\end{align*}
where in the third to last equality, we used the constraint equation \eqref{SPRKAdjointDAE6} and in the last equality, we used the constraint equation \eqref{DiscreteDAEVariationalEquations3}.
\end{proof}

\textbf{Proof of Proposition \ref{DiscreteAugmentedAdjointDAEQuadraticProp}.}
\begin{proof}
The proof uses computations analogous to those used in the proofs of Propositions \ref{DiscreteAugmentedQuadraticInvariantProp} and \ref{DiscreteDAEQuadraticInvariantProp}. In particular, starting from the simplest case of the nonaugmented adjoint ODE system, Proposition \ref{DiscreteAugmentedQuadraticInvariantProp} considers the case of augmenting the Hamiltonian, whereas Proposition \ref{DiscreteDAEQuadraticInvariantProp} considers the case of replacing the ODE with a DAE. The case at hand combines both and the proof involves a combination of both computations.
\end{proof}

\section{Proof of Naturality of Adjoint System Discretization}\label{AppendixNaturalityProof}
In this appendix, we prove the statement in Section \ref{NaturalityDiscretizationSection} that (for suitable choices of) discretization, reduction, and forming the adjoint all commute when applied to an index 1 DAE. The definitions and choices of these processes were made in Section \ref{NaturalityDiscretizationSection}. To prove that the diagram commutes, we prove that each face of the diagram commutes. We again include the relevant diagram which we wish to show commutes below.

\[\small\begin{tikzcd}[column sep=-4ex,row sep=10ex]
	\txt{Index 1 DAE} && \txt{ODE} \\
	& \txt{Discrete DAE} && \txt{Discrete ODE} \\
	\txt{Presymplectic Adjoint\\ DAE System} && \txt{Symplectic Adjoint\\ ODE System} \\
	& \txt{Presymplectic Galerkin\\ Hamiltonian Variational \\ Integrator} && \txt{Symplectic Galerkin \\ Hamiltonian Variational \\ Integrator}
	\arrow["{\text{Reduce}}", from=1-1, to=1-3]
	\arrow["{\text{Reduce}}"{pos=0.25}, from=3-1, to=3-3]
	\arrow["{\text{Adjoint}}"'{pos=0.6}, from=1-1, to=3-1]
	\arrow["{\text{Adjoint}}"{pos=0.6}, from=1-3, to=3-3]
	\arrow["{\text{Reduce}}"{pos=0.3}, from=2-2, to=2-4,crossing over]
	\arrow["{\text{Reduce}}"', from=4-2, to=4-4]
	\arrow["{\text{Adjoint}}"{pos=0.6}, from=2-2, to=4-2,crossing over]
	\arrow["{\text{Adjoint}}"{pos=0.6}, from=2-4, to=4-4]
	\arrow["{\text{Discretize}}"', from=1-1, to=2-2]
	\arrow["{\text{Discretize}}"{pos=0.6}, from=1-3, to=2-4]
	\arrow["{\text{Discretize}}"', from=3-1, to=4-2]
	\arrow["{\text{Discretize}}"', from=3-3, to=4-4]
\end{tikzcd}\]

\textbf{Back Face.} We have already proved that the back face commutes (i.e., that reduction and forming the adjoint commute when starting with an index 1 DAE), as discussed in Section \ref{DAEIndexPCASection}. One can then interpret the above diagram as an extension of this result with an extra dimension corresponding to  discretization.

\textbf{Right Face.} This was proven in \citet{Sa2016}. One can then interpret the above diagram as an extension of the result in \citet{Sa2016} by adding the reduction operation.

\textbf{Bottom Face.} Consider the augmented adjoint DAE system corresponding to the DAE \eqref{DAEa}-\eqref{DAEb}, which we take to have index 1, i.e., $\partial\phi/\partial u$ is pointwise invertible. We consider the augmented case because the nonaugmented case can be obtained by taking $L \equiv 0$. We show that reducing the system first and then applying a symplectic Galerkin Hamiltonian variational integrator is equivalent to applying a presymplectic Galerkin Hamiltonian variational integrator, with the same partitioned Runge--Kutta coefficients, and then reducing.

We start with the former approach. The symplectic adjoint ODE system given by reduction, as discussed in Section \ref{DAEIndexPCASection}, is the Hamiltonian system corresponding to the Hamiltonian
$$ H(q',p') = \langle p',f'(q'))\rangle + L'(q'), $$
where we have solved $u = u(q')$ and defined $f'(q') \equiv f(q', u(q')), L'(q') \equiv L(q', u(q'))$. Applying the symplectic Galerkin Hamiltonian variational integrator construction yields the integrator
\renewcommand{\theequation}{\thesection.\arabic{equation}}
\begin{subequations}
\begin{align}
q_1 &= q_0 + \Delta t \sum_i b_i f'(Q^i) \label{DiscreteReductionHamiltonian1} \\
    &= q_0 + \Delta t \sum_i b_i f(Q^i, u(Q^i)), \nonumber \\
Q^i &= q_0 + \Delta t \sum_j a_{ij} f'(Q^j) \label{DiscreteReductionHamiltonian2} \\
    &= q_0 + \Delta t \sum_j a_{ij} f(Q^j, u(Q^j)), \nonumber \\
p_1 &= p_0 - \Delta t \sum_i b_i ([Df'(Q^i)]^*P^i + dL'(Q^i)) , \label{DiscreteReductionHamiltonian3} \\
P^i &= p_0 - \Delta t \sum_j \tilde{a}_{ij} ([Df'(Q^j)]^*P^j + dL'(Q^j)). \label{DiscreteReductionHamiltonian4}
\end{align}
\end{subequations}
Note that the derivative $Df'$ can be equivalently expressed as
\begin{align*}
Df'(Q^i) &= D_1f(Q^i, u(Q^i)) + D_2f(Q^i, u(Q^i)) Du(Q^i),
\end{align*}
where $D_i$ denotes differentiation with respect to the $i^{th}$ argument. We switch to indexing the derivative operator here, so we do not have to make the distinction between total derivatives $D_q$ and partial derivatives $\partial_q$. Similarly, we can express $dL'$ as follows. First, note that we have been implicitly identifying the row vector $dL'$ with the column vector given by its transpose $\nabla L'$. Thus, $dL'$ in equations \eqref{DiscreteReductionHamiltonian3}-\eqref{DiscreteReductionHamiltonian4} should really be written as $\nabla L'$. Thus,
$$ dL'(Q^i) \cong \nabla L'(Q^i) = \nabla_1L(Q^i, u(Q^i)) + [Du(Q^i)]^*\nabla_2 L(Q^i, u(Q^i)). $$
Now, we show that the second approach is equivalent to the above system. The starting point is the presymplectic Galerkin Hamiltonian variational integrator, equations \eqref{SPRKAugmentAdjointDAE1}-\eqref{SPRKAugmentAdjointDAE6}. From \eqref{SPRKAugmentAdjointDAE5}, we can solve for $U^i$ in terms of $Q^i$ as $U^i = u(Q^i)$. Plugging this into \eqref{SPRKAugmentAdjointDAE1}-\eqref{SPRKAugmentAdjointDAE2} gives precisely \eqref{DiscreteReductionHamiltonian1}-\eqref{DiscreteReductionHamiltonian2}. Thus, we just need to see that, after solving the constraint \eqref{SPRKAugmentAdjointDAE6} for $\Lambda^i$, the two momenta equations \eqref{SPRKAugmentAdjointDAE3}-\eqref{SPRKAugmentAdjointDAE4} are equivalent to \eqref{DiscreteReductionHamiltonian3}-\eqref{DiscreteReductionHamiltonian4}. Solving \eqref{SPRKAugmentAdjointDAE6} for $\Lambda^i$ gives
$$ \Lambda^i = - ([D_2\phi(Q^i,u(Q^i))]^*)^{-1} [D_2f(Q^i,u(Q^i))]^*P^i - ([D_2\phi(Q^i,u(Q^i))]^*)^{-1}  \nabla_2L(Q^i,u(Q^i)). $$
Multiplying both sides by $[D_1\phi(Q^i,u(Q^i))]^*$ yields
\begin{align*}
[D_1\phi(Q^i,u(Q^i))]^*\Lambda^i &= - [D_1\phi(Q^i,u(Q^i))]^*([D_u\phi(Q^i,u(Q^i))]^*)^{-1} [D_2f(Q^i,u(Q^i))]^*P^i \\
& \qquad - [D_1\phi(Q^i,u(Q^i))]^*([\nabla_u\phi(Q^i,u(Q^i))]^*)^{-1}  \nabla_2L(Q^i,u(Q^i))  \\
&= [Du(Q^i)]^*[D_2f(Q^i,u(Q^i))]^* P^i + [Du(Q^i)]^* \nabla_2L(Q^i,u(Q^i)) \\
&= [D_2f(Q^i,u(Q^i))Du(Q^i)]^*P^i + [Du(Q^i)]^* \nabla_2L(Q^i,u(Q^i)),
\end{align*}
where in the second equality, we used $D_1\phi(Q^i,u(Q^i)) = - D_2\phi(Q^i, u(Q^i)) Du(Q^i)$ from the implicit function theorem. Plugging this expression and $U^i = u(Q^i)$ into \eqref{SPRKAugmentAdjointDAE3}-\eqref{SPRKAugmentAdjointDAE4} yields \eqref{DiscreteReductionHamiltonian3}-\eqref{DiscreteReductionHamiltonian4}, noting the above expressions for $Df', dL'$.

\begin{remark}
Note that, in the above, we used the implicit function theorem to obtain the local function $u = u(q)$. This is sufficient to prove that the two processes are the same for a single integration step, assuming that the timestep $\Delta t$ is sufficiently small and the vector field $f$ and constraint $\phi$ are sufficiently regular, so that $q_0$, $q_1$, and all of the internal stages $Q^i$ are in the neighborhood where the local function is defined. For each subsequent time step, one generally needs a different local function. This does not matter in practice since one works directly with the presymplectic integrator and solves the constraints iteratively.
\end{remark}

\textbf{Top Face.} We want to prove that, starting from an index 1 DAE, the processes of discretization and reduction commute, where the discretization of the ODE and DAE have the same Runge--Kutta coefficients. 

We start first with reduction then discretization. Starting from the index 1 DAE $\dot{q} = f(q,u)$, $\phi(q,u) = 0$, we apply the reduction operation, which gives the ODE $\dot{q} = f(q, u(q))$. Applying a Runge--Kutta discretization gives
\begin{align*}
q_1 &= q_0 + \Delta t \sum_i b_i f(Q^i,u(Q^i)), \\
Q^i &= q_0 + \Delta t \sum_j a_{ij} f(Q^j, u(Q^j)). \\
\end{align*}

On the other hand, we can discretize the DAE and then reduce. We discretize the DAE $\dot{q} = f(q,u)$, $\phi(q,u) = 0$ by applying a Runge--Kutta discretization with the same coefficients as before,
\begin{align*}
q_1 &= q_0 + \Delta t \sum_i b_i f(Q^i,U^i), \\
Q^i &= q_0 + \Delta t \sum_j a_{ij} f(Q^j, U^j), \\
0 &= \phi(Q^i, U^i).
\end{align*}
To reduce this system, we solve the constraint equations $U^i = u(Q^i)$ and substitute these into the two evolution equations, which yields the same system obtained from first reducing and then discretizing.

\textbf{Front Face.} The starting point for this loop is a discrete DAE system, which arises as a Runge--Kutta discretization of an index 1 DAE, i.e., it is given by the discrete system
\begin{subequations}
\begin{align}
q_1 &= q_0 + \Delta t \sum_i b_i f(Q^i,U^i), \label{DiscreteDAESystem1} \\
Q^i &= q_0 + \Delta t \sum_j a_{ij} f(Q^j, U^j), \label{DiscreteDAESystem2} \\
0 &= \phi(Q^i, U^i). \label{DiscreteDAESystem3}
\end{align}
\end{subequations}
From here, we wish to show that reducing and forming the discrete  adjoint system commute. 

First, we recall the notion of a discrete adjoint system. Suppose we are given a generally nonlinear system of equations, $F(x_1) = x_0$, where $x_1 \in V$ is unknown, $x_0 \in W$ is given, and $F: V \rightarrow W$ (where $V$ and $W$ are vector spaces). To define the adjoint system, we first consider the variational equations associated with this nonlinear system given by its linearization,
$$ DF(x_1)\delta x_1 = \delta x_0, $$
where $DF(x_1)$ is a linear map $V \rightarrow W$ and $\delta x_0 \in W$ is given. Suppose that we are interested in computing the quantity $\langle s_1, \delta x_1\rangle$ for a given vector $s_1 \in V^*$. In the setting of adjoint sensitivity analysis, the quantity $\langle s_1, \delta x_1\rangle$ is the sensitivity of the terminal cost function. We define the associated adjoint equation as
$$ [DF(x_1)]^*s_0 = s_1. $$
For a solution $s_0 \in W^*$ of this system, one has
$$ \langle s_1, \delta x_1\rangle = \langle [DF(x_1)]^*s_0, \delta x_1\rangle = \langle s_0, DF(x_1)\delta x_1\rangle = \langle s_0, \delta x_0\rangle. $$
Thus, to compute $\langle s_1, \delta x_1\rangle$, one could solve the variational equation for $\delta x_1$ and pair it with $s_1$ which is given, or, alternatively, solve the adjoint equation for $s_0$ and pair it with $\delta x_0$ which is given, since these linear systems are solvable by assumption. We define the adjoint system associated with the equation $F(x_1) = x_0$ as this equation combined with the associated adjoint equation, i.e., as the combined system
\begin{align*}
F(x_1) &= x_0, \\
[Df(x_1)]^* s_0 &= s_1.
\end{align*}
Following \citet{Ib2006}, we will utilize an alternative characterization of the adjoint system. We define the discrete adjoint action
$$ \mathbb{S}(x_1, s_0) \equiv \langle s_0, F(x_1)\rangle. $$
Then, observe that $\mathbb{S}$ is a generating function for the adjoint system $(x_1, s_0) \mapsto (x_0,s_1)$, in the sense that
\begin{align*}
x_0 &= \frac{\delta}{\delta s_0} \mathbb{S}(x_1,s_0) = F(x_1), \\
s_1 &= \frac{\delta}{\delta x_1} \mathbb{S}(x_1,s_0) = [Df(x_1)]^*s_0.
\end{align*}
This characterization serves two purposes. First, it will simplify the calculation of the adjoint system for the case at hand. Furthermore, it resembles the process of forming the adjoint at the continuous level: starting from the (discrete or continuous) differential(-algebraic) equation at hand, one forms the (discrete or continuous) adjoint action and applies the variational principle to obtain the adjoint system. To obtain the augmented adjoint system, we add a discrete Lagrangian $\mathbb{L}: V \rightarrow \mathbb{R}$ to the action (as a convention, we subtract the discrete Lagrangian). We define the augmented discrete adjoint action to be
$$ \mathbb{S}_{\mathbb{L}}(x_1,s_0) \equiv \langle s_0,F(x_1)\rangle - \mathbb{L}(x_1). $$
The map that this generates defines the augmented discrete adjoint system,
\begin{align*}
x_0 &= \frac{\delta}{\delta s_0} \mathbb{S}_{\mathbb{L}}(x_1,s_0) = F(x_1), \\
s_1 &= \frac{\delta}{\delta x_1} \mathbb{S}_{\mathbb{L}}(x_1,s_0) = [Df(x_1)]^*s_0 - d\mathbb{L}(x_1).
\end{align*}
Observe that this definition of an augmented discrete adjoint system is natural in the sense that,
$$ \langle s_1, \delta x_1\rangle = \langle [Df(x_1)]^*s_0 + d\mathbb{L}(x_1), \delta x_1\rangle = \langle s_0,x_0 \rangle - \langle dL(x_1),\delta x_1\rangle,  $$
which resembles the continuous analogue of the adjoint sensitivity result for a running cost function.

Now, we use this notion of a discrete adjoint system for the problem at hand. We begin first with reduction and then forming the adjoint system. Applying the reduction operation to the discrete DAE system \eqref{DiscreteDAESystem1}-\eqref{DiscreteDAESystem3}, given by solving $\phi(Q^i,U^i) = 0$ for $U^i = u(Q^i)$,
\begin{subequations}
\begin{align}
q_1 &= q_0 + \Delta t \sum_i b_i f(Q^i,u(Q^i)), \label{DiscreteDAEReduced1} \\
Q^i &= q_0 + \Delta t \sum_j a_{ij} f(Q^j, u(Q^j)), \label{DiscreteDAEReduced2}
\end{align}
\end{subequations}
Let us define $Q^i = q_0 + \Delta t \sum_j a_{ij} V^j$. We think of the internal stages $Q$ as functions of the internal stages $V$, which are the internal stage proxies for $\dot{q}$. Our discrete system \eqref{DiscreteDAEReduced1}-\eqref{DiscreteDAEReduced2} can then be defined by $x_1 = \{V^i\}_{i=1}^s$, $x_0 = \{0\}_{i=1}^s$, where $s$ is the number of internal stages, and
$$ x_0 = F(x_1) \equiv \begin{pmatrix} V^1 - f(Q^i(V), u(Q^i(V))) \\ \vdots \\ V^s - f(Q^s(V), u(Q^s(V))) \end{pmatrix}. $$
Observe that $F = 0$ only gives the internal stage equations \eqref{DiscreteDAEReduced2}. We do this for simplicity, since we will assume $c_s = 1$ as is typical for a Runge--Kutta discretization of a DAE as previously discussed and hence, equation \eqref{DiscreteDAEReduced1} is redundant, since $a_{sj} = b_j$.

We define $F$ and $x_1$ in terms of $V$ instead of $Q$ because when we form the adjoint action, we pair the components of $F$ with the dual variable $s_0$. In order to interpret $s_0$ as representing the momenta internal stages $P^i$, it should be paired with the proxy for the tangent vector $V$, instead of $Q$. We now form the discrete adjoint action. We define the dual variable for the adjoint system to be $s_0 = \{ \Delta t b_i P^i \}_{i=1}^s$. The normalization factor $\Delta t b_i$ is used so that the discrete action is the quadrature approximation of the continuous action. This is just a convention, but we would have to reinterpret the components of $s_0$ if we did not choose this convention. Finally, we define the discrete Lagrangian to be the quadrature approximation of the continuous Lagrangian $L'(q) \equiv L(q,u(q))$, i.e., $\mathbb{L}(x_1) = \Delta t\sum_i b_i L'(Q^i(V))$. This is the natural choice because the discrete sensitivity of a running cost function is $\Delta t \sum_i b_i \langle dL'(Q^i(V)), \delta Q^i(V)\rangle,$ which equals $\langle d\mathbb{L}(x_1),\delta x_1\rangle$ with the above choice of $\mathbb{L}$. The augmented discrete adjoint action is then
\begin{align*}
\mathbb{S}_{\mathbb{L}}(\{V^i\},\{b_iP^i\}) &= \mathbb{S}_{\mathbb{L}}(x_1,s_0) = \langle s_0, F(x_1)\rangle - \mathbb{L}(x_1) \\ 
&= \Delta t \sum_i b_i \Big( \langle P^i, V^i - f(Q^i(V),u(Q^i(V))) \rangle - L'(Q^i(V)) \Big). 
\end{align*}
To define the discrete adjoint system, we have to give $s_1$, which we take to be $s_1 = \{ \Delta t b_i p_1 \}_{i=1}^{s}$, where $p_1$ is given. Thus, the augmented discrete adjoint system is given by
\begin{align*}
0 &= \frac{\delta}{\delta P^k} \mathbb{S}_{\mathbb{L}} = V^k - f(Q^k(V), u(Q^k(V))), \\
\Delta t b_kp_1 &= \frac{\delta}{\delta V^k} \mathbb{S}_{\mathbb{L}} \\
 &= \Delta t b_k P^k - \Delta t^2\sum_i b_i a_{ik} \Big([D_1f(Q^i(V), u(Q^i(V)))]^*P^i \\ 
 & \qquad \qquad \qquad \qquad + [D_2f(Q^i(V), u(Q^i(V))) Du(Q^i(V))]^* P^i + dL'(Q^i(V)) \Big).
\end{align*}
The first set of equations above, combined with the definition of $Q$ in terms of $V$, gives \eqref{DiscreteDAEReduced2}. For the second set of equations, we first divide through by $\Delta t b_k$ and rearrange to obtain
\begin{align*}
P^k &= p_1 + \Delta t\sum_i \frac{b_i a_{ik}}{b_k} \Big([D_1f(Q^i(V), u(Q^i(V)))]^*P^i \\ 
 & \qquad \qquad \qquad \qquad + [D_2f(Q^i(V), u(Q^i(V))) Du(Q^i(V))]^* P^i + dL'(Q^i(V)) \Big). \\
 &= p_1 + \sum_i (b_i - \tilde{a}_{ki}) \Big([D_1f(Q^i(V), u(Q^i(V)))]^*P^i \\ 
 & \qquad \qquad \qquad \qquad + [D_2f(Q^i(V), u(Q^i(V))) Du(Q^i(V))]^* P^i + dL'(Q^i(V)) \Big).
\end{align*}
Note that this is the usual symplectic partitioned Runge--Kutta expansion for the internal stages $P^i$, expressed in terms of $p_1$ instead of $p_0$. Thus, the full adjoint system, combined with the redundant $k=s$ stages, yields a symplectic partitioned Runge--Kutta method. 

Now, in the other direction, we first form the adjoint system corresponding to the discrete DAE system and subsequently reduce. We begin by forming the adjoint system. We form the discrete action analogously to before, but now the discrete system \eqref{DiscreteDAESystem1}-\eqref{DiscreteDAESystem3} also has constraints which we must incorporate into $F$, since we have not yet reduced the system. We take $x_1 = \{ \{V^i\}, \{U^i\} \}_{i=1}^s$ and $s_0 = \{ \{\Delta t b_i P^i\}, \{\Delta t b_i \Lambda^i\} \}_{i=1}^s$. We define $F$ as
$$ x_0 = F(x_1) \equiv \begin{pmatrix} V^1 - f(Q^i(V), U^i) \\ \vdots \\ V^s - f(Q^s(V), U^i)) \\ -\phi(Q^1(V),U^1) \\ \vdots \\ -\phi(Q^s(V),U^s) \end{pmatrix}. $$
Note again that $Q$ is a function of $V$ as $Q^i = q_0 + \Delta t \sum_j a_{ij}V^j$. It is not a priori a function of $U$ because the condition $V^i = f(Q^i(V),U^i)$ has not yet been enforced. Rather, it is a consequence of the variational principle, which formally matters when one computes the variation of the discrete action. Define the discrete Lagrangian $\mathbb{L}(x_1) = \sum_i b_i L(Q^i(V),U^i).$ We form the augmented discrete adjoint action
\begin{align*}
\mathbb{S}_{\mathbb{L}}(\{V^i\},\{b_iP^i\}) &= \mathbb{S}_{\mathbb{L}}(x_1,s_0) = \langle s_0, F(x_1)\rangle - \mathbb{L}(x_1) \\ 
&= \Delta t \sum_i b_i \Big( \langle P^i, V^i - f(Q^i(V),U^i) \rangle - \langle \Lambda^i,\phi(Q^i,U^i)\rangle - L(Q^i(V), U^i) \Big). 
\end{align*}
We use this as a generating function to compute the adjoint system as before. The computation is analogous so we will just state the result, 
\begin{align*} 
V^k &= f(Q^k(V), U^k), \\
P^k &= p_1 + \Delta t \sum_i (b_i - \tilde{a}_{ki}) \left( [D_1f(Q^i(V),U^i)]^*P^i + [D_1\phi(Q^i(V),U^i)]^*\Lambda^i + D_qL(Q^i(V),U^i) \right),  \\
0 &= \phi(Q^i(V),U^i),  \\
0 &= [D_2f(Q^i(V),U^i)]^*P^i + [D_2\phi(Q^i(V),U^i)]^*\Lambda^i + D_2L(Q^i(V),U^i). 
\end{align*}
Finally, we reduce by solving the last two equations for $U^i$, $\Lambda^i$ as functions of $Q^i(V^i)$, $P^i$. Finally, an implicit function theorem computation analogous to the proof of the bottom face shows that this is the same as the system obtained by first reducing and then forming the discrete adjoint. 

\textbf{Left Face.} The proof for the left face is formally similar to the right face, but since we have already computed both directions, we will include it for completeness. Starting from an index 1 DAE, forming the adjoint and then discretizing just give the presymplectic Galerkin Hamiltonian variational integrator \eqref{SPRKAugmentAdjointDAE1}-\eqref{SPRKAugmentAdjointDAE6}. In the other direction, we first discretize the DAE and then take the adjoint which we did in the proof of the front face. Expressed in terms of $Q$, instead of $V$, this is 
\begin{align*} 
Q^k &= q_0 + \Delta t \sum_j a_{ij}f(Q^j, U^j), \\
P^k &= p_1 + \Delta t \sum_i (b_i - \tilde{a}_{ki}) \left( [D_1f(Q^i,U^i)]^*P^i + [D_1\phi(Q^i,U^i)]^*\Lambda^i + D_qL(Q^i,U^i) \right),  \\
0 &= \phi(Q^i,U^i),  \\
0 &= [D_2f(Q^i,U^i)]^*P^i + [D_2\phi(Q^i,U^i)]^*\Lambda^i + D_2L(Q^i,U^i). 
\end{align*}
Returning to the system given by first forming the adjoint and then discretizing, \eqref{SPRKAugmentAdjointDAE1}-\eqref{SPRKAugmentAdjointDAE6}, one substitutes \eqref{SPRKAugmentAdjointDAE3} into \eqref{SPRKAugmentAdjointDAE4} to write the internal stages for $P^i$ in terms of $p_1$, and this gives the above system.

\section{An Intrinsic Type II Variational Principle for Adjoint Systems}\label{typeIIvp appendix}

We show that the adjoint system \eqref{AugmentedAdjointSystem} arises from an intrinsic Type II variational principle. In coordinates, the type II variational principle corresponds to fixed initial position $q(t_0)=q_0$ and fixed final momenta $p(t_1)=p_1$, which are the boundary conditions used in adjoint sensitivity analysis, as discussed in Section \ref{AdjointSensitivitySection}.

Consider the augmented adjoint system 
\begin{align*}
\dot{q} &= \partial H_L/\partial p = f(q), \\
\dot{p} &= -\partial H_L/\partial q = - [Df(q)]^*p - dL(q),
\end{align*}
where $H_L$ is the augmented Hamiltonian. Recall that $H_L$ is intrinsically defined by $H_L = i_{\widehat{f}}\Theta + \pi_{T^*M}^*L$, where $\Theta$ is the tautological one-form on $T^*M$, $\pi_{T^*M}: T^*M \rightarrow M$ is the cotangent bundle projection, and $L: M \rightarrow \mathbb{R}$.

We would like to show that the above system arises from a Type II variational principle. We consider the action
$$ S[\psi] = \int_{t_0}^{t_1} \psi^* \left( \Theta - H_L \, dt \right), $$
where $\psi: (t_0,t_1) \rightarrow T^*M $ is a curve on $T^*M$. 

\begin{figure}[H]
\begin{center}
\includegraphics[width=100mm]{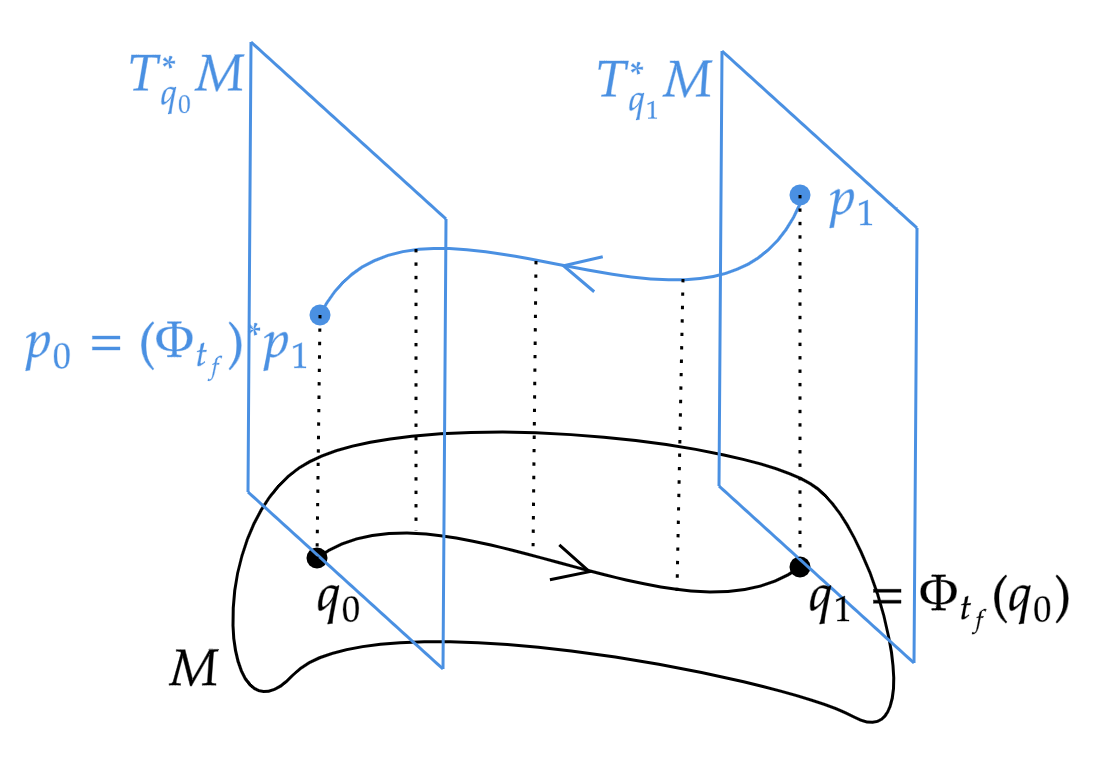}
\caption{Type II boundary conditions for adjoint systems}
\label{typeiibcs}
\end{center}
\end{figure}

We would to place Type II boundary conditions, $q(t_0) = q_0$ and $p(t_1) = p_1$, on the variational principle. However, Type II boundary conditions for Hamiltonian systems, in general, suffer the drawback that they do not make intrinsic sense on a manifold, since one cannot specify a covector $p(t_1) = p_1$ without specifying the basepoint $q(t_1)$. Fortunately, for Hamiltonian systems which are adjoint systems, Type II boundary conditions do make intrinsic sense, due to the fact they cover an ODE on the base manifold $M$. To see this, if we fix the boundary condition $q(t_0) = q_0$, the time $t_1-t_0$ flow of $f$, assuming it exists for this time, fixes the basepoint $q(t_1) = \Phi_{t_1-t_0}(q(t_0))$. In terms of the curve $\psi$, this means that once we fix $\pi_{T^*M}(\psi(t_0)) = q_0$, we have $\psi(t_1) \in T_{q(t_1)}^*M$, where $q(t_1) = \Phi_{t_1-t_0}(q(t_0))$. Thus, it then makes sense to specify a boundary condition on $\psi(t_1) \in T_{q(t_1)}^*M$ of the form $\psi(t_1) = p_1$, for any $p_1 \in T_{q(t_1)}^*M.$ Figure \ref{typeiibcs} illustrates Type II boundary conditions for an adjoint system; the flow of $f$ on the base manifold evolves the initial condition $q_0$ forward to $q_1$ and subsequently, the vertical component of the lifted vector field $X_{H_L}$ evolves the final momenta $p_1$, based at $q_1$, backwards to the initial momenta $p_0$. As discussed in Section \ref{AdjointSensitivitySection}, $p_1$ can be chosen by taking $p_1 = dC|_{q_1}$ to compute the sensivity of a terminal cost function $C:M\rightarrow \mathbb{R}$ with a non-augmented Hamiltonian $H_L=H$ or by taking $p_1=0$ to compute the sensivity of a running cost function $L$ with an augmented Hamiltonian $H_L = H+L$.

\begin{remark}
It is interesting to note that the reason for which Type I boundary conditions for adjoint systems are generally inconsistent (namely, that they cover an ODE on the base manifold) is precisely the reason that one can make intrinsic sense of Type II boundary conditions for adjoint systems. That is, Type II boundary conditions are consistent while Type I boundary conditions are generally inconsistent precisely because an adjoint system is a Hamiltonian system which covers an ODE on the base manifold. Conversely, every Hamiltonian system on $T^*M$ which covers an ODE on the base manifold $M$ is locally an adjoint system. To see this, if a Hamiltonian system covers an ODE on the base manifold, then Hamilton's equation in the position variable $\dot{q} = \partial H/\partial p$ must equal $f(q)$ for some vector field $f$ on $M$. Thus, we have $\partial H/\partial p = f(q)$. Integrating this equation yields a coordinate expression for the Hamiltonian
$$ H(q,p) = \langle p,f(q)\rangle + L(q), $$
where the ``constant of integration" (constant with respect to the $p$ variable) $L(q)$ is some arbitrary function of $q$. This is precisely the form of the Hamiltonian for an augmented adjoint system.
\end{remark}

To state an intrinsic Type II variational principle for adjoint systems, we regard the integrand of the above action (before pulling back by $\psi$) as a contact form on the extended phase space $I \times T^*M$. Namely, given an interval $I = (t_0,t_1) \subset \mathbb{R}$, $t_0 \neq t_1$, let $\pi_{I \times T^*M}: I \times T^*M \rightarrow T^*M$ denote the projection onto the second factor. Then, define the contact form
$$ \Theta_H = \pi_{I \times T^*M}^* \Theta - H dt, $$
where we have identified $H: T^*M \rightarrow \mathbb{R}$ with its pullback through $\pi_{I \times T^*M}$. In coordinates, $\Theta_H(q,p) = pdq - Hdt.$ Additionally, we define the presymplectic form $\Omega_H = -d\Theta_H$. Furthermore, we identify curves on $T^*M$, of the form $\psi: I \rightarrow T^*M$, with curves on $I \times T^*M$ which cover the identity on $I$; in coordinates, this identification reads $\psi(t) = (t, q(t), p(t))$. The above action can then be expressed
$$ S[\psi] = \int_I \psi^*\Theta_H. $$
To enforce Type II boundary conditions $\pi_{T^*Q} (\psi(t_0)) = q_0 \in M$ and $\psi(t_1) = p_1 \in T_{q_1}^*M$ where $q_1 = \Phi_{t_1-t_0}(q_0)$, we define the space of admissible variations with respect to these boundary conditions as the space of vector fields $X$ on $T^*M$ (identified with vertical vector fields on $I \times T^*M \rightarrow T^*M$) such that $(T\pi_{T^*M}X)(q_0) = 0$ and $X(\psi_1) = 0$, where $\psi_1 = (q_1,p_1) \in T_{q_1}^*M$. Intuitively, the first condition states that an admissible variation does not vary the initial position $q(t_0) = q_0$, whereas the second condition states that an admissible variation does not vary the final momenta $\psi_1$. 

\begin{prop}\label{TypeIIVPProposition}
Fix an interval $I = (t_0,t_1) \subset \mathbb{R}$, $t_0\neq t_1$. Consider the above augmented Hamiltonian, where we assume that the time $t_1-t_0$ flow of the vector field $f$ exists. Let $q_0 \in M$ and let $p_1 \in T_{q_1}^*M$ where $q_1 = \Phi_{t_1-t_0}(q_0)$. Then, the augmented adjoint system with Type II boundary conditions
\begin{align*}
\dot{q} &= f(q), \\
\dot{p} &= - [Df(q)]^*p - dL(q), \\
q(t_0) &= q_0, \\
p(t_1) &= p_1,
\end{align*}
is intrinsically given by the variation principle: enforce the stationarity of the action
$$ S[\psi] = \int_I \psi^*\Theta_H $$
with respect to admissible variations.
\begin{proof}
Let $\varphi_\epsilon$ denote the time-$\epsilon$ flow of an admissible variation $X$. Then, the variation principle for the action with respect to admissible variations is given by
\begin{align*}
0 &= dS[\psi]\cdot X = \frac{d}{d\epsilon}\Big|_{\epsilon = 0} S[\varphi_\epsilon \circ \psi] = \int_I \psi^* \frac{d}{d\epsilon}\Big|_0 \varphi_{\epsilon}^*\Theta_H = \int_I \psi^* \mathcal{L}_X \Theta_H \\
&= - \int_I \psi^* (i_X \Omega_H) + \int_I \psi^* d(i_X\Theta_H) = -\int_I \psi^* (i_X\Omega_H) + \int_I d(\psi^*i_X\Theta_H).
\end{align*}
Observe that the boundary term $\int_I d(\psi^*i_X\Theta_H) = (\psi^*i_X\Theta_H)(t_1) - (\psi^*i_X\Theta_H)(t_0)$ vanishes by the fact that $X$ is an admissible variation since $(\psi^*i_X\Theta_H)(t) = \langle p(t), (T\pi_{T^*M}X)(q(t))\rangle.$ Hence, the stationarity condition is given by
$$ \int_I \psi^* (i_X\Omega_H) = 0. $$
By the fundamental lemma of the calculus of variations, we have $\psi^* (i_X\Omega_H) = 0$, whose coordinate expression is precisely the adjoint system.
\end{proof}
\end{prop}

\begin{remark}
In our definition of the space of admissible variations, we set the conditions that the variation at $q_0$ is purely vertical, $(T\pi_{T^*M}X)(q_0) = 0$, whereas at $q_1$, we enforced that the variation is zero, $X(q_1,p_1)=0$. In coordinates where 
$$ X = \delta q \frac{\partial}{\partial q} + \delta p \frac{\partial}{\partial p}, $$
the first condition reads $\delta q_0 = 0$ and the second condition reads $\delta q_1 = 0$, $\delta p_1 = 0$. It would thus seem that we are enforcing an overdetermined set of three boundary conditions $q(t_0) = q_0$, $q(t_1) = q_1$, $p(t_1) = p_1$. However, the resolution is that the variations $\delta q_0$ and $\delta q_1$ are not independent; fixing one to zero sets the other one to zero, by virtue of the fact that the adjoint system covers an ODE on $M$. Thus, with the chosen variational principle, we are only setting two independent boundary conditions, $q(t_0) = q_0, p(t_1) = p_1$. 

Furthermore, in the above proof, by looking at the coordinate expression of the boundary term,
$$ \psi^*i_X\Theta_H\Big|_{t_0}^{t_1} = \langle p(t_1), \delta q(t_1)\rangle - \langle p(t_0), \delta q(t_0)\rangle, $$
we see that we only used $\delta q_0 = 0$, $\delta q_1 = 0$. We did not need that $\delta p_1 = 0$ for the boundary terms to vanish. However, without setting $\delta p_1 = 0$, we only have the system 
\begin{align*}
\dot{q} &= f(q), \\
\dot{p} &= - [Df(q)]^*p - dL(q), \\
q(t_0) &= q_0.
\end{align*}
Hence, this system is underdetermined; any curve $p(t)$ in the fibers of $T^*M$ satisfying 
$$\dot{p} = - [Df(q)]^*p - dL(q)$$
would suffice. Thus, to uniquely fix the system, we must also supply a boundary condition of the form $p(t_1) = p_1$. Thus, even though the condition $\delta p_1 = 0$ is not strictly necessary in the variational principle to derive the equations of motion, it is necessary to fix the curve $p(t)$ in the fibers that define the adjoint system with Type II boundary conditions.
\end{remark}

Analogously the adjoint DAE system \eqref{Augmented Adjoint DAE System Coord 1}-\eqref{Augmented Adjoint DAE System Coord 4}, for index 1 DAEs, can be derived by an intrinsic Type II variational principle, by considering variations $V$ of the action
$$S[\psi] = \int_{t_0}^{t_1} pdq - (H(q,u,p,\lambda) + L(q,u)) dt$$
such that $T\Pi_{q}V|_{t_0} = 0$ and $T\Pi_{(q,p)}V|_{t_1}=0$ where $\Pi_q: (q,u,p,\lambda) \mapsto q$ and $\Pi_{(q,p)}: (q,u,p,\lambda)\mapsto (q,p)$ are the canonical bundle projections on $\overline{T^*M}_d \oplus \Phi^*$.

\nocite{*}

\bibliographystyle{plainnat}
\bibliography{adjointdirac.bib}

\begin{thebibliography}{43}
\providecommand{\natexlab}[1]{#1}
\providecommand{\url}[1]{\texttt{#1}}
\expandafter\ifx\csname urlstyle\endcsname\relax
  \providecommand{\doi}[1]{doi: #1}\else
  \providecommand{\doi}{doi: \begingroup \urlstyle{rm}\Url}\fi

\bibitem[Absil et~al.(2008)Absil, Mahony, and Sepulchre]{AbMaSe2008}
P.-A. Absil, R.~Mahony, and R.~Sepulchre.
\newblock \emph{Optimization Algorithms on Matrix Manifolds}.
\newblock Princeton University Press, 2008.

\bibitem[Aguiar et~al.(2021)Aguiar, Camponogara, and Foss]{AgCaFo2021}
M.A. Aguiar, E.~Camponogara, and B.~Foss.
\newblock An augmented {L}agrangian for optimal control of {DAE} systems:
  Algorithm and properties.
\newblock \emph{IEEE Transactions on Automatic Control}, 66\penalty0
  (1):\penalty0 261--266, 2021.

\bibitem[Barbero-Li{\~n}{\'a}n and {Mart{\'i}n de Diego}(2021)]{BaMa2021}
M.~Barbero-Li{\~n}{\'a}n and D.~{Mart{\'i}n de Diego}.
\newblock Retraction maps: a seed of geometric integrators.
\newblock \emph{arXiv}, 2106.00607, 2021.

\bibitem[Barbero-Li{\~n}{\'a}n and {Mart{\'i}n de Diego}(2022)]{BaMa2022}
M.~Barbero-Li{\~n}{\'a}n and D.~{Mart{\'i}n de Diego}.
\newblock Presymplectic integrators for optimal control problms via retraction
  maps.
\newblock \emph{arXiv}, 2203.00790, 2022.

\bibitem[Benning et~al.(2019)Benning, Celledoni, Ehrhardt, Owren, and
  Sch\"onlieb]{DeCeEhOwSc2019}
M.~Benning, E.~Celledoni, {M.~J.} Ehrhardt, B.~Owren, and C.-B. Sch\"onlieb.
\newblock Deep learning as optimal control problems: Models and numerical
  methods.
\newblock \emph{J. Comput. Dyn.}, 6\penalty0 (2):\penalty0 171--198, 2019.

\bibitem[Berglund(2007)]{Be2007}
N.~Berglund.
\newblock Perturbation theory of dynamical systems.
\newblock \emph{DEA}, 2007.
\newblock \doi{cel-00130703}.

\bibitem[Biegler(2010)]{Bi2010}
L.~T. Biegler.
\newblock \emph{Nonlinear Programming: Concepts, Algorithms, and Applications
  to Chemical Processes}.
\newblock Society for Industrial and Applied Mathematics, Philadelphia, PA,
  2010.

\bibitem[Bullo and Lewis(2014)]{BuLe2014}
F.~Bullo and A.~D. Lewis.
\newblock Supplementary chapters for \emph{{G}eometric {C}ontrol of
  {M}echanical {S}ystems}, aug 2014.
\newblock URL \url{http://motion.mee.ucsb.edu/book-gcms/}.

\bibitem[Burby and Klotz(2020)]{BuKl2020}
J.W. Burby and T.J. Klotz.
\newblock Slow manifold reduction for plasma science.
\newblock \emph{Communications in Nonlinear Science and Numerical Simulation},
  89:\penalty0 105289, 2020.

\bibitem[Cacuci(1981)]{Ca1981}
D.~G. Cacuci.
\newblock Sensitivity theory for nonlinear systems. {I}. {N}onlinear functional
  analysis approach.
\newblock \emph{J. Math. Phys.}, 22\penalty0 (12):\penalty0 2794--2802, 1981.

\bibitem[Cao et~al.(2003)Cao, Li, Petzold, and Serban]{CaLiPeSe2003}
Y.~Cao, S.~Li, L.~Petzold, and R.~Serban.
\newblock Adjoint sensitivity analysis for differential-algebraic equations:
  The adjoint {DAE} system and its numerical solution.
\newblock \emph{SIAM J. Sci. Comput.}, 24\penalty0 (3):\penalty0 1076--1089 (14
  pages), 2003.

\bibitem[Cari{\~n}ena et~al.(1987)Cari{\~n}ena, Ibort, Gomis, and
  Rom{\'a}n-Roy]{CaIbGoRo1987}
J.~F. Cari{\~n}ena, L.~A. Ibort, J.~Gomis, and N.~Rom{\'a}n-Roy.
\newblock Applications of the canonical-transformation theory for presymplectic
  systems.
\newblock \emph{Il Nuovo Cimento B (1971-1996)}, 98\penalty0 (2):\penalty0
  172--196, 1987.

\bibitem[Chen and Trenn(2021)]{ChTr2021}
Y.~Chen and S.~Trenn.
\newblock An approximation for nonlinear differential-algebraic equations via
  singular perturbation theory.
\newblock \emph{arXiv}, 2103.12146, 2021.

\bibitem[de~Le{\'o}n et~al.(2004)de~Le{\'o}n, Cort{\'e}s, de~Diego, and
  Mart{\'i}nez]{dLCdDM2004}
Manuel de~Le{\'o}n, Jorge Cort{\'e}s, {David Mart{\'i}n} de~Diego, and Sonia
  Mart{\'i}nez.
\newblock General symmetries in optimal control.
\newblock \emph{Rep. Math. Phys.}, 53\penalty0 (1):\penalty0 55--78, 2004.

\bibitem[Delgado-T\'{e}llez and Ibort(2003)]{DeIb2003}
M.~Delgado-T\'{e}llez and A.~Ibort.
\newblock A panorama of geometrical optimal control theory.
\newblock \emph{Extracta mathematicae}, 18\penalty0 (2):\penalty0 129--151,
  2003.

\bibitem[Echeverr\'{i}a-Enr\'{i}quez et~al.(2003)Echeverr\'{i}a-Enr\'{i}quez,
  Mar\'{i}n-Solano, Mu{\~n}oz-Lecanda, and Rom\'{a}n-Roy]{EcMaMuRo2003}
A.~Echeverr\'{i}a-Enr\'{i}quez, J.~Mar\'{i}n-Solano, M.C. Mu{\~n}oz-Lecanda,
  and N.~Rom\'{a}n-Roy.
\newblock Geometric reduction in optimal control theory with symmetries.
\newblock \emph{Reports on Mathematical Physics}, 52\penalty0 (1):\penalty0
  89--113, 2003.

\bibitem[Giles and Pierce(2000)]{GiPi2000}
M.~B. Giles and N.~A. Pierce.
\newblock An introduction to the adjoint approach to design.
\newblock \emph{Flow Turbul. Combust.}, 65\penalty0 (3):\penalty0 393--415,
  2000.

\bibitem[Gotay and Nester(1979)]{GoNe1979}
M.~J. Gotay and J.~M. Nester.
\newblock Presymplectic {L}agrangian systems. {I} : the constraint algorithm
  and the equivalence theorem.
\newblock \emph{Annales de l'I.H.P. Physique théorique}, 30\penalty0
  (2):\penalty0 129--142, 1979.

\bibitem[Gotay et~al.(1978)Gotay, Nester, and Hinds]{GoNeHi1978}
M.~J. Gotay, J.~M. Nester, and G.~Hinds.
\newblock Presymplectic manifolds and the {D}irac--{B}ergmann theory of
  constraints.
\newblock \emph{Journal of Mathematical Physics}, 19\penalty0 (11):\penalty0
  2388--2399, 1978.

\bibitem[Griewank(2003)]{Gr2003}
A.~Griewank.
\newblock A mathematical view of automatic differentiation.
\newblock In \emph{Acta Numer.}, volume~12, pages 321--398. Cambridge
  University Press, 2003.

\bibitem[Ibragimov(2006)]{Ib2006}
N.~H. Ibragimov.
\newblock Integrating factors, adjoint equations and {L}agrangians.
\newblock \emph{Journal of Mathematical Analysis and Applications},
  318\penalty0 (2):\penalty0 742--757, 2006.

\bibitem[Ibragimov(2007)]{Ib2007}
N.~H. Ibragimov.
\newblock A new conservation theorem.
\newblock \emph{J. Math. Anal. Appl.}, 333\penalty0 (1):\penalty0 311--328,
  2007.

\bibitem[Leok and Ohsawa(2011)]{LeOh2008}
M.~Leok and T.~Ohsawa.
\newblock Variational and geometric structures of discrete {D}irac mechanics.
\newblock \emph{Found. Comput. Math.}, 11\penalty0 (5):\penalty0 529--562,
  2011.

\bibitem[Leok and Zhang(2011)]{LeZh2011}
M.~Leok and J.~Zhang.
\newblock Discrete {H}amiltonian variational integrators.
\newblock \emph{IMA J. Numer. Anal.}, 31\penalty0 (4):\penalty0 1497--1532,
  2011.

\bibitem[Li and Petzold(2004)]{LiPe2004}
S.~Li and L.~Petzold.
\newblock Adjoint sensitivity analysis for time-dependent partial differential
  equations with adaptive mesh refinement.
\newblock \emph{Journal of Computational Physics}, 198\penalty0 (1):\penalty0
  310--325, 2004.

\bibitem[Li and Petzold(2003)]{LiPe2003}
S.~Li and L.~R. Petzold.
\newblock Solution adapted mesh refinement and sensitivity analysis for
  parabolic partial differential equation systems.
\newblock In L.~T. Biegler, M.~Heinkenschloss, O.~Ghattas, and B.~van
  Bloemen~Waanders, editors, \emph{Large-Scale PDE-Constrained Optimization},
  pages 117--132. Springer, Berlin, 2003.

\bibitem[Mattsson and S{\"o}derlind(1993)]{MaSo1993}
S.~E. Mattsson and G.~S{\"o}derlind.
\newblock Index reduction in differential-algebraic equations using dummy
  derivatives.
\newblock \emph{SIAM Journal on Scientific Computing}, 14\penalty0
  (3):\penalty0 677--692, 1993.

\bibitem[Nguyen et~al.(2016)Nguyen, Georges, and Besan\c{c}on]{NgGeBe2016}
V.~T. Nguyen, D.~Georges, and G.~Besan\c{c}on.
\newblock State and parameter estimation in 1-{D} hyperbolic {PDEs} based on an
  adjoint method.
\newblock \emph{Automatica}, 67\penalty0 (C):\penalty0 185--191, May 2016.
\newblock ISSN 0005-1098.

\bibitem[Pierce and Giles(2000)]{PiGi2000}
N.~A. Pierce and M.~B. Giles.
\newblock Adjoint recovery of superconvergent functionals from {PDE}
  approximations.
\newblock \emph{SIAM Rev.}, 42\penalty0 (2):\penalty0 247--264, 2000.

\bibitem[Reid et~al.(2001)Reid, Lin, and Wittkopf]{ReLiWi2001}
G.~J. Reid, P.~Lin, and A.~D. Wittkopf.
\newblock Differential elimination–completion algorithms for {DAE} and
  {PDAE}.
\newblock \emph{Studies in Applied Mathematics}, 106\penalty0 (1):\penalty0
  1--45, 2001.

\bibitem[Roche(1989)]{Ro1989}
M.~Roche.
\newblock Implicit {R}unge--{K}utta methods for differential algebraic
  equations.
\newblock \emph{SIAM Journal on Numerical Analysis}, 26\penalty0 (4):\penalty0
  963--975, 1989.

\bibitem[Ross(2005)]{Ro2005}
I.~M. Ross.
\newblock A roadmap for optimal control: The right way to commute.
\newblock \emph{Ann. NY Acad. Sci.}, 1065\penalty0 (1):\penalty0 210--231,
  2005.

\bibitem[Ross and Fahroo(2001)]{RoFa2001}
M.~Ross and F.~Fahroo.
\newblock A pseudospectral transformation of the convectors of optimal control
  systems.
\newblock \emph{IFAC Proc. Ser.}, 34\penalty0 (13):\penalty0 543--548, 2001.

\bibitem[Sanz-Serna(2016)]{Sa2016}
J.~M. Sanz-Serna.
\newblock Symplectic {R}unge--{K}utta schemes for adjoint equations, automatic
  differentiation, optimal control, and more.
\newblock \emph{SIAM Review}, 58\penalty0 (1):\penalty0 3--33, 2016.

\bibitem[Schmitt and Leok(2017)]{ScLe2017}
J.~M. Schmitt and M.~Leok.
\newblock {Properties of {H}amiltonian variational integrators}.
\newblock \emph{IMA Journal of Numerical Analysis}, 38\penalty0 (1):\penalty0
  377--398, 2017.

\bibitem[Seiler(1998)]{Se1998}
W.~M. Seiler.
\newblock Numerical analysis of constrained {H}amiltonian systems and the
  formal theory of differential equations.
\newblock \emph{Math. Comput. Simul.}, 45\penalty0 (5–6):\penalty0 561–576,
  1998.

\bibitem[Sirkes and Tziperman(1997)]{SoTz1997}
Z.~Sirkes and E.~Tziperman.
\newblock Finite difference of adjoint or adjoint of finite difference?
\newblock \emph{Mon. Weather Rev.}, 125\penalty0 (12):\penalty0 3373--3378,
  1997.

\bibitem[Tran and Leok(2022)]{TrLe2022}
B.~Tran and M.~Leok.
\newblock Multisymplectic {H}amiltonian variational integrators.
\newblock \emph{International Journal of Computer Mathematics (Special Issue on
  Geometric Numerical Integration, Twenty-Five Years Later)}, 99\penalty0
  (1):\penalty0 113--157, 2022.

\bibitem[Vankerschaver et~al.(2012)Vankerschaver, Yoshimura, and
  Leok]{VaYoLeMa2012}
J.~Vankerschaver, H.~Yoshimura, and M.~Leok.
\newblock The {H}amilton--{P}ontryagin principle and multi-{D}irac structures
  for classical field theories.
\newblock \emph{J. Math. Phys.}, 53\penalty0 (7):\penalty0 072903 (25 pages),
  2012.

\bibitem[Wang et~al.(2012)Wang, Duraisamy, Alonso, and Iaccarino]{WaDuAlIa2012}
Q.~Wang, K.~Duraisamy, J.~J. Alonso, and G.~Iaccarino.
\newblock Risk assessment of scramjet unstart using adjoint-based sampling.
\newblock \emph{AIAA J.}, 50\penalty0 (3):\penalty0 581--592, 2012.

\bibitem[Yano and Ishihara(1973)]{YaIs1973}
K.~Yano and S.~Ishihara.
\newblock \emph{Tangent and cotangent bundles: differential geometry}.
\newblock Pure Appl. Math., No. 16. Marcel Dekker, Inc., New York, 1973.

\bibitem[Yoshimura and Marsden(2006{\natexlab{a}})]{YoMa2006a}
H.~Yoshimura and {J.~E.} Marsden.
\newblock Dirac structures in {L}agrangian mechanics {P}art {I}: {I}mplicit
  {L}agrangian systems.
\newblock \emph{J. Geom. Phys.}, 57\penalty0 (1):\penalty0 133--156,
  2006{\natexlab{a}}.

\bibitem[Yoshimura and Marsden(2006{\natexlab{b}})]{YoMa2006b}
H.~Yoshimura and {J.~E.} Marsden.
\newblock Dirac structures in {L}agrangian mechanics {P}art {II}: {V}ariational
  structures.
\newblock \emph{J. Geom. Phys.}, 57\penalty0 (1):\penalty0 209--250,
  2006{\natexlab{b}}.

\end{thebibliography}

\end{document}